\documentclass[12pt]{article}
\RequirePackage[OT1]{fontenc}
\RequirePackage{amsthm,amsmath}
\RequirePackage[authoryear]{natbib}

\usepackage{color, dsfont}
\usepackage{url}
\usepackage{amsmath,amsfonts,amsthm,amssymb,amsbsy}

\newtheorem{corollary}{Corollary}[section]

\newtheorem{example}{Example}[section]
\newtheorem{lemma}{Lemma}[section]

\newtheorem{remark}{Remark}[section]
\newtheorem{theorem}{Theorem}[section]

\begin{document}

	
\title{\large\bf D-optimal Designs for Multinomial Logistic Models}
\author{Xianwei Bu$^{1}$, Dibyen Majumdar$^{2}$ and Jie Yang$^{2}$\\
	$^1$AbbVie Inc. and $^2$University of Illinois at Chicago}
	
\maketitle

\begin{abstract}
We consider optimal designs for general multinomial logistic models, which cover baseline-category, cumulative, adjacent-categories, and continuation-ratio logit models, with proportional odds, non-proportional odds, or partial proportional odds assumption. We derive the corresponding Fisher information matrices in three different forms to facilitate their calculations, determine the conditions for their positive definiteness, and search for optimal designs. We conclude that, unlike the designs for binary responses, a feasible design for a multinomial logistic model may contain less experimental settings than parameters, which is of practical significance. We also conclude that even for a minimally supported design, a uniform allocation, which is typically used in practice, is not optimal in general for a multinomial logistic model. We develop efficient algorithms for searching D-optimal designs. Using examples based on real experiments, we show that the efficiency of an experiment can be significantly improved if our designs are adopted.
\end{abstract}


\section{Introduction}\label{sec:intro}

Experiments with categorical responses arise naturally in a rich variety of scientific disciplines \citep{christensen2013}. While there is a sizable and growing literature for optimal designs with binary response, the literature for experiments with more than two categories is limited.
On the other hand, there are many applications where three or more categories arise naturally. For instance, the response of an odor removal study \citep{ytm2016} has three levels, {\tt serious odor}, {\tt medium odor}, {\tt almost no odor} for manufactured bio-plastics; and the outcome of trauma clinical trial \citep{chuang1997} has five categories known as the Glasgow Outcome Scale \citep{jennett1975}, {\tt death}, {\tt vegetative} {\tt state}, {\tt major disability}, {\tt minor disability}, and {\tt good} {\tt recovery}. In some experiments responses are hierarchical, such as the outcome of emergence of house flies \citep{atkinson1999} with categories, {\tt died before opening of the pupae}, {\tt died before complete emergence}, and {\tt completely emerged}.
Other examples include a wine bitterness study \citep{randall1989}, a polysilicon deposition study \citep{wu2008}, and a toxicity study \citep{agresti2013}. Our aim is to obtain foundational results for the identification of optimal and efficient designs for experiments with three or more responses and explore the properties of these designs.

When the response is binary, generalized linear models have been used widely (McCullagh and Nelder, 1989; Dobson and Barnett, 2008) for analyzing the experimental data. 
For optimal designs of experiments with generalized linear models for {\it univariate} responses, there is a growing body of literature \citep{khuri2006, atkinson2007, stufken2012}. In this case, the minimum number of distinct experimental settings required by a nondegenerate Fisher information matrix is equal to the number of parameters \citep{fedorov1972, ym2015}. Moreover, for the widely studied {\it D-optimal approximate designs} the experimental units are uniformly allocated when a minimally supported design, that is, a design with the minimum number of distinct experimental settings, is adopted \citep{ym2015,ymm2016}.

For responses with three or more categories, these results do not directly apply. The models used in the literature are special cases of the {\it multivariate} generalized linear model \citep{pmcc1980}. According to the relationship among categories, categorical responses can be of three types: nominal, ordinal, and hierarchical \citep{atkinson1999}. In practice, a multinomial distribution is typically used to model the responses for a specified experimental setting. In the statistical literature, four kinds of logit models have been commonly used to link the categorical probabilities to experimental settings or values of covariates, the {\it baseline-category logit model} for nominal responses \citep{agresti2013, atkinson1999}, the {\it cumulative logit model} for ordinal responses \citep{pmcc1980, christensen2013}, the {\it adjacent-categories logit model} for ordinal responses \citep{agresti2005, agresti2013}, and the {\it continuation-ratio logit model} for hierarchical responses \citep{agresti2013, atkinson1999}.
Also, three different assumptions on parameter structures have been proposed for the four logit models, {\it proportional odds} ({\it po}) assuming the same parameters for different categories \citep{pmcc1980}, {\it non-proportional odds} ({\it npo}) allowing parameters to change across categories \citep{agresti2013}, and {\it partial proportional odds} ({\it ppo}) containing both {\it po} and {\it npo} components \citep{peterson1990} as special cases. The four logit models and three odds assumptions generate $12$ different models for multinomial responses including, for example, a cumulative logit model with proportional odds. All $12$ models can be fitted using {\tt SAS} \citep{stokes2012} or {\tt R} \citep{yee2015}, and AIC or BIC criterion may be used for model selection.

Despite this rich complexity of models, the relevant results in the design literature for multinomial responses are limited to special classes. \cite{atkinson1999} constructed a general framework of optimal designs for multinomial logistic models with non-proportional odds. \cite{perevozskaya2003} discussed a special class of cumulative logit models with proportional odds. \cite{ytm2016} obtained results for the cumulative link model, which is an extension of the cumulative logit models with proportional odds. Although the cumulative logit model with proportional odds is the most popular model in practice for ordinal responses, \cite{agresti2010} found strong evidence against the assumption of proportional odds for the trauma clinical trial data. For the emergence of house flies data \citep{atkinson1999}, a continuation-ratio logit model with non-proportional odds fits better. It is also known that for certain applications, partial proportional odds may have a better performance than {\it po} or {\it npo} models \citep{lall2002}.

The goal of this work is a comprehensive approach to the study of optimal designs for multinomial responses. We work in a general framework, which covers all of the 12 models. The traditional starting point for obtaining an optimal design, given a model, is the Fisher information matrix for an arbitrary design. This information matrix is then ``maximized'' using a criterion function; for instance, the criterion we use, D-optimality, maximizes the determinant of the Fisher information matrix. Our first step is to derive explicit representations of the Fisher information matrix, thereby providing the foundation for the quest of all optimal designs. Our next step is to derive conditions for the positive definiteness of the information matrix, which is necessary for the non-degeneracy of the design and essential for formulating the criteria functions. This also allows us to determine the minimum number of experimental conditions needed for positive definiteness, which is necessary to understand the cost-efficiency of the study. Then we proceed to study D-optimal designs. The Fisher information matrix for the multinomial logistic models, as in all models except linear models, depends on the (unknown) model parameters. There are two approaches in the optimal design literature to deal with this. The first approach, local optimality, which we pursue, uses assumed values of the parameters. Difficulties arise when no prior information, such as prior experimental result, is available; we will address this issue in Section~\ref{section:exactdesign}. The second approach, that requires specification of a prior distribution of the parameters is Bayesian optimality \citep{chaloner1995}. Using the D-criterion, these approaches lead to locally D-optimal and Bayesian D-optimal designs. Instead of Bayesian optimality, however, we study a variant, EW optimality \citep{atkinson2007, ymm2016, ytm2016}, which maximizes the determinant of the expected Fisher information matrix under a prior on the parameters; this is a good surrogate of Bayesian D-optimal designs, with substantially reduced computational cost.

Our work shows that the optimal designs for multinomial responses with three or more categories are remarkably different from the ones for binary responses in two major aspects: (i) the required minimum number of experimental settings is less than the number of parameters; (ii) even among minimally supported designs, unlike D-optimal designs for binary responses, uniform allocation is not D-optimal, except for {\it regular npo} models (defined in Section~\ref{section:minimallysupported}). Indeed, we find that uniform designs can be quite inefficient. Theoretically this work reveals significant new features of optimal designs for general multinomial models. It also provides a way to accurately compute efficiency of designs of experiments based on multinomial models, and shows that widely used designs, like uniform designs, may not be efficient in these models. More generally, this work indicates that as we go from univariate to multivariate responses, some widely-applicable tenets of optimal design theory, like uniform allocation on a minimal set of experimental conditions, may no longer hold. 

Similar conclusions for minimally supported designs were observed by \cite{ytm2016} for cumulative link models and proportional odds. Our results confirm and extend these to the general setup.

The rest of this paper is organized as follows: In Section~\ref{sec:fisher}, we formulate the general multinomial logistic model and its Fisher information matrix; in Section~\ref{sec:positive}, we derive a necessary and sufficient condition for the Fisher information matrix to be positive definite, which determines the minimal number of required experimental settings; in Section~\ref{sec:determinant}, we formulate the determinant of the Fisher information matrix as a homogeneous polynomial of allocations of experimental units and simplify its structure for D-optimality; in Section~\ref{sec:specifiedX}, we develop numerical algorithms for searching D-optimal designs, as well as some analytical results for minimally supported designs; we conclude with discussions in Section~\ref{sec:discussion}. Although we focus on D-optimality, our basic results on Fisher information matrix are useful for other criteria as well, such as A-optimality, E-optimality, etc. We also study two examples extensively, the trauma clinical trial \citep{chuang1997} and the emergence of house flies \citep{atkinson1999}. These are used to illustrate the results and explore the properties of the designs, such as efficiency and robustness.

One important point to note is that throughout this paper the major results are expressed in terms of {\it ppo} models only, as this includes both {\it po} and {\it npo} models as special cases. The special cases are of considerable interest, however, especially in many applications, and all of the specific formulae for {\it po} and {\it npo} models are provided in the Supplementary Materials, so they are readily available to the reader. The Supplementary Materials also contain results that are of further interest.

\section{Multinomial Logistic Model and Its Fisher Information Matrix}\label{sec:fisher}

This section is dedicated to the description of the model and the Fisher information matrix, the starting points of the study of optimal designs. We consider an experiment with $d\geq 1$ factors and $m\geq 2$ distinct experimental settings ${\mathbf x}_i = (x_{i1}, \ldots, x_{id})^T$, $i=1, \ldots, m$, that could be used in the experiment. For the $i$th experimental setting, $n_i\geq 0$ categorical responses are collected i.i.d.~from a discrete distribution with $J\geq 2$ categories, with $n_i=0$ indicating no experimental unit assigned to this experimental setting. When $n_i > 0$, the responses associated with the $i$th experimental setting are summarized into a multinomial response
${\mathbf Y}_i=(Y_{i1},\cdots,Y_{iJ})^T \sim {\rm Multinomial}(n_i; \pi_{i1},\cdots,\pi_{iJ})$,
where $\pi_{ij}$ is the probability that the response falls into the $j$th category at the $i$th experimental setting. Note that $\pi_{ij}$’s are functions of the experimental settings and the model parameters (to be introduced next); however, for simplicity we do not show the dependence in the notation of $\pi_{ij}$. Throughout this paper, we assume $\pi_{ij}>0$ for all $i=1, \ldots, m$ and $j=1, \ldots, J$, which is necessary for a multinomial logistic model (in Section~\ref{sec:designspace} we examine the implications of this on the {\it design space}).

The general structure of models that are considered is a linear regression of log odds on two components: one specific to the category and the other common to all categories. As mentioned in the Introduction, we will work with the most general model – the partial proportional odds ({\it ppo}) model. For this, 
we write the four logit models (baseline-category, cumulative, adjacent-categories, and continuation-ratio) in terms of {\it ppo} structure as follows:
\begin{eqnarray*}
	\log\left(\frac{\pi_{ij}}{\pi_{iJ}}\right) &=& {\mathbf h}_j^T({\mathbf x}_i)\boldsymbol\beta_j+{\mathbf h}_c^T({\mathbf x}_i)\boldsymbol\zeta~, \>\mbox{baseline-category}\\
	\log\left(\frac{\pi_{i1}+\cdots + \pi_{ij}}{\pi_{i,j+1}+\cdots + \pi_{iJ}}\right) &=& {\mathbf h}_j^T({\mathbf x}_i)\boldsymbol\beta_j+{\mathbf h}_c^T({\mathbf x}_i)\boldsymbol\zeta~, \>\mbox{cumulative}\\
	\log\left(\frac{\pi_{ij}}{\pi_{i,j+1}}\right) &=& {\mathbf h}_j^T({\mathbf x}_i)\boldsymbol\beta_j+{\mathbf h}_c^T({\mathbf x}_i)\boldsymbol\zeta~, \>\mbox{adjacent-categories}\\
	\log\left(\frac{\pi_{ij}}{\pi_{i,j+1} + \cdots + \pi_{iJ}}\right) &=& {\mathbf h}_j^T({\mathbf x}_i)\boldsymbol\beta_j+{\mathbf h}_c^T({\mathbf x}_i)\boldsymbol\zeta~, \>\mbox{continuation-ratio}
\end{eqnarray*}
where $i=1, \ldots, m$, $j=1, \ldots, J-1$, ${\mathbf h}_j^T(\cdot) = (h_{j1}(\cdot), \ldots, h_{jp_j}(\cdot))$ are known functions to determine the $p_j$ predictors associated with the $p_j$ unknown parameters $\boldsymbol\beta_j = (\beta_{j1}, \ldots, \beta_{jp_j})^T$ for the $j$th response category, and ${\mathbf h}_c^T(\cdot) = (h_{1}(\cdot), \ldots, h_{p_c}(\cdot))$ are known functions to determine the $p_c$ predictors associated with the $p_c$ unknown parameters $\boldsymbol\zeta = (\zeta_{1}, \ldots, \zeta_{p_c})^T$ that are common for all categories.  As special cases of {\it ppo}, ${\mathbf h}_j^T({\mathbf x}_i) \equiv 1$ leads to {\it po} models, and ${\mathbf h}_c^T({\mathbf x}_i) \equiv 0$ leads to {\it npo} models. The corresponding expressions for {\it po} and {\it npo} models are provided in the Supplementary Materials (Sections~\ref{subsec:pomodels} and \ref{subsec:npomodels}).

Following \cite{pmcc1995} and \cite{atkinson1999}, we rewrite these four logit models into a unified form
\begin{equation}\label{logitunifiedmodel}
{\mathbf C}^T\log({\mathbf L}{\boldsymbol\pi}_i)={\boldsymbol\eta}_i={\mathbf X}_i{\boldsymbol\theta}, \qquad     i=1,\cdots,m
\end{equation}
where ${\boldsymbol\pi}_i = (\pi_{i1}, \ldots, \pi_{iJ})^T$ satisfying $\sum_{j=1}^J \pi_{ij}=1$, ${\boldsymbol\eta}_i = (\eta_{i1}, \ldots, \eta_{iJ})^T$,
\[
{\mathbf C}^T=
\left(\begin{array}{ccc}
{\mathbf I}_{J-1} & -{\mathbf I}_{J-1} & {\mathbf 0}_{J-1}\\
{\mathbf 0}_{J-1}^T & {\mathbf 0}_{J-1}^T & 1
\end{array}\right)
\]
is a $J\times(2J-1)$ constant matrix, where ${\mathbf I}_k$ is the identity matrix of order $k$ and ${\mathbf 0}_k$ is a vector of $k$ zeros, and
${\mathbf L}$ is a $(2J-1)\times J$ constant matrix taking different forms across the four logit models (see Appendix~\ref{subsec:L} for details). The model matrix is
\begin{equation}\label{eqn:Xi_ppo}
{\mathbf X}_i= \begin{pmatrix}
{\mathbf h}_1^T({\mathbf x}_i) &   & & {\mathbf h}_c^T({\mathbf x}_i)\\
&   \ddots &  & \vdots \\
&   & {\mathbf h}_{J-1}^T({\mathbf x}_i) & {\mathbf h}_c^T({\mathbf x}_i)\\
\boldsymbol0^T_{p_1} & \cdots & \boldsymbol0^T_{p_{J-1}} & \boldsymbol0^T_{p_c}\\
\end{pmatrix}_{J \times p}
\end{equation}
and the parameter vector $\boldsymbol\theta=(\boldsymbol\beta_{1}^T,\cdots,\boldsymbol\beta_{J-1}^T,\boldsymbol\zeta^T)^T$ consists of $p=p_1+\cdots+p_{J-1}+p_c$ unknown parameters in total. 
Note that $\pi_{i1}+\cdots+\pi_{iJ}=1$ implies that $\eta_{iJ}=0$ and thus the last row of ${\mathbf X}_i$ is all 0's. We keep $\eta_{iJ}$ and the last row of ${\mathbf X}_i$ in (\ref{logitunifiedmodel}) for convenience following \cite{pmcc1995}. Model~(\ref{logitunifiedmodel}) covers all four logit models and all three odds structures ({\it po}, {\it npo}, and {\it ppo}).

\begin{example}\label{example:yee2015} {\rm \citeauthor{yee2015} (2015) considered a {\it ppo} model with $d=4$ factors, $J=3$ response categories, and parameter vector $\boldsymbol\theta = (\beta_{11}, \beta_{12},$ $\beta_{13},$ $\beta_{21},$ $\beta_{22},$ $\beta_{23}, \zeta_{1}, \zeta_{2})^T$.
		At the $i$th experimental setting ${\mathbf x}_i=(x_{i1}, x_{i2}, x_{i3}, x_{i4})^T$,
		\[
		{\mathbf X}_i =
		\begin{pmatrix}
		1 & x_{i1} & x_{i2} &  0 & 0 & 0 & x_{i3} & x_{i4} \\
		0 & 0 & 0 & 1 & x_{i1} & x_{i2}  & x_{i3} & x_{i4}\\
		0 & 0 & 0 & 0 & 0 & 0 & 0  & 0
		\end{pmatrix}.
		\]
		We will revisit this example in Section~\ref{sec:determinant}.
		\hfill{$\Box$}
}\end{example}

Using matrix differentiation formulae (see, for example, \citeauthor{seber2008} (2008, Chapter~17)), we obtain the Fisher information matrix for model~\eqref{logitunifiedmodel} as follows:

\begin{theorem}\label{thm:fisher}
	Consider the multinomial logistic model \eqref{logitunifiedmodel} with independent observations. The Fisher information matrix
	\begin{equation}\label{eqn:F}
	{\mathbf F}=\sum_{i=1}^{m}n_i{\mathbf F}_i
	\end{equation}
	where
	\begin{equation}\label{eqn:Fi}
	{\mathbf F}_i=(\frac{\partial \boldsymbol\pi_i}{\partial \boldsymbol\theta^T})^T {\rm diag}(\boldsymbol\pi_i)^{-1}\frac{\partial \boldsymbol\pi_i}{\partial \boldsymbol\theta^T}
	\end{equation}
	with $\partial \boldsymbol\pi_i/\partial \boldsymbol\theta^T=({\mathbf C}^T{\mathbf D}_i^{-1}{\mathbf L})^{-1}{\mathbf X}_i$ and ${\mathbf D}_i={\rm diag}({\mathbf L} \boldsymbol\pi_i)$.
\end{theorem}
Theorem~\ref{thm:fisher} is a special case of \cite{pmcc1995} who built a more general framework for multiple categorical responses. We provide independent proofs in the Supplementary Materials (Section~\ref{sec:proofs}), as well as a result, Lemma~\ref{lemma:pii}, for use later on. Our results apply to more general models than \cite{atkinson1999}.

\begin{remark}\label{remark:thm2.1}{\rm
The Fisher information matrix ${\mathbf F}$ plays a key role in optimal design theory. For example, a D-optimal design maximizes the determinant of ${\mathbf F}$, an A-optimal design minimizes the trace of ${\mathbf F}^{-1}$, and an E-optimal design maximizes the minimum eigenvalue of ${\mathbf F}$. 
		Given experimental settings ${\mathbf x}_1, \ldots, {\mathbf x}_m$ and the parameter vector $\boldsymbol\theta$, one can calculate ${\mathbf F}_1, \ldots, {\mathbf F}_m$ using \eqref{eqn:Fi}. Then Theorem~\ref{thm:fisher} provides a convenient way for calculating ${\mathbf F}$ as a function of the allocation $(n_1, \ldots, n_m)$. The nonsingularity of the key matrix ${\mathbf C}^T{\mathbf D}_i^{-1}{\mathbf L}$ was guaranteed by \citeauthor{pmcc1995} (1995, Theorem~1). To facilitate calculations, we derive explicit forms of $({\mathbf C}^T{\mathbf D}_i^{-1}{\mathbf L})^{-1}$ for all the four logit models in the Supplementary Materials (Section~\ref{subsec:cij}).
}\end{remark}

\section{Positive Definiteness of the Fisher Information Matrix}\label{sec:positive}

Positive definiteness of the information matrix is essentially a necessary condition for the existence of unbiased estimators of parameters with finite variance \citep{stoica2001}. In this section our focus is the determination of conditions for the positive definiteness. We also determine the minimal number of experimental settings required for positive definiteness, which is a basic question, since it deals with the cost of experimentation. In order to do that we first derive a more amenable representation of ${\mathbf F}$.   

To study the structure of ${\mathbf F}$, we first denote $({\mathbf C}^T{\mathbf D}_i^{-1}{\mathbf L})^{-1} = ({\mathbf c}_{i1},$ $\ldots,$ ${\mathbf c}_{iJ})$, with ${\mathbf c}_{ij}$ denoting the $J\times 1$ column $j$. For the $i$th experimental setting, we define $u_{st}({\boldsymbol\pi}_i)={\mathbf c}_{is}^T {\rm diag}(\boldsymbol\pi_i)^{-1} {\mathbf c}_{it}$, for $s,t=1, \ldots, J-1$ and then a $J\times J$ matrix
\begin{equation}\label{eqn:Ui}
{\mathbf U}_i=\left(
\begin{array}{cccc}
u_{11}({\boldsymbol\pi}_i) & \cdots & u_{1,J-1}({\boldsymbol\pi}_i) & 0\\
\vdots & \ddots & \vdots & \vdots\\
u_{J-1,1}({\boldsymbol\pi}_i) & \cdots & u_{J-1,J-1}({\boldsymbol\pi}_i) & 0\\
0 & \cdots & 0 & 1
\end{array}\right)
\end{equation}
Then we can rewrite the Fisher information at the $i$th experimental setting ${\mathbf F}_i$ into a simpler form as a corollary of Theorem~\ref{thm:fisher}:

\begin{corollary}\label{col:Fi}
	Under the setup of Theorem~\ref{thm:fisher},
	${\mathbf F}_i = {\mathbf X}_i^T {\mathbf U}_i {\mathbf X}_i$~.
\end{corollary}

In order to reformulate the Fisher information matrix ${\mathbf F}$ into a form that facilitates the discussion of the positive definiteness of ${\mathbf F}$, we further define an $m(J-1)\times m(J-1)$ matrix ${\mathbf U} = ({\mathbf U}_{st})_{s,t=1,\ldots, J-1}$ with ${\mathbf U}_{st}= {\rm diag}\{n_1 u_{st}({\boldsymbol\pi}_1), \ldots, n_m u_{st}({\boldsymbol\pi}_m)\}$, and a $p\times m(J-1)$ matrix
\begin{equation}\label{eq:H_matrix}
{\mathbf H}=
\left(\begin{array}{ccc}
{\mathbf H}_{1} & &\\  & \ddots & \\& &{\mathbf H}_{J-1} \\{\mathbf H}_{c} & \cdots & {\mathbf H}_{c}
\end{array}\right)
\end{equation}
where ${\mathbf H}_{j}= ({\mathbf h}_j({\mathbf x}_1), \cdots, {\mathbf h}_j({\mathbf x}_m))$ and ${\mathbf H}_{c} = ({\mathbf h}_c({\mathbf x}_1), \cdots, {\mathbf h}_c({\mathbf x}_m))$.

With the aid of Corollary~\ref{col:Fi}, we obtain the theorem below, whose proof is relegated to the Supplementary Materials (Section~\ref{sec:proofs}).

\begin{theorem}\label{thm:fisher_all}
	Consider the multinomial logistic model (\ref{logitunifiedmodel}) with independent observations. The Fisher information matrix
	${\mathbf F}= {\mathbf H}{\mathbf U}{\mathbf H}^T$~.
\end{theorem}

It can be verified that the matrix ${\mathbf U}$ is positive definite if $n_i > 0$ for all $i=1, \ldots,m$ (see Section~\ref{subsec:positiveU}). For general cases, if $n_i=0$ for some $i$, one can denote $k = \#\{i:n_i>0\}$ and ${\mathbf U}_{st}^* = {\rm diag}\{n_i u_{st}({\boldsymbol\pi}_i) : n_i>0\}$, then the reduced $k(J-1)\times k(J-1)$ matrix ${\mathbf U}^* = ({\mathbf U}^*_{st})_{s,t=1, \ldots, J-1}$ is still positive definite. In this case, one can remove all columns of ${\mathbf H}$ associated with $n_i=0$ and denote the leftover as ${\mathbf H}^*$,  which is a $p\times k(J-1)$ matrix.

\begin{theorem}\label{thm:Fpd}
	The Fisher information matrix ${\mathbf F}$ is positive definite if and only if ${\mathbf H}^*$ is of full row rank. Furthermore, if $n_i>0$ for all $i=1, \ldots, m$, then ${\mathbf F}$ is positive definite if and only if ${\mathbf H}$ is of full row rank.
\end{theorem}

Detailed discussion on the row rank of the matrix ${\mathbf H}$ is relegated to the Supplementary Materials (Section~\ref{subsec:rowrankH}).
As a direct conclusion, we obtain the main result for the positive definiteness of the Fisher information matrix ${\mathbf F}$.

\begin{theorem}\label{thm:minimalm}
	Consider the multinomial logistic model (\ref{logitunifiedmodel}) with $m$ distinct experimental settings ${\mathbf x}_i$ and $n_i>0$ experimental units, $i=1, \ldots, m$. The Fisher information matrix ${\mathbf F}$ is positive definite if and only if
	\begin{itemize}
		\item[(1)] $m\geq \max\{p_1, \ldots, p_{J-1},p_c+p_H\}$, where $p_H = dim (\cap_{j=1}^{J-1} {\cal M}({\mathbf H}_j^T ))$ and ${\cal M}({\mathbf H}_j^T )$ stands for the column space of ${\mathbf H}_j^T$; and
		\item[(2)] ${\mathbf x}_i$'s keep ${\mathbf H}_j$ of full row rank $p_j$, $j=1, \ldots, J-1$; and keep ${\mathbf H}_c$ of full row rank $p_c$; as well as ${\cal M}({\mathbf H}_c^T)\cap (\cap_{j=1}^{J-1} {\cal M}({\mathbf H}_j^T))=\{0\}$.
	\end{itemize}
	As a special case, if ${\mathbf H}_1 = \cdots ={\mathbf H}_{J-1}$, then ${\mathbf F}$ is positive definite if and only if $m\geq p_c + p_1$ and the extended matrix $({\mathbf H}_1^T, {\mathbf H}_c^T)$ is of full rank $p_c + p_1$.
\end{theorem}

Since the determination of ${\cal M}({\mathbf H}_c^T)\cap (\cap_{j=1}^{J-1} {\cal M}({\mathbf H}_j^T) )$ in Theorem~\ref{thm:minimalm} is not straightforward, we provide a formula for its calculation in Appendix~\ref{subsec:dim}.

Theorem~\ref{thm:minimalm} implies that the number $m$ of distinct experimental settings could be as low as $\max\{p_1, \ldots, p_{J-1},p_c+p_H\}$, denoted as $k_{\min}$, which is strictly less than the number of parameters $p = p_1 + \cdots + p_{J-1} + p_c$ if $J\geq 3$. If the model is constructed with the same set of predictors for different categories, for example, a main-effects model, then ${\mathbf H}_1 = \cdots = {\mathbf H}_{J-1}$ and $k_{\min}$ is just $p_c+p_1$. This confirms the finding by \citeauthor{ytm2016} (2017, Theorem~3) on cumulative link models and extends it to the general case covering all the 12 multinomial logistic models. 
We examine examples in Section~\ref{sec:specifiedX}, and revisit the issue of minimal support in Section~\ref{section:minimallysupported}. For the applications that we have examined, $k_{\min}$ experimental settings ensure the positive definiteness of ${\mathbf F}$. We provide more examples in the Supplementary Materials (Section~\ref{subsec:moreexamples}).

\section{Determinant of the Fisher Information Matrix}\label{sec:determinant}

Until now, we have derived basic results that are needed for the identification of optimal designs irrespective of the specific criterion. We also derived a result on the minimal number of experimental settings. Now we turn our focus to D-optimal designs. The D-criterion for optimal designs maximizes the determinant of the Fisher information matrix ${\mathbf F}$. In order to study D-optimal designs and their properties, we start with a characterization of the structure of $|{\mathbf F}|$. Recall that $n$ is the total number of observations with $n_i$ of them assigned to the $i$th experimental setting ${\mathbf x}_i$. We further define an $mJ\times mJ$ diagonal matrix ${\mathbf W} = {\rm diag}\{w_1{\rm diag}({\boldsymbol\pi}_1)^{-1}, \ldots,$ $w_m {\rm diag}({\boldsymbol\pi}_m)^{-1}\}$ with proportions $w_i = n_i/n$, and an $mJ\times p$ matrix
\begin{equation}
{\mathbf G} = \left(\begin{array}{cccc}
{\mathbf c}_{11}{\mathbf h}_1^T({\mathbf x}_1) & \cdots & {\mathbf c}_{1,J-1} {\mathbf h}_{J-1}^T ({\mathbf x}_1) & \sum_{j=1}^{J-1} {\mathbf c}_{1j} \cdot {\mathbf h}_c^T({\mathbf x}_1)\\
{\mathbf c}_{21}{\mathbf h}_1^T({\mathbf x}_2) & \cdots & {\mathbf c}_{2,J-1} {\mathbf h}_{J-1}^T ({\mathbf x}_2) & \sum_{j=1}^{J-1} {\mathbf c}_{2j} \cdot {\mathbf h}_c^T({\mathbf x}_2)\\
\cdots & \cdots & \cdots & \cdots \\
{\mathbf c}_{m1}{\mathbf h}_1^T({\mathbf x}_m) & \cdots & {\mathbf c}_{m,J-1} {\mathbf h}_{J-1}^T ({\mathbf x}_m) & \sum_{j=1}^{J-1} {\mathbf c}_{mj} \cdot {\mathbf h}_c^T({\mathbf x}_m)
\end{array}\right)
\end{equation}
Then we can reformulate ${\mathbf F}$ into the form of ${\mathbf G}^T {\mathbf W} {\mathbf G}$ so that $|{\mathbf F}|$ becomes a homogeneous polynomial of allocations.

\begin{theorem}\label{thm:Fdiagonal}
	Consider the multinomial logistic model (\ref{logitunifiedmodel}) with independent observations. The Fisher information matrix
	${\mathbf F} = n {\mathbf G}^T {\mathbf W} {\mathbf G}$~.
\end{theorem}

In order to find D-optimal designs, we need to maximize $|{\mathbf G}^T {\mathbf W} {\mathbf G}|$. Since ${\mathbf W}$ is diagonal, we obtain the following theorem as a direct consequence of Theorem~1.1.2 of \cite{fedorov1972} or Lemma~3.1 of \cite{ym2015}:

\begin{theorem}\label{thm:FdetGeneral}
	The determinant of the Fisher information matrix is
	\begin{equation}\label{eq:Fdetcalpha}
	|{\mathbf F}|=n^p|{\mathbf G}^T {\mathbf W} {\mathbf G}|
	=n^p\sum_{\alpha_1\geq 0, \ldots, \alpha_m\geq 0\ :\ \sum_{i=1}^m\alpha_i=p} c_{\alpha_1,\ldots, \alpha_m} \cdot w_1^{\alpha_1} \cdots w_m^{\alpha_m}
	\end{equation}
	with $c_{\alpha_1,\ldots,\alpha_m}=$
	\begin{equation}\label{eq:calpha}
	\sum_{(i_1, \ldots, i_p) \in \Lambda (\alpha_1, \ldots, \alpha_m)} |{\mathbf G}[i_1, \ldots, i_p]|^2 \prod_{k:\alpha_k>0} \prod_{l:(k-1)J<i_l \leqslant kJ} \pi_{k,i_l-(k-1)J}^{-1} \geq 0
	\end{equation}
	where $\alpha_1, \ldots, \alpha_m$ are nonnegative integers,
	$\Lambda (\alpha_1, \ldots, \alpha_m)  = \{(i_1,$ $\ldots,$ $i_p) \mid 1\leq i_1 < \cdots < i_p\leq mJ;
	\# \{l:(k-1)J<i_l \leqslant kJ \}=\alpha_k, k=1, \ldots, m  \}$,
	and ${\mathbf G}[i_1, \ldots, i_p]$ is the submatrix consisting of the $i_1\mbox{th}, \ldots, i_p\mbox{th}$ rows of ${\mathbf G}$.
\end{theorem}

It follows from Theorem~\ref{thm:FdetGeneral}, that the determinant of the Fisher information matrix is an order-$p$ homogeneous polynomial of the proportions $w_1, \ldots, w_m$, with coefficients $c_{\alpha_1, \ldots, \alpha_m}$, which must be nonnegative. As a matter of fact, the structure of the determinant can be further simplified quite significantly. This is given in Lemma~\ref{lemma:max>=J}, Theorem~\ref{thm:alpha>0}, and Corollaries~\ref{col:calpha>0H*} and  \ref{col:calpha>0minimalm}, in the Supplementary Materials (Section~\ref{subsec:alpha}).

\medskip\noindent
{\bf Example~\ref{example:yee2015}.} {\it (continued)}\quad
In this example, the number of factors is $d=4$, and the experimental settings are ${\mathbf x}_i = (x_{i1}, x_{i2}, x_{i3}, x_{i4})^T$, $i=1, \ldots, m$. Since $p_1=p_2=3$, $p_c=2$, and the number of parameters $p=p_1+p_2+p_c=8$, the minimal number of experimental settings is $m=p_1+p_c=5$ by Theorem~\ref{thm:minimalm}. We consider the simplest case $m=5$. That is,
\[
{\mathbf H}^T_1 = {\mathbf H}^T_2 =
\left(\begin{array}{ccc}
1 & x_{11} & x_{12}\\
\vdots & \vdots & \vdots\\
1 & x_{51} & x_{52}
\end{array}\right),
\>\>
{\mathbf H}^T_c =
\left(\begin{array}{cc}
x_{13} & x_{14}\\
\vdots & \vdots\\
x_{53} & x_{54}
\end{array}\right)
\]
By \eqref{eq:Fdetcalpha} of Theorem~\ref{thm:FdetGeneral}, $|{\mathbf F}|$ is an order-$8$ homogeneous polynomial of allocations to the $5$ experimental settings, which may contain up to $(8+5-1)!/(8!(5-1)!)=465$ terms. However, Lemma~\ref{lemma:max>=J} implies $c_{\alpha_1,\ldots, \alpha_5} \neq 0$ only if $\alpha_i \in \{0,1,2\}$. On the other hand, Corollary~\ref{col:calpha>0minimalm} says $c_{\alpha_1,\ldots, \alpha_5} \neq 0$ only if $\#\{i\mid \alpha_i>0\} \geq p_1+p_c=5$, that is, $\alpha_i >0$ for each $i$. Therefore $|{\mathbf F}|$ contains only $5!/(3!2!)=10$ nonzero terms with $\alpha_i \in \{1,2\}$. That is,
\[
|{\mathbf F}| = n^8\cdot \prod_{i=1}^5 w_i \cdot \sum_{1\leq i_1 < i_2 < i_3\leq 5} e_{i_1, i_2, i_3} w_{i_1} w_{i_2} w_{i_3}
\]
for some coefficients $e_{i_1,i_2,i_3}$~. Actually, in terms of the notation $c_{\alpha_1,\ldots,\alpha_m}$ in (\ref{eq:calpha}) of Theorem \ref{thm:FdetGeneral}, $e_{i_1, i_2, i_3} = c_{\alpha_1, \ldots, \alpha_5}$ with $\alpha_i = 2$ for $i \in \{i_1, i_2, i_3\}$ and $1$ otherwise. For example, $e_{1,2,3}=c_{2,2,2,1,1}$~. 
\hfill{$\Box$}

\begin{remark}\label{remark:Fdeterminant}{\rm
		The reformulation of ${\mathbf F}$ in Theorem~\ref{thm:Fdiagonal} enables us to conclude that $|{\mathbf F}|$ is an order-$p$ homogeneous polynomial of allocations  or proportions in Theorem~\ref{thm:FdetGeneral}. It establishes the foundation for an efficient numerical algorithm for searching D-optimal designs (Section~\ref{sec:specifiedX}). On the other hand, the simplification of $|{\mathbf F}|$ based on Lemma~\ref{lemma:max>=J} and Corollary~\ref{col:calpha>0minimalm} allows us to obtain D-optimal designs analytically for minimally supported designs, which is critical for investigating their theoretical properties (Section~\ref{section:minimallysupported}).
		
		With the aid of Theorem~\ref{thm:Fdiagonal}, Theorem~\ref{thm:FdetGeneral} is derived in a more straightforward way. It is broader and with more detailed structures than Theorem~2 in \cite{ytm2016}, which focused on {\it po} models only. Lemma~\ref{lemma:max>=J} and Corollary~\ref{col:calpha>0minimalm} also extend the corresponding results in \citeauthor{ytm2016} (2017, Lemmas~S.4 and S.5 in their Supplementary Materials).
		\hfill{$\Box$}
}\end{remark}

\section{D-optimal Designs}\label{sec:specifiedX}

This section is dedicated to D-optimal designs. Our main focus is on locally D-optimal designs, in which the information matrix is computed at an assumed value of the model parameter $\boldsymbol\theta$. Later, in Section~\ref{sec:ewdesign}, we study EW D-optimal designs under a prior distribution of the parameters. We consider both {\it approximate designs}, which specify proportions of observations at experimental settings without regard to the number of observations $n$, as well as {\it exact designs} where allocation of the $n$ observations to the locations is determined. To start with, we fix a given set of experimental settings ${\mathbf x}_i$, $i=1, \ldots, m$, and consider D-optimal designs $(n_1, \ldots, n_m)$ or $(w_1, \ldots, w_m)$ that maximizes the determinant of the Fisher information matrix $|{\mathbf F}|$. Then we consider the optimization problem over experimental settings as well, using a grid-point search algorithm. While we are thus treating the factor space as {\it discrete} in analytical and computational results, in Section~\ref{sec:gridpoint} we will address the problem of determining D-optimal designs in a {\it continuous} factor space.

\subsection{Design space for multinomial logistic models}\label{sec:designspace}

The experimental settings ${\mathbf x}_i$'s appear in the multinomial logistic model through the categorical probabilities $\pi_{ij}$'s, which are functions of the parameter vector $\boldsymbol\theta$ and the ${\mathbf x}_i$’s.
Our first result examines the restriction imposed on the ${\mathbf x}_i$’s by the requirement $0 < \pi_{ij} < 1, j=1, \ldots, J$. This gives us the collection of all feasible experimental settings, known as the {\it design space}, 
$${\mathcal X} = \left\{ {\mathbf x} = (x_1, \ldots, x_d)^T \mid 0<\pi_j<1, j=1, \ldots, J\right\}$$

\begin{theorem}\label{thm:designspace}
	Let $a_j={\mathbf h}_j^T({\mathbf x})\boldsymbol\beta_j+{\mathbf h}_c^T({\mathbf x})\boldsymbol\zeta$, $j=1, \ldots, J-1$. The design space is
	$
	{\cal X} = \{{\mathbf x} = (x_1, \ldots, x_d)^T \mid a_j \in (-\infty, \infty),\> j=1, \ldots, J-1\}
	$
	for baseline-category, adjacent-categories, and continuation-ratio logit models; and
	$
	{\cal X} = \{{\mathbf x} = (x_1, \ldots, x_d)^T \mid  -\infty < a_1 < a_2 <\cdots < a_{J-1} < \infty\}
	$
	for cumulative logit models.
\end{theorem}

Theorem~\ref{thm:designspace} essentially places no restriction on the design space for models other than cumulative logit models, regardless of the odds structure. Its proof and an illustrative example are provided in the Supplementary Materials (Sections~\ref{sec:proofs} and \ref{subsec:moreexamples}, respectively).

\subsection{D-optimal approximate designs}\label{section:approximatedesign}

Given distinct experimental settings ${\mathbf x}_i \in {\cal X}$, $i=1, \ldots, m$, we look for a D-optimal approximate design ${\mathbf w} = (w_1, \ldots, w_m)^T$ that maximizes the quantity $|{\mathbf G}^T {\mathbf W} {\mathbf G}|$ defined in Theorem~\ref{thm:Fdiagonal}.
The collection of all feasible approximate designs $S=\{(w_1, \ldots, w_m)^T$ $\in$ $\mathbb{R}^m \mid w_i \geq 0, i=1, \ldots, m; \sum_{i=1}^m w_1 = 1\}$ is a bounded closed convex set.
The objective function
$$f({\mathbf w})=|{\mathbf G}^T {\mathbf W} {\mathbf G}|$$
is an order-$p$ homogeneous polynomial by Theorem~\ref{thm:FdetGeneral}. Therefore, a D-optimal approximate design that maximizes $f({\mathbf w})$ must exist. For typical applications, we need designs coming from $S_+ = \{{\mathbf w} \in S\mid f({\mathbf w}) > 0\}$ to avoid degenerate cases. Due to Theorem~\ref{thm:fisher} and the log-concavity of the determinant on positive semi-definite matrices, we know $f({\mathbf w})$ is log-concave \citep{silvey1980, ytm2016} and $S_+$ is convex. A useful result as a corollary of Theorem~\ref{thm:Fpd}, with a proof provided in the Supplementary Materials (Section~\ref{sec:proofs}), is the following:

\begin{corollary}\label{col:fwu>0}  
$S_+$ is nonempty if and only if $f({\mathbf w}_u) > 0$, where ${\mathbf w}_u = (1/m, \ldots, 1/m)^T$ is the uniform allocation. In this case, $f({\mathbf w}) > 0$ for any ${\mathbf w} = (w_1, \ldots, w_m)^T$ such that $0 < w_i < 1, i=1, \ldots, m$.
\end{corollary}

In order to avoid trivial cases, we assume $f({\mathbf w}_u)>0$ from now on.
Following \citeauthor{ytm2016} (2017, Section~3), we define 
\[
f_i(z) = f\left(\frac{w_1(1-z)}{1-w_i}, \ldots, \frac{w_{i-1}(1-z)}{1-w_i}, z, \frac{w_{i+1}(1-z)}{1-w_i}, \ldots, \frac{w_m(1-z)}{1-w_i}\right)
\]
with $0 \leq z \leq 1$ and ${\mathbf w} = (w_1, \ldots, w_m)^T \in S_+$. As a special case of the general equivalence theorem \citep{kiefer1974, pukelsheim1993, atkinson2007, stufken2012, fedorov2014, ymm2016, ytm2016}, ${\mathbf w}$ is D-optimal if and only if $f_i(z)$ attains its maximum at $z=w_i$ for each $i=1, \ldots, m$. Actually, $f_i(z) = n^{-p}|{\mathbf F}|$ for the approximate design ${\mathbf w} + ({\mathbf e}_i - {\mathbf w})\cdot (z-w_i)/(1-w_i)$, where ${\mathbf e}_i$ is an $m\times 1$ vector with the $i$th coordinate $1$ and all others $0$. Parallel to Theorem~6 in \cite{ytm2016}, we obtain an explicit formula for $f_i(z)$ as an order-$p$ polynomial of $z$ (see Theorem~\ref{thm:f_i(z)} in the Supplementary Materials). 

In order to find D-optimal designs numerically, we use a lift-one algorithm (Section~\ref{sec:algorithm} in the Supplementary Materials), which is essentially the same as the one in \cite{ytm2016} for cumulative link models. The lift-one algorithm is of general-equivalence-theorem type \citep{ymm2016}. Its convergence to a global maximum is guaranteed \citep{ym2015}.

We use the following real experiment as an illustration. In the application (Zocchi and Atkinson, 1999) a uniform design and a continuation-ratio logit model with {\it npo} component was used for the analysis. As noted earlier, the theoretical results for the {\it npo} model, as special cases of the results for the {\it ppo} model, are given in the Supplementary Materials (Section~\ref{subsec:npomodels}).

\begin{example}\label{example:flies} {\bf Emergence of house flies} {\rm
	\cite{atkinson1999} discussed an experiment on emergence of house flies originally reported by \cite{itepan1995}.  In this experiment, $n=3,500$ pupae were grouped evenly into seven sets and exposed to seven doses of radiation (in units Gy), $x_i = 80, 100, 120, 140,$ $160,$ $180,$ $200$, respectively. After a period of time, for each set of $n_i=500$ pupae, the summarized responses are the number $y_{i1}$ of flies that died before the opening of the pupae ({\tt unopened pupae}), the number $y_{i2}$ of flies out of {\tt opened pupae} but died before complete emergence, and the number $y_{i3} = n_i - y_{i1} - y_{i2}$ of flies out of {\tt opened pupae} and completely emerged. Following \cite{atkinson1999}, we assume that the reponses of the $n_i=500$ pupae in the $i$th set are independent and follow the same distribution. Then the summary responses ($y_{i1}, y_{i2}, y_{i3}$) follow a multinomial distribution and have a clearly nested or hierarchical structure (see Table~1 of \cite{atkinson1999} for the experimental data).
	We confirm that the continuation-ratio logit model with {\it npo} component fits the data the best in terms of AIC and BIC (see Table~\ref{tab:compare5} in the Supplementary Materials); this model was adopted by \cite{atkinson1999} as follows
	\[
	\log\left(\frac{\pi_{i1}}{\pi_{i2} + \pi_{i3}}\right) = \beta_{11}+\beta_{12} x_i+\beta_{13} x_i^2,
	\quad
	\log\left(\frac{\pi_{i2}}{\pi_{i3}}\right) = \beta_{21}+\beta_{22} x_i
	\]
	The model has 5 parameters with fitted values $\hat{\boldsymbol\beta} = (\hat{\beta}_{11}, \hat{\beta}_{12}, \hat{\beta}_{13}, \hat{\beta}_{21}, \hat{\beta}_{22})^T$ $=$ $(-1.935,$ $-0.02642,$ $0.0003174, -9.159, 0.06386)^T$.

	\begin{table}[ht]
		\caption{Exact and Approximate Allocations for House Flies Experiment}
		\centering 
		\footnotesize
		\begin{tabular}{l|rrrrrrr} 
			\hline\hline
			Dose of radiation (Gy) &80  &100 &120 &140 &160 &180 &200\\ \hline
			Original allocation & 500 & 500 & 500 & 500 & 500 & 500 & 500 \\
			D-optimal exact & 1091  & 0 & 1021 & 374 & 1014 & 0 & 0  \\ \hline
			Original proportion & .1429  & .1429 & .1429 & .1429 & .1429 & .1429 & .1429\\
			D-optimal approximate & .3116  & 0 & .2917 & .1071 & .2896 & 0 & 0 \\
			Bayesian D-optimal & .3159  & .0000 & .2692 & .1160 & .2990 & .0000 & .0000 \\
			EW D-optimal & .3120  & 0 & .2911 & .1087 & .2882 & 0 & 0 \\
			\hline\hline 
		\end{tabular}
		\label{tab:approx_opt1}
		\normalsize
	\end{table}

	Considering a followup experiment with the fitted parameter values as the assumed values, we obtain the locally D-optimal approximate design (Table~\ref{tab:approx_opt1}) using the lift-one algorithm. The efficiency of the original uniform allocation is $(|{\mathbf F}_{original}|/|{\mathbf F}_{D-opt}|)^{1/5} = 83.1\%$. The D-optimal approximate design only requires 4 rounds (doses) of experiments. It is interesting to note that in this example the D-optimal design does not allocate observations to high doses of radiation; we will revisit this issue in Section~\ref{sec:discussion}. 
	The D-optimal exact design, Bayesian and EW D-optimal designs listed in Table~\ref{tab:approx_opt1} will be described in Sections~\ref{section:exactdesign} and \ref{sec:ewdesign}, respectively. We will return to this example later.
}\hfill{$\Box$}
\end{example}

\subsection{D-optimal exact designs}\label{section:exactdesign}

In practice, a design should specify an integer-valued allocation of the $n$ experimental units to the $m$ distinct experimental settings ${\mathbf x}_i$, $i=1, \ldots, m$. An integer-valued allocation ${\mathbf n} = (n_1, \ldots, n_m)^T$ that maximizes $|\sum_{i=1}^m n_i {\mathbf F}_i|$ as defined in (\ref{eqn:F}) with $\sum_{i=1}^m n_i = n$ is known as a {\it D-optimal exact design}.

Although different rounding algorithms have been proposed in the literature to obtain an exact allocation from an optimal approximate allocation, an algorithm searching for optimal exact designs directly is still needed (see \citeauthor{ytm2016}~(2017, Section~4) and reference therein). For simplicity, we denote the objective function as $f({\mathbf n}) = f(n_1, \ldots, n_m) = |\sum_{i=1}^m n_i {\mathbf F}_i|$ if there is no ambiguity. Following \cite{ymm2016} and \cite{ytm2016}, we define $f_{ij}(z) = f(n_1, \ldots, n_{i-1}, z, n_{i+1}, \ldots, n_{j-1}, n_i + n_j - z, n_{j+1}, \ldots, n_m)$ with $z= 0, 1, \ldots, n_i+n_j$ given $1\leq i < j \leq m$ and ${\mathbf n} = (n_1, \ldots, n_m)^T$.
Similar to Theorem~9 in \cite{ytm2016}, we can obtain an explicit formula for $f_{ij}(z)$ (see Theorem~\ref{thm:f_ij(z)} in the Supplementary Materials) from Theorem~\ref{thm:FdetGeneral}, Lemma~\ref{lemma:max>=J} and Corollary~\ref{col:calpha>0minimalm}.

Given $n>0$, we assume that there exists an exact allocation ${\mathbf n} = (n_1, \ldots, n_m)^T$ such that $f({\mathbf n}) > 0$. 
An exchange algorithm (Section~\ref{sec:algorithm} in the Supplementary Materials), which is essentially the same as the one in \cite{ytm2016}, is used for obtaining an exact design. Note that the integer-valued allocation found by the exchange algorithm is not guaranteed to be D-optimal, especially when $n$ is relatively small compared with $m$ \citep{ymm2016}.

The following real experiment shows that a minimally supported design can be D-optimal for a cumulative logit model with {\it npo} assumption.

\begin{example}\label{example:trauma}  {\bf Trauma clinical trial} {\rm
	\cite{chuang1997} studied a dataset of trauma patients that has five ordered response categories, {\tt death}, {\tt vegetative} {\tt state}, {\tt major disability}, {\tt minor disability}, and {\tt good} {\tt recovery}, describing their clinical outcomes. These five categories are often called the Glasgow Outcome Scale (GOS) in the literature of critical care \citep{jennett1975}. 
	\cite{agresti2010} found strong evidence against the assumption of proportional odds for the trauma data. We fit the data with different logit models, as well as different odds assumptions, and confirm that the cumulative logit model with {\it npo} assumption fits the trauma data the best in terms of AIC and BIC (see Table~\ref{tab:compare1} in the Supplementary Materials). Thus the model is,
	\begin{equation}\label{eq:traumamodel}
	\log\left(\frac{\pi_{i1}+\cdots + \pi_{ij}}{\pi_{i,j+1}+\cdots + \pi_{i5}}\right) = \beta_{j1}+\beta_{j2} x_i~,\quad i=1,2,3,4; j=1,2,3,4
	\end{equation}
	The fitted parameter values are $\hat{\boldsymbol\beta} = (\hat{\beta}_{11}, \hat{\beta}_{12}, \hat{\beta}_{21}, \hat{\beta}_{22},
	\hat{\beta}_{31},
	\hat{\beta}_{32}, \hat{\beta}_{41},  \hat{\beta}_{42})^T = (-0.865,$ $-0.113,$ $-0.094, -0.269, 0.706, -0.182, 1.909, -0.119)^T$. The clinical trial with 802 patients assigned 210, 190, 207, 195 patients to the four treatment groups, {\tt Placebo}, {\tt Low dose}, {\tt Medium dose}, and {\tt High dose}, respectively, which is roughly uniform on dosage. If we treat the fitted parameter values as the assumed values for local optimality, using the exchange algorithm, we find the D-optimal exact design (401, 0, 0, 401), which is minimally supported and uniform on its supporting points. If a followup experiment is considered and the fitted parameter values are the true values, the efficiency of the original allocation compared with our D-optimal design is only $74.7\%$.
}\hfill{$\Box$}
\end{example}

\medskip\noindent
{\bf Example~\ref{example:flies}.} {\bf Emergence of house flies} ({\it continued})
In this experiment, the allocation of experimental units used by the scientists is 500 for each of the seven doses of radiation. Using the exchange algorithm, we obtain the D-optimal exact allocation (Table~\ref{tab:approx_opt1}), which is similar to the D-optimal approximate design. Actually, the proportions of the D-optimal exact design $(1091, 1021, 374, 1014)/3500 = (0.3117, 0.2917, 0.1069, 0.2897)$ on the four support points are roughly the same as the ones in the D-optimal approximate design. If the fitted parameter values are the true values, the efficiency of the allocation used in the experiment is $83.1\%$ compared with the D-optimal exact allocation.
\hfill{$\Box$}

\begin{remark}\label{remark:twostage}{\rm
Our D-optimal exact designs are ``locally" D-optimal since they require assumed parameter values. When there is no pilot study available for a good parameter estimate, D-optimal designs can still be used in a multistage design (see \cite{khuri2006} for a review). In the trauma clinical trial, for example, we may consider a two-stage design with 802 patients in total. At Stage~1, we may assign 240 patients (about 30\%) randomly and uniformly to the four treatment groups. After collecting the outcomes of the patients from Stage~1, we obtain the estimated coefficients and then find the corresponding D-optimal exact design for the rest 562 patients, known as Stage~2. The parameter estimates from the two-stage design are based on the outcomes of all the 802 patients. In terms of root mean squared error  (RMSE) from the assumed parameter values, the two-stage design achieves $0.149$ for the intercepts $\beta_{11}, \ldots, \beta_{41}$ and $0.053$ for the slopes $\beta_{12}, \ldots, \beta_{42}$, on average over 100 simulations. The corresponding average RMSEs from the original design are $0.190$ and $0.070$, respectively. The reductions in RMSEs by using our D-optimal exact designs in the two-stage design are both significant.       
}\hfill{$\Box$}
\end{remark}

\subsection{EW D-optimal designs}\label{sec:ewdesign}

The D-optimal approximate and exact designs discussed thus far are locally D-optimal designs using  assumed parameter values. Bayesian D-optimality \citep{chaloner1995}, which maximizes $E(\log|{\mathbf F}|)$ with a given prior distribution on the unknown parameters, provides an alternative approach. A drawback of the Bayesian approach is its computational intensity since the objective function deals with multiple integrals. An alternative solution is the EW D-optimality \citep{atkinson2007}, also see \cite{ymm2016, ytm2016}, which maximizes $\log |E({\mathbf F})|$ or $|E({\mathbf F})|$ instead.
Among all of the criteria proposed by \citeauthor{atkinson2007} (2007, Table~18.1) as surrogates for Bayesian D-optimality, including $-\log E (|{\mathbf F}|^{-1})$, $-\log |E({\mathbf F}^{-1})|$, and $\log E|{\mathbf F}|$, EW D-optimal design requires the minimum computation. \cite{ymm2016} showed that an EW D-optimal design could be highly efficient in terms of the Bayesian criterion in comparison to the Bayesian D-optimal design, while computational time is essentially the same as a locally D-optimal design. \cite{ytm2016} also used EW-criterion for cumulative link models with proportional odds and confirmed its high-efficiency.

By Theorem~\ref{thm:fisher} and Corollary~\ref{col:Fi}, the Fisher information matrix in our case is ${\mathbf F} = \sum_{i=1}^m n_i {\mathbf X}_i^T {\mathbf U}_i {\mathbf X}_i$, where ${\mathbf U}_i$ consists of $u_{st}({\boldsymbol\pi}_i)$, the only components involving model parameters. In order to calculate $E({\mathbf F})$ with respect to a prior distribution on parameters, we first calculate $E(u_{st}({\boldsymbol\pi}_i))$ and then use the results and algorithms developed for locally D-optimal designs to get EW D-optimal designs.

We provide formulae in both the proof of Theorem~\ref{thm:designspace} and Section~\ref{sec:pifromxi} (in matrix form) in the Supplementary Materials for calculating $\pi_{ij}$'s given ${\mathbf X}_i$'s and the parameter values. We also provide formulae in Appendix~\ref{subsec:piij} for calculating $u_{st}({\boldsymbol\pi}_i) = {\mathbf c}_{is}^T {\rm diag}({\boldsymbol\pi}_i)^{-1} {\mathbf c}_{it}$'s given $\pi_{ij}$'s.

By Theorem~\ref{thm:designspace}, for baseline-category, adjacent-categories, and continuation-ratio logit models, there is essentially no restriction on values of parameters for a given design point ${\mathbf x} = (x_1, \ldots, x_d)^T$. The corresponding $E(u_{st}({\boldsymbol\pi}_i))$ can be obtained by integrating over the same range of parameter values for different design points. However, for cumulative logit models, where the parameters should satisfy $a_1 < a_2 < \cdots < a_{J-1}$ with $a_j={\mathbf h}_j^T({\mathbf x})\boldsymbol\beta_j+{\mathbf h}_c^T({\mathbf x})\boldsymbol\zeta$, the domain of integration is typically not rectangular.

\medskip\noindent
{\bf Example~\ref{example:trauma}.} {\bf Trauma clinical trial} ({\it continued})\quad In this example, the cumulative logit model \eqref{eq:traumamodel} with non-proportional odds was adopted. Given the predetermined set ${\mathcal X} = \{1,2,3,4\}$ consisting of $m=4$ design points, the feasible parameter space is $\boldsymbol\Theta = \{ \boldsymbol\theta = (\beta_{11}, \beta_{12}, \beta_{21}, \beta_{22}, \beta_{31}, \beta_{32}, \beta_{41},  \beta_{42})^T \mid \beta_{11} + \beta_{12} x < \beta_{21} + \beta_{22} x < \beta_{31} + \beta_{32} x < \beta_{41} + \beta_{42} x, \mbox{ for }x\in {\mathcal X}\}$, which is not rectangular. \cite{atkinson1999} used a multivariate normal prior with its variance-covariance matrix estimated from an initial dataset. For the computations in this example, we bootstrap the 802 observations from the initial dataset for 1000 times and denote the corresponding fitted parameters as $\boldsymbol\theta_1, \ldots, \boldsymbol\theta_{1000}$. Then an estimate of the Bayesian criterion $\phi({\mathbf w}) = E(\log|{\mathbf F}({\mathbf w}, \boldsymbol\theta)|)$ for design ${\mathbf w} = (w_1, \ldots, w_m)^T$ is
\[
\hat\phi({\mathbf w}) = \frac{1}{1000}\sum_{j=1}^{1000} \log|{\mathbf F}({\mathbf w}, \boldsymbol\theta_j)| = \frac{1}{1000} \sum_{j=1}^{1000} \log |\sum_{i=1}^m w_i {\mathbf F}_i(\boldsymbol\theta_j)|
\]
The Bayesian D-optimal design listed in Table~\ref{tab:trauma_opt1} maximizes $\hat\phi({\mathbf w})$, while the EW D-optimal design maximizes $|\sum_{i=1}^m w_i \hat{E}({\mathbf F}_i)|$ with estimated entry-by-entry expectation $\hat{E}({\mathbf F}_i) = \sum_{j=1}^{1000} {\mathbf F}_i(\boldsymbol\theta_j)/1000$. The two designs are essentially the same in this case.
\hfill{$\Box$}

\begin{table}[ht]
\caption{D-optimal Designs for Trauma Clinical Trial}
\centering 
\footnotesize
\begin{tabular}{l|rrrr} 
	\hline\hline
	Design point $x$ & 1  & 2 & 3 & 4\\
	Treatment group & Placebo & Low & Medium & High\\ \hline
	Original allocation & 210  & 190 & 207 & 195\\
	D-optimal exact & 401  & 0 & 0 & 401\\ \hline
	Original proportion & 0.2618  & 0.2369 & 0.2581 & 0.2431\\
	D-optimal approximate & 0.5  & 0 & 0 & 0.5\\
	Bayesian D-optimal & 0.4997  & $1.0\times 10^{-8}$ & $2.5\times 10^{-8}$ & 0.5003\\
	EW D-optimal & 0.5  & 0 & 0 & 0.5\\
	\hline\hline 
\end{tabular}
\label{tab:trauma_opt1}
\normalsize
\end{table}

\medskip\noindent
{\bf Example~\ref{example:flies}.} {\bf Emergence of house flies} ({\it
continued})\quad Similar to the trauma example, we bootstrap the original observations 1000 times and obtain the corresponding Bayesian and EW D-optimal designs (see Table~\ref{tab:approx_opt1}). In this case, EW D-optimal design is close to the Bayesian design. In order to check the robustness of these designs towards misspecified parameter values, we consider the same set of bootstrapped parameter vectors $\boldsymbol\theta_1, \ldots, \boldsymbol\theta_{1000}$. For $j=1, \ldots, 1000$, we apply our lift-one algorithm to find the corresponding D-optimal allocation ${\mathbf p}_j$ with $\boldsymbol\theta_j$ as the assumed parameter value. The efficiencies of a target design ${\mathbf p}$ with respect to ${\mathbf p}_j$ is defined as $(|{\mathbf F}({\mathbf p}, \boldsymbol\theta_j)|/|{\mathbf F}({\mathbf p}_j, \boldsymbol\theta_j)|)^{1/p}$ with $p=5$ parameters in this case. The summary statistics of efficiencies in Table~\ref{tab:flies_efficiency} show that the Bayesian and EW D-optimal designs are highly robust in terms of parameter misspecification; the EW design is slightly better. Both of them are much better than the uniform design used in the study.
\hfill{$\Box$}

\begin{table}[ht]
\caption{Efficiencies of Designs for House Flies Experiment}
\centering 
\footnotesize
\begin{tabular}{l|rrrrr} 
	\hline\hline
	Design & Min & 1st Quartile & Median & 3rd Quartile & Max\\ \hline
	Bayesian D-optimal & 0.9912 & 0.9985 & 0.9989 & 0.9992 & 0.9998 \\
	EW D-optimal       & 0.9934 & 0.9991 & 0.9996 & 0.9998 & 1.0000 \\
	Uniform design     & 0.7762 & 0.8181 & 0.8304 & 0.8445 & 0.8861 \\
	\hline\hline 
\end{tabular}
\label{tab:flies_efficiency}
\normalsize
\end{table}

\subsection{Grid point search for continuous factors}\label{sec:gridpoint}

Both the lift-one and exchange algorithms are appropriate for searching optimal designs with a predetermined finite set of experimental settings. When a factor is continuous, one common practice is to partition the continuous region of the factor levels into finite subintervals and consider only the grid points as its discrete levels.
It works reasonably well with a moderate number of continuous factors (see \cite{yangmin2013} for a discussion on the efficiency of optimal designs based on grid-point experimental settings).

Once a set of grid points is chosen for each continuous factor, the design problem becomes an allocation problem on a finite set of design points. Hence, all of the previous algorithms and results can be applied. We use the example of house flies for the purpose of illustration.

\medskip\noindent
{\bf Example~\ref{example:flies}.} {\bf Emergence of house flies} ({\it
continued}) The sole factor in this experiment, dose of radiation, is by nature a continuous factor. In the emergence of house flies experiment, seven levels ranging equidistantly from 80 to 200 were used, that is, with grid distance 20. Suppose a followup experiment is considered and dose levels with grid distance 5 are feasible (which technically depends on the sensitivity of the radiation device). Then there are 25 dose levels available as 80, 85, 90, $\dots$, 195, 200. The D-optimal approximate design given the 25 grid-5 design concentrates on five design points 80, 120, 125, 155, 160 with optimal allocation $0.3163, 0.1429, 0.2003,$ $0.1683,$ $0.1723$, respectively. If we further consider the set of 121 grid-1 dose levels, the D-optimal design is supported on five design points 80, 122, 123, 157, 158 with allocations $0.3163, 0.0786, 0.2636, 0.2206, 0.1209$. The optimal design seems to converge  to a three-point design as the grid points become finer and finer. Actually, if we reallocate the grid-1 D-optimal design, denoted as ${\mathbf p}_1$, into a minimally supported design ${\mathbf p}_{1m}$ at dose levels 80, 123, 157 with weights $0.3163, 0.3422, 0.3415$, then the efficiency of ${\mathbf p}_{1m}$ compared with ${\mathbf p}_{1}$ is as high as $99.99\%$.
\hfill{$\Box$}

\subsection{Minimally supported designs}\label{section:minimallysupported}

An important question in design theory is, what is the least number of design points we need to keep ${\mathbf F}$ positive definite? In other words, what is the number of distinct experimental settings for a {\it minimally supported design}?
Theorem~\ref{thm:minimalm} provides a lower bound, which has been denoted as $k_{\min}$. For typical models proposed in the literature, for example, main-effects models (see Example~\ref{ex:maineffects}), $k_{\min}$ is the minimal number of experimental settings since the conditions in Theorem~\ref{thm:minimalm} are satisfied. Answers to general cases rely on the specific forms of the predictor functions ${\mathbf h}_j$ and ${\mathbf h}_c$ if applicable.

Another question is whether a uniform allocation is D-optimal among minimally supported designs. The answer is known to be ``Yes" for $J=2$ since all of the four logit models are equivalent to the usual logistic model for binary responses.

Nevertheless, for {\it po} models, \cite{ytm2016} showed analytically that uniform allocations are not D-optimal among minimally supported designs for cumulative link. We use Example~\ref{example:d1J3m2po} in the Supplementary Materials to show that the conclusion is the same for {\it po} models with baseline-category, adjacent-categories, and continuation-ratio logit links as well. Actually, the objective functions take the same form for the four different logit models. Since {\it po} models are special cases of {\it ppo} models, we conclude that uniform allocations are not D-optimal in general for {\it ppo} models as well.

As for {\it npo} models with $J\geq 3$, in general, uniform allocations are not D-optimal either (see Section~\ref{subsec:uniformDoptimal} in the Supplementary Materials).
Nevertheless, for ``regular'' {\it npo} models, that is, when $p_1= \cdots = p_{J-1}$, uniform allocations are D-optimal among minimally supported designs (see Corollary~\ref{col:npouniform} in the Supplementary Materials), which theoretically confirms the D-optimal design for trauma clinical trial in Example~\ref{example:trauma}.

\section{Discussion}\label{sec:discussion}

A criticism of optimal designs is that their rigid reliance on the assumed model does not allow exploration of different models; this is mainly due to the frugality of the number of experimental settings. Nevertheless the results can be used to identify efficient designs that allow flexibility. For instance, in the trauma clinical trial example (Example~\ref{example:trauma}), the D-optimal exact design assigns 401 observations to the placebo (1) and high dose (4), but none to the low (2) or medium (3) levels. If instead we allocate 361 observations to each of 1 and 4, and 40 observations to each of 2 and 3, we gain some flexibility with a loss of efficiency (in comparison to the D-optimal design) of only 5\%. This is further discussed in the Supplementary Materials (Table~\ref{tab:trauma.list}). 

In Example~\ref{example:flies}, the support of the D-optimal design did not include some of the extreme points of the experimental region. This is a departure from the nature of D-optimal designs for univariate responses, and is something that will be studied in the future.

For EW and Bayesian D-optimal designs, the choice of the prior on the parameter space $\boldsymbol\Theta$ is critical. Suppose the parameters are $\theta_1, \ldots, \theta_p$ with individual ranges $\theta_i \in I_i$.  As a common practice for generalized linear models, a uniform prior or independent normal prior may be assumed for a rectangular domain $I_1 \times \cdots \times I_p$ \citep{ymm2016}. One issue with multinomial logit models is that the feasible domain of $\boldsymbol\Theta$ may not be rectangular, at least for cumulative logit models. For the two examples in Section~\ref{sec:ewdesign}, we use bootstrapping to obtain an empirical prior. \cite{atkinson1999} adopted a multivariate normal prior with mean vector and covariance structure estimated from previous experimental data. For the two examples we study, our results are essentially the same when we use the multivariate normal prior. However, when we use uniform or independent normal prior, the EW D-optimal design is not as good as the Bayesian D-optimal design in terms of robustness. We will further investigate the choice of prior and criteria other than EW optimality.

When a pilot study or experimental data is not available for locally or EW optimal designs, sequential design or multistage design may be used to obtain an initial guess or a reasonable prior for unknown parameters. In Remark~\ref{remark:twostage}, we provide a two-stage design for the trauma clinical trial as an example. Further investigation needs to be done towards best multistage designs for multinomial logistic models. 

When confronted with model uncertainty, {\it npo} versus {\it ppo} for instance, our work may provide the experimenter the option to choose a design, for instance, that is highly efficient for both {\it ppo} and {\it npo} models. Further investigations along this line would be practically useful.

In the literature, model parameters other than  $\boldsymbol\theta=(\boldsymbol\beta_1^T, \cdots,\boldsymbol\beta_{J-1}^T,\boldsymbol\zeta)^T$ for linear predictors $\eta_{ij}={\mathbf h}_j^T({\mathbf x}_i)\boldsymbol\beta_j+{\mathbf h}_c^T({\mathbf x}_i)\boldsymbol\zeta$ were also used, for example, in \cite{perevozskaya2003}. It can be verified that D-optimal designs will not be affected by the choice of parameters (see Section~\ref{subsec:reparametrization} in the Supplementary Materials).


\appendix

\section{Appendix section}\label{app}

\subsection{Constant matrix ${\mathbf L}$ in equation~\eqref{logitunifiedmodel}}\label{subsec:L}
\footnotesize
\[
{\mathbf L}_{baseline}= \begin{pmatrix}
1&& & & 0\\
&1& & & 0\\
&& \ddots& & \vdots\\
&& & 1& 0\\
0&0&\cdots &0 &1\\
0&0&\cdots &0 &1\\
\vdots&\vdots&\ddots &\vdots &\vdots\\
0&0&\cdots &0 &1\\
1&1&\cdots &1 &1
\end{pmatrix},\>
{\mathbf L}_{cumulative}= \begin{pmatrix}
1&0&\cdots & 0& 0\\
1&1&0 &\cdots & 0\\
\vdots&\vdots& \ddots& & \vdots\\
1&1&\cdots & 1& 0\\
0&1&\cdots &1 &1\\
0&0&1&\cdots  &1\\
\vdots&\vdots&\ddots &\ddots &\vdots\\
0&0&\cdots &0 &1\\
1&1&\cdots &1 &1
\end{pmatrix}
\]
\[
{\mathbf L}_{continuation}= \begin{pmatrix}
1&& & & 0\\
&1& & & 0\\
&& \ddots& & \vdots\\
&& & 1& 0\\
0&1&\cdots &\cdots &1\\
0&0&1&\cdots  &1\\
\vdots&\vdots&\ddots &\ddots &\vdots\\
0&0&\cdots &0 &1\\
1&1&\cdots &1 &1
\end{pmatrix},\>
{\mathbf L}_{adjacent}= \begin{pmatrix}
1&& & & 0\\
&1& & & 0\\
&& \ddots& & \vdots\\
&& & 1& 0\\
0&1&& & \\
0&&1& &\\
\vdots&&& \ddots \\
0&&& & 1\\
1&1&\cdots &1 &1
\end{pmatrix}
\]
\normalsize

\subsection{Formula towards ${\cal M}(\mathbf H_c^T)\cap ( \cap_{j=1}^{J-1}{\cal M}(\mathbf H_j^T))$ in Theorem~\ref{thm:minimalm}}\label{subsec:dim}

\begin{theorem}\label{thm:rankintersection}
Suppose ${\mathbf H}_i$ is of $p_i\times m$ with rank $r_i$, $i=1, \ldots, n$. Denote $rank\left(({\mathbf H}_{i_1}^T,\right.$ $\ldots,$ $\left.{\mathbf H}_{i_k}^T)\right) = r_{i_1,\ldots, i_k}$ for any $2\leq k\leq n$ and $1\leq i_1<\cdots <i_k\leq n$. Then
\begin{equation}\label{eq:nintersections}
dim\left(\bigcap_{i=1}^n {\cal M}({\mathbf H}_i^T)\right)
= \sum_{i=1}^n r_i - \sum_{i_1<i_2} r_{i_1,i_2} + \sum_{i_1<i_2<i_3} r_{i_1,i_2,i_3} - \cdots + (-1)^{n-1} r_{1,2,\ldots,n}
\end{equation}
\end{theorem}

The proof of Theorem~\ref{thm:rankintersection} is relegated to the Supplementary Materials. Note that ${\cal M}(\mathbf H_c^T)\cap ( \cap_{j=1}^{J-1}{\cal M}(\mathbf H_j^T)) = \{0\}$ if and only if its dimension is 0.

\subsection{Formulae for calculating $u_{st}({\boldsymbol\pi}_i)$'s given $\pi_{ij}$'s}\label{subsec:piij}

\begin{theorem}\label{thm:ustpi}
	Consider the multinomial logit model~\eqref{logitunifiedmodel}. For $i=1, \ldots, m$,
	\begin{itemize}
		\item[(i)] $u_{st}({\boldsymbol\pi}_i) = u_{ts}({\boldsymbol\pi}_i)$, $s,t=1, \ldots, J$;
		\item[(ii)] $u_{sJ}({\boldsymbol\pi}_i) = 0$ for $s=1, \ldots, J-1$ and $u_{JJ}({\boldsymbol\pi}_i) = 1$;
		\item[(iii)] For $s=1, \ldots, J-1$,
		\[
		u_{ss}({\boldsymbol\pi}_i) =\left\{\begin{array}{cl}
		\pi_{is}(1-\pi_{is}) & \mbox{for baseline-category},\\
		\gamma_{is}^2(1-\gamma_{is})^2(\pi_{is}^{-1} + \pi_{i,s+1}^{-1}) & \mbox{for cumulative},\\
		\gamma_{is}(1-\gamma_{is}), & \mbox{for adjacent-categories},\\
		\pi_{is}(1-\gamma_{is})(1-\gamma_{i,s-1})^{-1} & \mbox{for continuation-ratio};
		\end{array}\right.
		\]
		\item[(iv)] For $1\leq s < t \leq J-1$,
		\[
		u_{st}({\boldsymbol\pi}_i) =\left\{\begin{array}{cl}
		-\pi_{is} \pi_{it} & \mbox{for baseline-category},\\
		-\gamma_{is}\gamma_{it}(1-\gamma_{is})(1-\gamma_{it})\pi_{it}^{-1} & \mbox{for cumulative}, t-s=1,\\
		0 & \mbox{for cumulative}, t-s>1,\\
		\gamma_{is}(1-\gamma_{it}) & \mbox{for adjacent-categories},\\
		0 & \mbox{for continuation-ratio};
		\end{array}\right.
		\]
	\end{itemize}
	where $\gamma_{ij} = \pi_{i1} + \cdots + \pi_{ij}$, $j=1, \ldots, J-1$; $\gamma_{i0}\equiv 0$ and $\gamma_{iJ}\equiv 1$.
\end{theorem}

\section*{Supplementary Materials}

It contains:
(1) list of notations;
(2) formulae of matrix differentiation;
(3) explicit forms of $({\mathbf C}^T{\mathbf D}_i^{-1}{\mathbf L})^{-1}$ for all the four logit models;
(4) positive definiteness of ${\mathbf U}$;
(5) row rank of ${\mathbf H}$ matrix;
(6) results on the coefficient $c_{\alpha_1, \ldots, \alpha_m}$ for simplifying $|{\mathbf F}|$;
(7) expressions for {\it po} models;
(8) expressions for {\it npo} models;
(9) model selections for Trauma Clinical Trial and Emergence of House Flies;
(10) lift-one and exchange algorithms;
(11) formulae for calculating $\pi_{ij}$'s from ${\mathbf X}_i$'s;
(12) reparametrization and D-optimality;
(13) more discussion on D-optimality of uniform designs;
(14) more examples;
(15) proofs of Theorems~\ref{thm:fisher}, \ref{thm:fisher_all}, \ref{thm:Upd}, \ref{thm:vdet}, \ref{thm:rankH}, \ref{thm:rankintersection}, \ref{thm:Fdiagonal}, \ref{thm:alpha>0}, \ref{thm:designspace}, \ref{thm:f_ij(z)}, \ref{thm:J=2}, \ref{thm:ustpi}, Corollaries~\ref{lemma:ppo_m}, \ref{col:fwu>0}, \ref{col:npouniform}, and Lemmas~\ref{lemma:pii}, \ref{lemma:vi=CDL}, \ref{lemma:max>=J},
\ref{claim4.1}, \ref{claim4.2}, \ref{claim4.3}, \ref{claim4.4}, \ref{lemma:solutiontoflies}.


\clearpage
\setcounter{page}{1}
\def\thepage{S\arabic{page}}

\fontsize{10.95}{14pt plus.8pt minus .6pt}\selectfont
\vspace{0.8pc}
\centerline{\large\bf D-optimal Designs for Multinomial Logistic Models}
\vspace{.25cm}
\centerline{Xianwei Bu$^{1}$, Dibyen Majumdar$^{2}$ and Jie Yang$^{2}$}
\vspace{.4cm}
\centerline{\it  $^1$AbbVie Inc.~and $^2$University of Illinois at Chicago}
\vspace{.55cm}
 \centerline{\bf Supplementary Materials}
\vspace{.55cm}
\fontsize{9}{11.5pt plus.8pt minus .6pt}\selectfont
\par

\renewcommand{\thesection}{S}
\setcounter{equation}{0}
\setcounter{subsection}{0}
\renewcommand{\theequation}{S.\arabic{equation}}

\normalsize

\subsection{List of notations}\label{subsec:list}

\begin{list}{}{
	\setlength{\labelwidth}{0.5in}
    \setlength{\leftmargin}{0.6in}
    \setlength{\labelsep}{.3in}
    \setlength{\rightmargin}{\leftmargin}
    }
\item[${\mathbf 0}_k$\hfill] A vector of $k$ zeros
\item[$a_j$\hfill] ${\mathbf h}_j^T({\mathbf x})\boldsymbol\beta_j+{\mathbf h}_c^T({\mathbf x})\boldsymbol\zeta$, $j=1, \ldots, J-1$, given ${\mathbf x} = (x_1, \ldots, x_d)^T$
\item[$b_j$\hfill] Coefficients in representing $f_i(z)$, $j=0, \ldots, J-1$
\item[${\mathbf B}_J$\hfill] $J\times J$ constant matrix used for deriving the coefficients of $f_i(z)$, $(s^{t-1})_{st}$
\item[${\mathbf C}$\hfill] $J\times (2J-1)$ constant matrix, same for all the four logit models
\item[${\mathbf c}$\hfill] Vector used for deriving coefficients of $f_i(z)$, $(c_1, \ldots, c_{J-1})^T$
\item[${\mathbf c}_{ij}$\hfill] $J\times 1$ vectors such that $({\mathbf C}^T{\mathbf D}_i^{-1}{\mathbf L})^{-1} = ({\mathbf c}_{i1}, \ldots, {\mathbf c}_{iJ})$
\item[$c_j$\hfill] $(j+1)^p j^{J-1-p} f_i(1/(j+1)) - j^{J-1}f_i(0)$, $j=1, \ldots, J-1$
\item[$c_{\alpha_1, \ldots, \alpha_m}$\hfill] Coefficient of $w_1^{\alpha_1} \cdots w_m^{\alpha_m}$ in the determinant of ${\mathbf G}^T {\mathbf W} {\mathbf G}$
\item[$d$\hfill] Total number of design factors
\item[$d_s$\hfill] $d_s=(f_{ij}(s)-f_{ij}(0))/s$, $s=1, \ldots, q$, for coefficients in $f_{ij}$
\item[${\mathbf D}_i$\hfill] ${\rm diag}({\mathbf L}{\boldsymbol\pi}_i)$
\item[${\mathbf e}_i$\hfill] $m\times 1$ vector with the $i$th coordinate $1$ and all others $0$
\item[${\mathbf F}$\hfill] Fisher information matrix of the design, ${\mathbf F} = \sum_{i=1}^m n_i {\mathbf F}_i$
\item[$f$\hfill] $f({\mathbf w}) = f(w_1, \ldots, w_m) = |{\mathbf G}^T {\mathbf W} {\mathbf G}|$ which is proportional to $|{\mathbf F}|$; or $f({\mathbf n}) = f(n_1, \ldots, n_m) = |\sum_{i=1}^m n_i {\mathbf F}_i| = |{\mathbf F}|$
\item[${\mathbf F}_i$\hfill] Fisher information matrix at the $i$th design point
\item[$f_i$\hfill] $f_i(z) = f(w_1(1-z)/(1-w_i), \ldots, w_{i-1}(1-z)/(1-w_i), z,$ $w_{i+1}(1-z)/(1-w_i), \ldots, w_m(1-z)/(1-w_i))$ with $0 \leq z < 1$
\item[$f_{ij}$\hfill] $f_{ij}(z) = f(n_1, \ldots, n_{i-1}, z, n_{i+1}, \ldots, n_{j-1}, n_i+n_j-z, n_{j+1},$ $\ldots,$ $n_m)$ with $z= 0, 1, \ldots, n_i+n_j$
\item[${\mathbf G}$\hfill] Matrix component for Fisher information matrix such that ${\mathbf F}=n{\mathbf G}^T {\mathbf W} {\mathbf G}$, $mJ\times p$
\item[$g_s$\hfill] $g_0=f_{ij}(0)$ and $(g_1, \ldots, g_q)^T={\mathbf B}_q^{-1}(d_1,$ $\ldots,$ $d_q)^T$
\item[${\mathbf H}$\hfill] Matrix component for Fisher information matrix such that ${\mathbf F} = {\mathbf H} {\mathbf U} {\mathbf H}^T$, consisting of ${\mathbf H}_1, \ldots, {\mathbf H}_{J-1}$ and possibly ${\mathbf H_c}$, $p\times m(J-1)$
\item[${\mathbf H}_c$\hfill] Matrix for the common component of $J-1$ categories, $({\mathbf h}_c({\mathbf x}_1), \ldots, {\mathbf h}_c({\mathbf x}_m))$, $p_c\times m$
\item[${\mathbf h}_c({\mathbf x}_i)$\hfill] Vector of $p_c$ predictors associated with the $p_c$ parameters ${\boldsymbol\zeta} = (\zeta_1,$ $\ldots,$ $\zeta_{p_c})^T$ that are common for all of the response categories as known functions of the $i$th experimental setting, $(h_1({\mathbf x}_i), \ldots,$ $h_{p_c}({\mathbf x}_i))^T$
\item[${\mathbf H}_j$\hfill] Matrix for the $j$th category only, $({\mathbf h}_j({\mathbf x}_1), \ldots, {\mathbf h}_j({\mathbf x}_m))$, $p_j\times m$
\item[${\mathbf h}_j({\mathbf x}_i)$\hfill] Vector of $p_j$ predictors associated with the $p_j$ parameters $\boldsymbol\beta_j = (\beta_{j1}, \ldots, \beta_{jp_j})^T$ for the $j$th response category as known functions of the $i$th experimental setting, $(h_{j1}({\mathbf x}_i),$ $\ldots,$ $h_{jp_j}({\mathbf x}_i))^T$
\item[${\mathbf I}_k$\hfill] The identity matrix of order $k$
\item[$J$\hfill] Total number of response categories
\item[$k_{\min}$\hfill] Smallest possible $\#\{i\mid \alpha_i > 0\}$ such that $c_{\alpha_1, \ldots, \alpha_m} > 0$
\item[${\mathbf L}$\hfill] Constant $(2J-1)\times J$ matrix, different for the four logit models
\item[$m$\hfill] Total number of distinct experimental settings or design points
\item[${\cal M}({\mathbf H})$\hfill] Column space of matrix ${\mathbf H}$, that is, the linear subspace spanned by the columns of ${\mathbf H}$
\item[$n$\hfill] Total number of experimental units, $n=n_1 + \cdots + n_m$
\item[${\mathbf n}$\hfill] Allocation of experimental units, $(n_1, \ldots, n_m)^T$, $n_i\geq 0$, $\sum_i n_i =n$
\item[$n_i$\hfill] Number of replicates at the $i$th experimental setting
\item[$p$\hfill] Total number of parameters
\item[$p_c$\hfill] Number of common parameters for $J-1$ categories
\item[$p_H$\hfill] $dim \left(\cap_{j=1}^{J-1} {\cal M}({\mathbf H}_j^T )\right)$
\item[$p_j$\hfill] Number of parameters for the $j$th category only
\item[$q$\hfill] $\min\{2J-2, p-k_{\min}+2, p\}$, upper bound of order of $f_{ij}(z)$
\item[$S$\hfill] Collection of all feasible approximate allocations, $\{(w_1, \ldots,$ $w_m)^T \in \mathbb{R}^m \mid w_i \geq 0, i=1, \ldots, m; \sum_{i=1}^m w_1 = 1\}$
\item[$S_+$\hfill] Collection of approximate allocations, $\{{\mathbf w} \in S\mid f({\mathbf w}) > 0\}$
\item[${\mathbf U}$\hfill] Block matrix $({\mathbf U}_{st})_{s,t=1,\ldots, J-1}$, $m(J-1)\times m(J-1)$
\item[${\mathbf U}_{st}$\hfill] ${\rm diag}\{n_1 u_{st}({\boldsymbol\pi}_1), \ldots, n_m u_{st}({\boldsymbol\pi}_m)\}$, $m\times m$
\item[$u_{st}(\boldsymbol\pi_i)$\hfill] ${\mathbf c}_{is}^T{\rm diag}({\boldsymbol\pi}_i)^{-1} {\mathbf c}_{it}$ for $s,t = 1, \ldots, J-1$
\item[${\mathbf w}$\hfill] Real-valued allocation of experimental units, $(w_1, \ldots, w_m)^T$, $w_i \geq 0$, $\sum_i w_i = 1$
\item[${\mathbf W}$\hfill] ${\rm diag}\{w_1 {\rm diag}({\boldsymbol\pi}_1)^{-1}, \ldots, w_m {\rm diag}({\boldsymbol\pi}_m)^{-1}\}$, $mJ\times mJ$
\item[$w_i$\hfill] Proportion of experimental units assigned to the $i$th experimental setting, $n_i/n$
\item[${\mathbf w}_u$\hfill] Uniform allocation, $(1/m, \ldots, 1/m)^T$
\item[${\cal X}$\hfill] Design space, the collection of all design points yielding strictly positive categorical probabilities of response; or a predetermined set of design points considered
\item[${\mathbf x}_i$\hfill] The $i$th distinct experimental setting or design point, $(x_{i1},$ $\ldots,$ $x_{id})^T$
\item[${\mathbf X}_i$\hfill] Model matrix at the $i$th design point, $J\times p$, the last row is all $0$'s
\item[$\boldsymbol\beta_{j}$\hfill] Vector of parameters for the $j$th response category only, $(\beta_{j1}, \ldots, \beta_{jp_j})^T$
\item[$\gamma_{ij}$\hfill] The cumulative probability from the $1$st to $j$th categories at the $i$th experimental setting, $\gamma_{ij}=\pi_{i1}+\cdots + \pi_{ij}$
\item[$\boldsymbol\zeta$\hfill] Vector of common parameters for all of the response categories, $(\zeta_1, \ldots, \zeta_{p_c})^T$
\item[$\boldsymbol\eta_{i}$\hfill] Vector of linear predictors at the $i$th experimental setting, ${\boldsymbol\eta}_i = (\eta_{i1}, \ldots, \eta_{iJ})^T = {\mathbf X}_i \boldsymbol\theta$ with $\eta_{iJ}\equiv 0$
\item[$\boldsymbol\theta$\hfill] Vector of all parameters, $p\times 1$
\item[$\boldsymbol\Theta$\hfill] Parameter space, the collection of all feasible parameter vectors
\item[$\boldsymbol\pi_i$\hfill] Vector of response category probabilities at the $i$th experimental setting. $\boldsymbol\pi_i = (\pi_{i1}, \ldots, \pi_{iJ})^T$, $\pi_{i1} + \cdots + \pi_{iJ}=1$
\item[$\Lambda (\alpha_1, \ldots, \alpha_m)$\hfill]  $\{(i_1, \ldots, i_p) \mid 1\leq i_1 < \cdots < i_p\leq mJ;
  \# \{l:(k-1)J<i_l \leqslant kJ \}=\alpha_k, k=1, \ldots, m  \}$
\item[$\pi_{ij}$\hfill] Probability that the response falls into the $j$th category at the $i$th experimental setting
\item[$\phi$\hfill] Bayesian D-optimal criterion, $\phi({\mathbf p}) = E(\log|{\mathbf F}|)$
\end{list}

\subsection{Formulae of matrix differentiation}\label{matrixformula}

According to \citeauthor{seber2008} (2008, Chapter~17)),
\begin{eqnarray*}
\frac{\partial {\mathbf y}}{\partial {\mathbf x}^T} &=& \left(\frac{\partial y_i}{\partial x_j}\right)_{ij}\\
\frac{\partial {\mathbf A}{\mathbf x}}{\partial {\mathbf x}^T} &=& {\mathbf A}\\
\frac{\partial {\mathbf z}}{\partial {\mathbf x}^T} &=& \frac{\partial {\mathbf z}}{\partial {\mathbf y}^T}\cdot\frac{\partial {\mathbf y}}{\partial {\mathbf x}^T}\\
\frac{\partial \log {\mathbf y}}{\partial {\mathbf x}^T} &=& [{\rm diag}({\mathbf y})]^{-1}\frac{\partial {\mathbf y}}{\partial {\mathbf x}^T}
\end{eqnarray*}
where ${\mathbf x}=(x_i)_i,\ {\mathbf y}=(y_i)_i,\ {\mathbf z} = (z_i)_i$, and thus $\log{\mathbf y} = (\log y_i)_i$ are vectors, and ${\mathbf A}$ is a constant matrix.

\subsection{Explicit forms of $({\mathbf C}^T{\mathbf D}_i^{-1}{\mathbf L})^{-1}$ for all the four logit models}\label{subsec:cij}

There are the four different kinds of multinomial logistic models in the literature: {\it baseline-category logit model} for nominal responses, {\it cumulative logit model} for ordinal responses, {\it adjacent-categories logit model} for ordinal responses, and {\it continuation-ratio logit model} for hierarchical responses. According to Theorem~\ref{thm:fisher}, $({\mathbf C}^T{\mathbf D}_i^{-1}{\mathbf L})^{-1}$ is a key matrix that we must calculate.

Recall that $\pi_{i1} + \cdots + \pi_{iJ} = 1$, $i=1, \ldots, m$. Then
\begin{equation*}
({\mathbf C}^T{\mathbf D}_i^{-1}{\mathbf L})_{baseline}=\begin{pmatrix}
\textstyle\frac{1}{\pi_{i1}}& 0 & \cdots & 0 & -\textstyle\frac{1}{\pi_{iJ}}\\
0 & \textstyle\frac{1}{\pi_{i2}}& \ddots & \vdots & -\textstyle\frac{1}{\pi_{iJ}}\\
\vdots & \ddots & \ddots & 0 & \vdots\\
0 & \cdots & 0 & \textstyle\frac{1}{\pi_{i,J-1}} & -\textstyle\frac{1}{\pi_{iJ}}\\
1 & 1 & \cdots & 1 & 1\\
\end{pmatrix}_{J \times J}
\end{equation*}
\begin{equation*}
({\mathbf C}^T{\mathbf D}_i^{-1}{\mathbf L})_{cumulative}=\begin{pmatrix}
\frac{1}{\gamma_{i1}}& -\frac{1}{1-\gamma_{i1}} & -\frac{1}{1-\gamma_{i1}} & \cdots & -\frac{1}{1-\gamma_{i1}}\\
\frac{1}{\gamma_{i2}} & \frac{1}{\gamma_{i2}} & -\frac{1}{1-\gamma_{i2}} & \cdots & -\frac{1}{1-\gamma_{i2}}\\
\vdots & \vdots & \ddots & \ddots & \vdots\\
\frac{1}{\gamma_{i,J-1}} & \frac{1}{\gamma_{i,J-1}} & \cdots & \frac{1}{\gamma_{i,J-1}}  & -\frac{1}{1-\gamma_{i,J-1}}\\
1 & 1 & \cdots & 1 & 1\\
\end{pmatrix}_{J \times J}
\end{equation*}
\begin{equation*}
({\mathbf C}^T{\mathbf D}_i^{-1}{\mathbf L})_{continuation}=\begin{pmatrix}
\frac{1}{\pi_{i1}}& -\frac{1}{1-\gamma_{i1}} & -\frac{1}{1-\gamma_{i1}} & \cdots & -\frac{1}{1-\gamma_{i1}}\\
0 & \frac{1}{\pi_{i2}} & -\frac{1}{1-\gamma_{i2}} & \cdots & -\frac{1}{1-\gamma_{i2}}\\
\vdots & \ddots & \ddots & \ddots & \vdots\\
0 & \cdots & 0 & \frac{1}{\pi_{i,J-1}}  & -\frac{1}{1-\gamma_{i,J-1}}\\
1 & 1 & \cdots & 1 & 1\\
\end{pmatrix}_{J \times J}
\end{equation*}
\begin{equation*}
({\mathbf C}^T{\mathbf D}_i^{-1}{\mathbf L})_{adjacent}=\begin{pmatrix}
\textstyle\frac{1}{\pi_{i1}}& -\textstyle\frac{1}{\pi_{i2}} & 0 & \cdots & 0 \\
0 & \textstyle\frac{1}{\pi_{i2}}& -\textstyle\frac{1}{\pi_{i3}} & \ddots & \vdots\\
\vdots & \ddots & \ddots & \ddots & 0\\
0 & \cdots & 0 & \textstyle\frac{1}{\pi_{i,J-1}} & -\textstyle\frac{1}{\pi_{iJ}}\\
1 & 1 & \cdots & 1 & 1\\
\end{pmatrix}_{J \times J}
\end{equation*}
where
$\gamma_{ij} = \pi_{i1}+\cdots +\pi_{ij}$ is the cumulative categorical probability, $j=1, \ldots, J-1$.
The corresponding inverse matrices are
\begin{align*}
& ({\mathbf C}^T{\mathbf D}_i^{-1}{\mathbf L})^{-1}_{baseline}\\
&=\begin{pmatrix}
-\pi_{i1}^2+\pi_{i1} & -\pi_{i1}\pi_{i2} &  \cdots & -\pi_{i1}\pi_{i,J-1} & \pi_{i1}\\
-\pi_{i1}\pi_{i2} & -\pi_{i2}^2+\pi_{i2} & \cdots & -\pi_{i2}\pi_{i,J-1} & \pi_{i2}\\
\vdots & \vdots & \ddots & \vdots & \vdots\\
-\pi_{i1}\pi_{i,J-1} & -\pi_{i2}\pi_{i,J-1} & \cdots & -\pi_{i,J-1}^2+\pi_{i,J-1} & \pi_{i,J-1}\\
-\pi_{i1}\pi_{iJ} & -\pi_{i2}\pi_{iJ} & \cdots & -\pi_{i,J-1}\pi_{iJ} & \pi_{iJ}\\
\end{pmatrix}_{J \times J}\\
&\stackrel{\triangle}{=}\begin{pmatrix}
\mathbf c_{i1} & \mathbf c_{i2} & \cdots & \mathbf c_{iJ}
\end{pmatrix}_{baseline}
\end{align*}
where
$({\mathbf c}_{ij})_{baseline} = \pi_{ij}({\mathbf e}_j - {\boldsymbol\pi}_i)$, $j=1, \ldots, J-1$, $({\mathbf c}_{iJ})_{baseline} = {\boldsymbol\pi}_i$, and ${\mathbf e}_j$ here is the $J\times 1$ vector with the $j$th coordinate $1$ and all others $0$. Recall that ${\boldsymbol\pi}_i = (\pi_{i1}, \ldots, \pi_{iJ})^T$.
\begin{align*}
& ({\mathbf C}^T{\mathbf D}_i^{-1}{\mathbf L})^{-1}_{cumulative}\\
&=\begin{pmatrix}
\gamma_{i1}(1-\gamma_{i1}) & 0 & \cdots & 0 & \pi_{i1}\\
-\gamma_{i1}(1-\gamma_{i1}) & \gamma_{i2}(1-\gamma_{i2}) & \ddots & \vdots & \pi_{i2}\\
0 & -\gamma_{i2}(1-\gamma_{i2}) & \ddots & 0 & \vdots\\
\vdots & \ddots & \ddots & \gamma_{i,J-1}(1-\gamma_{i,J-1}) & \pi_{i,J-1}\\
0 & \cdots & 0 & -\gamma_{i,J-1}(1-\gamma_{i,J-1}) & \pi_{iJ}\\
\end{pmatrix}_{J \times J}\\
&\stackrel{\triangle}{=}\begin{pmatrix}
\mathbf c_{i1} & \mathbf c_{i2} & \cdots & \mathbf c_{iJ}
\end{pmatrix}_{cumulative}
\end{align*}
where
$(\mathbf c_{ij})_{cumulative}=\gamma_{ij}(1-\gamma_{ij})({\mathbf e}_j-{\mathbf e}_{j+1})$ with ${\mathbf e}_j$ defined as above;
and $(\mathbf c_{iJ})_{cumulative} = {\boldsymbol\pi}_i$~.
\begin{align*}
& ({\mathbf C}^T{\mathbf D}_i^{-1}{\mathbf L})^{-1}_{continuation}\\ &=\begin{pmatrix}
\pi_{i1}(1-\gamma_{i1}) & 0 & \cdots & 0 & \pi_{i1}\\
-\pi_{i1}\pi_{i2} & \frac{\pi_{i2}(1-\gamma_{i2})}{1-\gamma_{i1}} & \ddots & \vdots & \pi_{i2}\\
\vdots & \vdots & \ddots & 0 & \vdots\\
-\pi_{i1}\pi_{i,J-1} & -\frac{\pi_{i2}\pi_{i,J-1}}{1-\gamma_{i1}} & \cdots & \frac{\pi_{i,J-1}(1-\gamma_{i,J-1})}{1-\gamma_{i,J-2}} & \pi_{i,J-1}\\
-\pi_{i1}\pi_{iJ} & -\frac{\pi_{i2}\pi_{iJ}}{1-\gamma_{i1}} & \cdots & -\frac{\pi_{i,J-1}\pi_{iJ}}{1-\gamma_{i,J-2}}  & \pi_{iJ}\\
\end{pmatrix}_{J \times J}\\
&=\begin{pmatrix}
\mathbf c_{i1} & \mathbf c_{i2} & \cdots & \mathbf c_{iJ}
\end{pmatrix}_{continuation}
\end{align*}
where
$({\mathbf c}_{i1})_{continuation} = \pi_{i1}(1-\gamma_{i1}, -\pi_{i2}, \ldots, -\pi_{iJ})^T$,\\ $({\mathbf c}_{ij})_{continuation}$ $=$ $\frac{\pi_{ij}}{1-\gamma_{i,j-1}}(0,$ $\ldots,$ $0, 1-\gamma_{ij}, -\pi_{i,j+1}, \ldots, -\pi_{iJ})^T$ with ``$1-\gamma_{ij}$" being the $j$th coordinate, $j=2, \ldots, J-1$,
and $(\mathbf c_{iJ})_{continuation} = {\boldsymbol\pi}_i$~.
\begin{align*}
& ({\mathbf C}^T{\mathbf D}_i^{-1}{\mathbf L})^{-1}_{adjacent}\\ &=\begin{pmatrix}
(1-\gamma_{i1}) \pi_{i1} & (1-\gamma_{i2})\pi_{i1} & \cdots & (1-\gamma_{i,J-1})\pi_{i1} & \pi_{i1}\\
-\gamma_{i1}\pi_{i2} & (1-\gamma_{i2})\pi_{i2} & \cdots & (1-\gamma_{i,J-1})\pi_{i2} & \pi_{i2}\\
\vdots & \vdots & \ddots & \vdots & \vdots\\
-\gamma_{i1}\pi_{i,J-1} & -\gamma_{i2}\pi_{i,J-1} & \cdots & (1-\gamma_{i,J-1})\pi_{i,J-1} & \pi_{i,J-1}\\
-\gamma_{i1}\pi_{iJ} & -\gamma_{i2}\pi_{iJ} & \cdots & -\gamma_{i,J-1}\pi_{iJ} & \pi_{iJ}\\
\end{pmatrix}_{J \times J}\\
&=\begin{pmatrix}
\mathbf c_{i1} & \mathbf c_{i2} & \cdots & \mathbf c_{iJ}
\end{pmatrix}_{adjacent}
\end{align*}
where
$({\mathbf c}_{ij})_{adjacent}=((1-\gamma_{ij})\pi_{i1}, \ldots, (1-\gamma_{ij})\pi_{ij}, -\gamma_{ij}\pi_{i,j+1}, \ldots, -\gamma_{ij}\pi_{iJ})^T$, $j=1, \ldots, J-1$,
and $(\mathbf c_{iJ})_{adjacent} = {\boldsymbol\pi}_i$~.

For certain applications, we need to know $|{\mathbf C}^T{\mathbf D}_i^{-1}{\mathbf L}|$ (see, for example, Lemma~\ref{lemma:vi=CDL}). Since adding a multiple of one row
(column) to another row (column) does not change the determinant (see, for example, 4.28(f) in \citeauthor{seber2008} (2008, page 58)), we may (1) do row operations on $({\mathbf C}^T{\mathbf D}_i^{-1}{\mathbf L})_{baseline}$ and change it into an upper triangular matrix with diagonal entries $\pi_{i1}^{-1}, \ldots, \pi_{iJ}^{-1}$; (2) do row operations on $({\mathbf C}^T{\mathbf D}_i^{-1}{\mathbf L})^{-1}_{cumulative}$ and change it into an upper triangular matrix with diagonal entries $\gamma_{i1} (1-\gamma_{i1}), \ldots, \gamma_{i,J-1}(1-\gamma_{i,J-1}), 1$; (3) do column operations on $({\mathbf C}^T{\mathbf D}_i^{-1}{\mathbf L})_{adjacent}$ and change it into a lower triangular matrix with diagonal entries $\pi_{i1}^{-1}, \ldots, \pi_{iJ}^{-1}$; and (4) do column operations on $({\mathbf C}^T{\mathbf D}_i^{-1}{\mathbf L})_{continuation}$ and change it into a lower triangular matrix with diagonal entries $\pi_{i1}^{-1}, \ldots, \pi_{iJ}^{-1}$. Therefore,
\begin{equation}\label{eq:CDiLdet}
|{\mathbf C}^T{\mathbf D}_i^{-1}{\mathbf L}| =
\left\{\begin{array}{ll}
\prod_{j=1}^J \pi_{ij}^{-1} & \mbox{for baseline-category,}\\
& \mbox{ adjacent-categories,}\\
& \mbox{ and continuation-ratio logit models}\\
\prod_{j=1}^{J-1} \gamma_{ij}^{-1} (1-\gamma_{ij})^{-1} & \mbox{for cumulative logit models}
\end{array}\right.
\end{equation}
As a direct conclusion, $|{\mathbf C}^T{\mathbf D}_i^{-1}{\mathbf L}|> 0$ as long as $\pi_{ij}>0$ for all $j=1, \ldots, J$.

\subsection{Positive definiteness of ${\mathbf U}$}\label{subsec:positiveU}

In order to determine the positive definiteness of ${\mathbf F}$, we first investigate the $m(J-1)\times m(J-1)$ matrix ${\mathbf U}$ defined for Theorem \ref{thm:fisher_all}, which is symmetric since $u_{st}(\boldsymbol\pi_i) = u_{ts}(\boldsymbol\pi_i)$ and thus ${\mathbf U}_{st} = {\mathbf U}_{ts}$.

\begin{theorem}\label{thm:Upd}
	If $n_i > 0$ for all $i=1, \ldots, m$, then ${\mathbf U}$ is positive definite.
\end{theorem}

\begin{theorem}\label{thm:vdet}
	$|{\mathbf U}| = (\prod_{i=1}^{m} n_i)^{J-1}\cdot \prod_{i=1}^m (\prod_{j=1}^{J} \pi_{ij})^{-1} |{\mathbf C}^T {\mathbf D}_i^{-1} {\mathbf L}|^{-2}~.
	$
\end{theorem}

The proofs of Theorems~\ref{thm:Upd} and \ref{thm:vdet} are relegated to Section~\ref{sec:proofs}.
Note that Theorem~\ref{thm:Upd} is not a corollary of Theorem~\ref{thm:vdet} since nonsingularity itself does not mean positive definiteness.
Theorem~\ref{thm:vdet} implies that ${\mathbf U}$ is singular if $n_i=0$ for some $i=1, \ldots, m$.
Note that ${\mathbf F}$ can still be positive definite even if ${\mathbf U}$ is singular, as long as ${\mathbf H}$ is of full row rank.
In general, given an allocation $(n_1, \ldots, n_m)$ of the $n$ experimental units with $n_i \geq 0$ and $\sum_{i=1}^m n_i = n$, if we denote $k = \#\{i:n_i>0\}$ and ${\mathbf U}_{st}^* = {\rm diag}\{n_i u_{st}({\boldsymbol\pi}_i) : n_i>0\}$, then ${\mathbf U}^* = ({\mathbf U}^*_{st})_{s,t=1, \ldots, J-1}$
is a $k(J-1)\times k(J-1)$ matrix. After removing all columns of ${\mathbf H}$ associated with $n_i=0$, we denote the leftover as ${\mathbf H}^*$, which is a $p\times k(J-1)$ matrix. It can be verified that

\begin{lemma}\label{lemma:F=>F*}
	${\mathbf H} {\mathbf U} {\mathbf H}^T = {\mathbf H}^* {\mathbf U}^* \left({\mathbf H}^*\right)^T$.
\end{lemma}

\begin{lemma}\label{lemma:U*det}
	$|{\mathbf U}^*| = (\prod_{i:n_i>0} n_i)^{J-1}\cdot \prod_{i:n_i>0} (\prod_{j=1}^{J} \pi_{ij})^{-1} |{\mathbf C}^T {\mathbf D}_i^{-1} {\mathbf L}|^{-2}
	$.
\end{lemma}

Since ${\mathbf U}^*$ is simply ${\mathbf U}$ if all $n_i >0$, we have the following corollary of Theorem~\ref{thm:Upd}:

\begin{corollary}\label{col:U*pd}
	${\mathbf U}^*$ is positive definite.
\end{corollary}

\subsection{Row rank of ${\mathbf H}$ matrix}\label{subsec:rowrankH}

According to Theorem~\ref{thm:Fpd}, the positive definiteness of the Fisher information matrix ${\mathbf F}$ depends on the row rank of ${\mathbf H}$  or ${\mathbf H}^*$.
To simplify the notations, we assume $n_i > 0, i=1, \ldots, m$ throughout this section. In this case, ${\mathbf H}= {\mathbf H}^*$ and ${\mathbf U} = {\mathbf U}^*$. We also assume that
\begin{equation}\label{rowassumption}
m\geq p_j,\quad j=1, \ldots, J-1\quad\mbox{ and }\quad m\geq p_c\>\mbox{ if applicable}
\end{equation}
since ${\mathbf H}$ is of full row rank only if
$rank({\mathbf H}_j) = p_j$, $j=1, \ldots, J-1$ and $rank({\mathbf H}_c) = p_c$ if applicable.

Since ${\mathbf H}$ takes different forms for {\it ppo}, {\it npo}, and {\it po} models, we investigate its row rank case by case.

\begin{theorem}\label{thm:rankH}
	Consider the $p\times m(J-1)$ matrix ${\mathbf H}$ in Theorem~\ref{thm:fisher_all}.
	\begin{itemize}
		\item[(1)] For {\it npo} models,
		$rank({\mathbf H}) = rank({\mathbf H}_1) + \cdots + rank({\mathbf H}_{J-1})$.
		\item[(2)] For {\it po} models,
		$rank({\mathbf H}) = rank(\left({\mathbf 1}, {\mathbf H}_c^T\right)) + J-2$, where ${\mathbf 1}$ is a vector of all $1$'s.
		\item[(3)] For {\it ppo} models,
		$rank({\mathbf H})=rank({\mathbf H_{1}})+\cdots+rank({\mathbf H_{J-1}}) + rank({\mathbf H_{c}})-dim [{\cal M}(\mathbf H_c^T)\cap ( \cap_{j=1}^{J-1}{\cal M}(\mathbf H_j^T))]$, where ${\cal M}(\mathbf H_c^T)$ stands for the column space of $\mathbf H_c^T$ or the row space of $\mathbf H_c$.
	\end{itemize}
\end{theorem}

The proof of Theorem~\ref{thm:rankH} is relegated to Section~\ref{sec:proofs}. In order to apply it to {\it ppo} models, we need an efficient way to calculate $dim [{\cal M}(\mathbf H_c^T)\cap ( \cap_{j=1}^{J-1}{\cal M}(\mathbf H_j^T))]$. We provide a formula for calculating $dim(\bigcap_{j} {\cal M}({\mathbf H}_j^T))$ for general matrices, Theorem~\ref{thm:rankintersection} in the Appendix, and relegated its proof to Section~\ref{sec:proofs}.

Recall that $p_H$ $=$ $dim (\cap_{j=1}^{J-1} {\cal M}({\mathbf H}_j^T ))$.
As a direct conclusion of Theorem~\ref{thm:rankH}, we have

\begin{corollary}\label{lemma:ppo_m}
	For {\it ppo} models, $|F|>0$ only if $m \geqslant p_c + p_H$.
\end{corollary}

\subsection{Results on the coefficient $c_{\alpha_1, \ldots, \alpha_m}$ for simplifying $|{\mathbf F}|$}\label{subsec:alpha}

\begin{lemma}\label{lemma:max>=J}
	If $\max_{1\leq i\leq m} \alpha_i \geq J$, then $|{\mathbf G}[i_1, \ldots, i_p]| = 0$ for any $(i_1, \ldots, $ $i_p)$ $\in$ $\Lambda (\alpha_1,$ $\ldots,$ $\alpha_m)$. Therefore, $c_{\alpha_1, \ldots, \alpha_m}=0$ in this case.
\end{lemma}

\begin{theorem}\label{thm:alpha>0}
	The coefficient $c_{\alpha_1, \ldots, \alpha_m}$ as defined in (\ref{eq:calpha}) is nonzero only if the restricted Fisher information matrix ${\mathbf F}_{res}=\sum_{i : \alpha_i>0} {\mathbf F}_i$ is positive definite, where ${\mathbf F}_i$ is defined as in (\ref{eqn:Fi}).
\end{theorem}

The proofs for Lemma~\ref{lemma:max>=J} and Theorem~\ref{thm:alpha>0} are relegated to Section~\ref{sec:proofs}. Combining Theorems~\ref{thm:Fpd} and \ref{thm:alpha>0}, Theorems~\ref{thm:minimalm} and \ref{thm:alpha>0}, respectively, we obtain the following corollaries:

\begin{corollary}\label{col:calpha>0H*}
	The coefficient $c_{\alpha_1, \ldots, \alpha_m}$ is nonzero only if ${\mathbf H}_{\alpha_1, \ldots, \alpha_m}$ is of full row rank $p$, where ${\mathbf H}_{\alpha_1, \ldots, \alpha_m}$ is the submatrix of ${\mathbf H}$ after removing all columns associated with ${\mathbf x}_i$ for which $\alpha_i=0$.
\end{corollary}

\begin{corollary}\label{col:calpha>0minimalm}
	The coefficient $c_{\alpha_1, \ldots, \alpha_m} = 0$ if $\#\{i\mid \alpha_i>0\} \leq k_{\min} - 1$, where
	$k_{\min} = \max\{p_1, \ldots, p_{J-1}, p_c + p_H\}$.  If ${\mathbf H}_1 = \cdots = {\mathbf H}_{J-1}$, $k_{\min} = p_c + p_1$.
\end{corollary}

We provide an example (Example~\ref{example:ppopc}) in Section~\ref{subsec:moreexamples} to illustrate that $c_{\alpha_1, \ldots, \alpha_m}$ could be nonzero for {\it ppo} models with $\#\{i\mid \alpha_i > 0\} = p_c + p_H$~.

\subsection{Expressions for proportional odds ({\it po}) models}\label{subsec:pomodels}

As special cases of {\it ppo}, {\it po} models are degenerate cases of {\it ppo} models with ${\mathbf h}_j^T({\mathbf x}_i)$ replaced by $1$, $j=1, \ldots, J-1$, and thus $p_1  = \cdots = p_{J-1} = 1$.

In Section~\ref{sec:fisher}, the four logit models in the literature with proportional odds are:
\begin{eqnarray*}
	\log\left(\frac{\pi_{ij}}{\pi_{iJ}}\right) &=& \beta_j+{\mathbf h}_c^T({\mathbf x}_i)\boldsymbol\zeta~, \>\mbox{baseline-category}\\
	\log\left(\frac{\pi_{i1}+\cdots + \pi_{ij}}{\pi_{i,j+1}+\cdots + \pi_{iJ}}\right) &=& \beta_j+{\mathbf h}_c^T({\mathbf x}_i)\boldsymbol\zeta~, \>\mbox{cumulative}\\
	\log\left(\frac{\pi_{ij}}{\pi_{i,j+1}}\right) &=& \beta_j+{\mathbf h}_c^T({\mathbf x}_i)\boldsymbol\zeta~, \>\mbox{adjacent-categories}\\
	\log\left(\frac{\pi_{ij}}{\pi_{i,j+1} + \cdots + \pi_{iJ}}\right) &=& \beta_j+{\mathbf h}_c^T({\mathbf x}_i)\boldsymbol\zeta~, \>\mbox{continuation-ratio}
\end{eqnarray*}
where $i=1, \ldots, m$, $j=1, \ldots, J-1$, $\beta_j$ is an unknown parameter for the $j$th response category, ${\mathbf h}_c^T(\cdot) = (h_{1}(\cdot), \ldots, h_{p_c}(\cdot))$ are known functions to determine the $p_c$ predictors associated with the $p_c$ unknown parameters $\boldsymbol\zeta = (\zeta_1, \ldots, \zeta_{p_c})^T$ that are common for all categories. 

In equation~\eqref{logitunifiedmodel}, the corresponding model matrix is
\begin{equation}\label{eqn:Xi_po}
{\mathbf X}_i= \begin{pmatrix}
1  & 0  & \cdots & 0 & {\mathbf h}_c^T({\mathbf x}_i)\\
0  & 1  &\ddots  & \vdots & \vdots\\
\vdots  & \ddots & \ddots & 0 & {\mathbf h}_c^T({\mathbf x}_i)\\
0  &  \cdots & 0 & 1 & {\mathbf h}_c^T({\mathbf x}_i)\\
0  & 0  & \cdots & 0 & \boldsymbol0^T\\
\end{pmatrix}_{J \times p}
\end{equation}
and the parameter vector $\boldsymbol\theta=(\beta_1,\beta_2,\cdots,\beta_{J-1},\boldsymbol\zeta)^T$ consists of $p=J-1+p_c$ unknown parameters in total. The previous $\boldsymbol\beta_{j}$ reduces to $\beta_j$ serving as the cut-off point in this case.

In Section~\ref{sec:positive}, the $p\times m(J-1)$ matrix
\begin{equation}
{\mathbf H}=
\left(\begin{array}{ccc}
{\mathbf 1}^T & &\\  & \ddots & \\& &{\mathbf 1}^T \\{\mathbf H}_{c} & \cdots & {\mathbf H}_{c}
\end{array}\right)
\end{equation}
where ${\mathbf H}_{c} = ({\mathbf h}_c({\mathbf x}_1), \cdots, {\mathbf h}_c({\mathbf x}_m))$.

As a special case of Theorem~\ref{thm:minimalm}, 

\begin{theorem}\label{thm:minimalm_po}
	Consider the multinomial logistic model (\ref{logitunifiedmodel}) with $m$ distinct experimental settings ${\mathbf x}_i$ with $n_i>0$ experimental units, $i=1, \ldots, m$. For proportional odds models, the Fisher information matrix ${\mathbf F}$ is positive definite if and only if $m\geq p_c+1$ and the extended matrix $({\mathbf 1}, {\mathbf H}_c^T)$ is of full rank $p_c+1$.
\end{theorem}

In Section~\ref{sec:determinant}, for proportional odds models, the $mJ\times p$ matrix
\begin{equation}
{\mathbf G} = \left(\begin{array}{cccc}
{\mathbf c}_{11} & \cdots & {\mathbf c}_{1,J-1} & \sum_{j=1}^{J-1} {\mathbf c}_{1j} \cdot {\mathbf h}_c^T({\mathbf x}_1)\\
{\mathbf c}_{21} & \cdots & {\mathbf c}_{2,J-1} & \sum_{j=1}^{J-1} {\mathbf c}_{2j} \cdot {\mathbf h}_c^T({\mathbf x}_2)\\
\cdots & \cdots & \cdots & \cdots \\
{\mathbf c}_{m1} & \cdots & {\mathbf c}_{m,J-1} & \sum_{j=1}^{J-1} {\mathbf c}_{mj} \cdot {\mathbf h}_c^T({\mathbf x}_m)
\end{array}\right)
\end{equation}

As a special case of Corollary~\ref{col:calpha>0minimalm},

\begin{corollary}\label{col:calpha>0minimalm_po}
	The coefficient $c_{\alpha_1, \ldots, \alpha_m} = 0$ if $\#\{i\mid \alpha_i>0\} \leq k_{\min} - 1$, where
	$k_{\min} = p_c+1$ for {\it po} models.
\end{corollary}

As special cases of {\it ppo} models, {\it po} models imply $p_1 = \cdots = p_{J-1} = p_H =1$, and ${\mathbf H}_1 = \cdots = {\mathbf H}_{J-1}$ implies $p_1 = \cdots = p_{J-1} = p_H$. That is, $k_{\min}$'s are consistent across different odds models.

\subsection{Expressions for non-proportional odds ({\it npo}) models}\label{subsec:npomodels}

As special cases of {\it ppo}, ${\mathbf h}_c^T({\mathbf x}_i) \equiv 0$ leads to {\it npo} models. Therefore, $p_c=0$.

In Section~\ref{sec:fisher}, the four logit models in the literature with non-proportional odds are:
\begin{eqnarray*}
	\log\left(\frac{\pi_{ij}}{\pi_{iJ}}\right) &=& {\mathbf h}_j^T({\mathbf x}_i)\boldsymbol\beta_j~, \>\mbox{baseline-category}\\
	\log\left(\frac{\pi_{i1}+\cdots + \pi_{ij}}{\pi_{i,j+1}+\cdots + \pi_{iJ}}\right) &=& {\mathbf h}_j^T({\mathbf x}_i)\boldsymbol\beta_j~, \>\mbox{cumulative}\\
	\log\left(\frac{\pi_{ij}}{\pi_{i,j+1}}\right) &=& {\mathbf h}_j^T({\mathbf x}_i)\boldsymbol\beta_j~, \>\mbox{adjacent-categories}\\
	\log\left(\frac{\pi_{ij}}{\pi_{i,j+1} + \cdots + \pi_{iJ}}\right) &=& {\mathbf h}_j^T({\mathbf x}_i)\boldsymbol\beta_j~, \>\mbox{continuation-ratio}
\end{eqnarray*}
where $i=1, \ldots, m$, $j=1, \ldots, J-1$, ${\mathbf h}_j^T(\cdot) = (h_{j1}(\cdot), \ldots, h_{jp_j}(\cdot))$ are known functions to determine the $p_j$ predictors associated with the $p_j$ unknown parameters $\boldsymbol\beta_j = (\beta_{j1}, \ldots, \beta_{jp_j})^T$ for the $j$th response category.

In equation~\eqref{logitunifiedmodel}, the corresponding model matrix is
\begin{equation}\label{eqn:Xi_npo}
{\mathbf X}_i= \begin{pmatrix}
{\mathbf h}_1^T({\mathbf x}_i) &  \boldsymbol0^T & \cdots & \boldsymbol0^T\\
\boldsymbol0^T &  {\mathbf h}_2^T({\mathbf x}_i) &\ddots & \vdots \\
\vdots &  \ddots& \ddots &  \boldsymbol0^T \\
\boldsymbol0^T & \cdots & \boldsymbol0^T & {\mathbf h}_{J-1}^T({\mathbf x}_i) \\
\boldsymbol0^T & \cdots & \cdots & \boldsymbol0^T \\
\end{pmatrix}_{J \times p}
\end{equation}
and the parameter vector reduces to $\boldsymbol\theta=(\boldsymbol\beta_{1},\boldsymbol\beta_{2},\cdots,\boldsymbol\beta_{J-1})^T$,
which consists of $p=p_1+\cdots+p_{J-1}$ unknown parameters in total. Note that we always use $p$ to represent the total number of parameters.

In Section~\ref{sec:positive}, the $p\times m(J-1)$ matrix
\begin{equation}
{\mathbf H}=
\left(\begin{array}{ccc}
{\mathbf H}_{1} & &\\  & \ddots & \\& &{\mathbf H}_{J-1}
\end{array}\right)
\end{equation}
where ${\mathbf H}_{j}= ({\mathbf h}_j({\mathbf x}_1), \cdots, {\mathbf h}_j({\mathbf x}_m))$, $j=1,\ldots,J-1$.

As a special case of Theorem~\ref{thm:minimalm}, we have

\begin{theorem}\label{thm:minimalm_npo}
	Consider the multinomial logistic model (\ref{logitunifiedmodel}) with $m$ distinct experimental settings ${\mathbf x}_i$ with $n_i>0$ experimental units, $i=1, \ldots, m$. For non-proportional odds ({\it npo}) models, the Fisher information matrix ${\mathbf F}$ is positive definite if and only if $m\geq \max\{p_1, \ldots, p_{J-1}\}$ and ${\mathbf x}_i$'s keep ${\mathbf H}_j$ of full row rank $p_j$, $j=1, \ldots, J-1$.
\end{theorem}

In Section~\ref{sec:determinant}, for non-proportional odds models, the $mJ\times p$ matrix
\begin{equation}
{\mathbf G} = \left(\begin{array}{ccc}
{\mathbf c}_{11}{\mathbf h}_1^T({\mathbf x}_1) & \cdots & {\mathbf c}_{1,J-1} {\mathbf h}_{J-1}^T ({\mathbf x}_1)\\
{\mathbf c}_{21}{\mathbf h}_1^T({\mathbf x}_2) & \cdots & {\mathbf c}_{2,J-1} {\mathbf h}_{J-1}^T ({\mathbf x}_2)\\
\cdots & \cdots & \cdots \\
{\mathbf c}_{m1}{\mathbf h}_1^T({\mathbf x}_m) & \cdots & {\mathbf c}_{m,J-1} {\mathbf h}_{J-1}^T ({\mathbf x}_m)
\end{array}\right)
\end{equation}

As a special case of Corollary~\ref{col:calpha>0minimalm}, we have

\begin{corollary}\label{col:calpha>0minimalm_npo}
	The coefficient $c_{\alpha_1, \ldots, \alpha_m} = 0$ if $\#\{i\mid \alpha_i>0\} \leq k_{\min} - 1$, where $k_{\min} = \max\{p_1, \ldots, p_{J-1}\}$ for {\it npo} models.
\end{corollary}

As special cases of {\it ppo} models, {\it npo} models imply $p_c=0$ and $p_H \leq \min\{p_1, \ldots, p_{J-1}\}$. That is, $k_{\min}$'s are consistent across different odds models.

\subsection{Model selection}\label{sec:modelselection}

See Tables~\ref{tab:compare1} and \ref{tab:compare5}.

\begin{table}[ht]
	\caption{Model Comparison for Trauma Clinical Trial Data}
	\centering 
	\footnotesize
	\begin{tabular}{lrrrrrr} 
		\hline\hline  \\ [-1.5ex]
		&Cumulative &Cumulative  &Continuation &Continuation  &Adjacent &Adjacent \\ [-0.5ex]
		&{\it po} &{\it npo} &{\it po} &{\it npo} &{\it po} &{\it npo}\\ [0.5ex]
		\hline\\[-1.5ex]
		AIC &107.75  &{\bf 99.41} &108.98 &101.36 &107.67 &101.54 \\
		BIC &104.68  &{\bf 94.51} &105.91 &96.45 &104.60 &96.63 \\[0.5ex]
		\hline\hline 
	\end{tabular}
	\label{tab:compare1}
	\normalsize
\end{table}

\begin{table}[h]
	\caption{Model Comparison for Emergence of House Flies Data}
	\centering 
	\footnotesize
	\begin{tabular}{lrrrrrr} 
		\hline\hline  \\ [-1.5ex]
		&Cumulative &Cumulative  &Continuation  &Continuation &Adjacent &Adjacent\\ [-0.5ex]
		&{\it po} &{\it npo} &{\it po} &{\it npo} &{\it po} &{\it npo}\\ [0.5ex]
		\hline\\[-1.5ex]
		AIC &195.87 & 121.17 & 116.40 & {\bf 114.42} & 209.64 & 194.47 \\
		BIC &195.71 & 120.96 & 116.24 & {\bf 114.20} & 209.47 & 194.25 \\[0.5ex]
		\hline\hline 
	\end{tabular}
	\label{tab:compare5}
	\normalsize
\end{table}

\subsection{Lift-one and exchange algorithms}\label{sec:algorithm}

Following \citeauthor{ytm2016} (2017, Section~3), we define 
\[
f_i(z) = f\left(\frac{w_1(1-z)}{1-w_i}, \ldots, \frac{w_{i-1}(1-z)}{1-w_i}, z, \frac{w_{i+1}(1-z)}{1-w_i}, \ldots, \frac{w_m(1-z)}{1-w_i}\right)
\]
with $0 \leq z \leq 1$ and ${\mathbf w} = (w_1, \ldots, w_m)^T \in S_+$. Parallel to Theorem~6 in \cite{ytm2016}, we obtain the following result by Theorem~\ref{thm:FdetGeneral}:

\begin{theorem}\label{thm:f_i(z)}
	Given an approximate allocation ${\bf w}=(w_1, \ldots, w_m)^T \in S_+$ and an $i \in \{1, \ldots, m\}$, for $0 < z < 1$,
	\begin{equation}\label{eq:f_i(z)}
	f_i(z) = (1-z)^{p-J+1} \sum_{j=0}^{J-1} b_j z^j (1-z)^{J-1-j}
	\end{equation}
	\vskip -0.3in
	\begin{equation}\label{eq:f_i'(z)}
	f_i'(z) = (1-z)^{p-J} \sum_{j=1}^{J-1} b_j (j-pz) z^{j-1} (1-z)^{J-1-j} - p b_0 (1-z)^{p-1}
	\end{equation}
	where $b_0=f_i(0)$, $(b_{J-1}, \ldots, b_1)^T={\mathbf B}_{J-1}^{-1}{\mathbf c}$, ${\mathbf B}_{J-1} = (s^{t-1})_{s,t=1, \ldots, J-1}  $ is a $(J-1)\times (J-1)$ constant matrix, and ${\mathbf c}=(c_1, \ldots, c_{J-1})^T$ with $c_j = (j+1)^p j^{J-1-p} f_i(1/(j+1)) - j^{J-1}f_i(0)$, $j=1, \ldots, J-1$.
\end{theorem}

Theorem~\ref{thm:f_i(z)} shows that $f_i(z)$ is an order-$p$ polynomial of $z$. Since $f_i(1)=0$, the solution to maximization of $f_i(z), 0\leq z\leq 1$ can occur only at $z=0$ or $0< z < 1$ such that $f_i'(z)=0$, that is,
\begin{equation}\label{eq:f_i'=0}
\sum_{j=1}^{J-1} j b_j z^{j-1} (1-z)^{J-j-1} = p\sum_{j=0}^{J-1} b_j z^j (1-z)^{J-j-1},\>\>\> 0<z<1.
\end{equation}
This is an order-$(J-1)$ polynomial equation in $z$. For $J\leq 5$, \eqref{eq:f_i'=0} is a polynomial equation of order-$4$ or less, which can be solved analytically. For $J \geq 6$, a quasi-Newton algorithm can be applied for searching numerical solutions.

\bigskip\noindent {\it Lift-one algorithm for D-optimal allocation ${\mathbf w}=(w_1,$ $\ldots,$ $w_m)^T$:}
\begin{itemize}
	\item[$1^\circ$] Start with an arbitrary allocation ${\mathbf w}_0=(w_1,\ldots,w_m)^T$ satisfying $0<w_i<1$, $i=1,\ldots,m$ and compute $f\left({\mathbf w}_0\right)$.
	\item[$2^\circ$] Set up a random order of $i$ going through $\{1,2,\ldots,m\}$.
	\item[$3^\circ$] For each $i$, determine $f_i(z)$ according to Theorem~\ref{thm:f_i(z)}. In this step, $J$ determinants $f_i(0), f_i(1/2),$ $f_i(1/3),$ $\ldots,$ $f_i(1/J)$ are calculated.
	\item[$4^\circ$] Use quasi-Newton algorithm to find $z_*$ maximizing $f_i(z)$ with $0\leq z\leq 1$. If $f_i(z_*) \leq f_i(0)$, let $z_*=0$. Define ${\mathbf w}_*^{(i)} =(w_1(1-z_*)/(1-w_i),$ $\ldots,$ $w_{i-1}(1-z_*)/(1-w_i), z_*, w_{i+1}(1-z_*)/(1-w_i), \ldots, w_m(1-z_*)/(1-w_i))^T$. Note that $f({\mathbf w}_*^{(i)}) = f_i(z_*)$.
	\item[$5^\circ$] Replace ${\mathbf w}_0$ with ${\mathbf w}_*^{(i)}$, and $f\left({\mathbf w}_0\right)$ with $f({\mathbf w}_*^{(i)})$.
	\item[$6^\circ$] Repeat $2^\circ\sim 5^\circ$ until convergence, that is, $f({\mathbf w}_0)=f({\mathbf w}_*^{(i)})$ for each $i$.
\end{itemize}

\vspace{0.5cm}
Following \cite{ymm2016, ytm2016}, we define $$f_{ij}(z) = f(n_1, \ldots, n_{i-1}, z, n_{i+1}, \ldots, n_{j-1}, n_i + n_j - z, n_{j+1}, \ldots, n_m)$$ 
with $z= 0, 1, \ldots, n_i+n_j$ given $1\leq i < j \leq m$ and ${\mathbf n} = (n_1, \ldots, n_m)^T$. As a conclusion of 
Theorem~\ref{thm:FdetGeneral}, Lemma~\ref{lemma:max>=J} and Corollary~\ref{col:calpha>0minimalm}, we obtain the following result:

\begin{theorem}\label{thm:f_ij(z)}
	Suppose ${\mathbf n}=(n_1,\ldots,n_m)^T$ satisfies $f({\mathbf n})>0$ and $n_i+n_j\geq q$ for given $1\leq i< j \leq m$, where $q = \min\{2J-2, p-k_{\min}+2, p\}$. Then
	\begin{equation}\label{eq:f_ij(z)poly}
	f_{ij}(z)=\sum_{s=0}^{q} g_s z^s, \quad z=0, 1, \ldots, n_i+n_j
	\end{equation}
	where $g_0=f_{ij}(0)$, and $g_1, \ldots, g_q$ can be obtained using $(g_1, \ldots, g_q)^T={\mathbf B}_q^{-1}(d_1,$ $\ldots,$ $d_q)^T$ with ${\mathbf B}_q = (s^{t-1})_{s,t=1, \ldots, q}$ as a $q \times q$ constant matrix and $d_s=(f_{ij}(s)-f_{ij}(0))/s$.
\end{theorem}

\bigskip\noindent
{\it Exchange algorithm for D-optimal allocation $(n_1, \ldots, n_m)^T$ given $n>0$:}
\begin{itemize}
	\item[$1^\circ$] Start with an initial allocation ${\mathbf n}=(n_1,\ldots,n_m)^T$ such that $f({\mathbf n})>0$.
	\item[$2^\circ$] Set up a random order of $(i,j)$ going through all pairs $\{(1,2),$ $(1,3),$ $\ldots,$ $(1,m),$ $(2,3),$ $\ldots,(m-1, m)\}$.
	\item[$3^\circ$] For each $(i,j)$, let $c=n_i+n_j$~. If $c=0$, let ${\mathbf n}^*_{ij}={\mathbf n}$. Otherwise, there are two cases. {\it Case one:} $0<c\leq q$, we calculate $f_{ij}(z)$ for $z=0, 1, \ldots, c$ directly and find $z^*$ which maximizes $f_{ij}(z)$. {\it Case two:} $c > q$, we first calculate $f_{ij}(z)$ for $z=0, 1, \ldots, q$; secondly determine $g_0, g_1, \ldots, g_q$ in (\ref{eq:f_ij(z)poly}) according to Theorem~\ref{thm:f_ij(z)}; thirdly calculate $f_{ij}(z)$ for $z=q+1, \ldots, c$ based on (\ref{eq:f_ij(z)poly}); fourthly find $z^*$ maximizing $f_{ij}(z)$ for $z=0, \ldots, c$. For both cases, we define
	$${\mathbf n}^*_{ij} =\left(n_1,\ldots,n_{i-1},z^*,n_{i+1},\ldots,n_{j-1},c-z^*,n_{j+1},\ldots,n_m\right)^T$$
	Note that $f({\mathbf n}^*_{ij})=f_{ij}(z^*)\geq f({\mathbf n})>0$. If $f({\mathbf n}^*_{ij}) > f({\mathbf n})$, replace ${\mathbf n}$ with ${\mathbf n}^*_{ij}$~, and $f({\mathbf n})$ with $f({\mathbf n}^*_{ij})$.
	\item[$4^\circ$] Repeat $2^\circ\sim 3^\circ$ until convergence, that is, $f({\mathbf n}^*_{ij})= f({\mathbf n})$ in step~$3^\circ$ for all $(i,j)$.
\end{itemize}

\subsection{Formulae for calculating $\pi_{ij}$'s from ${\mathbf X}_i$'s}\label{sec:pifromxi}

Following the notations in model~(\ref{logitunifiedmodel}), ${\boldsymbol\eta}_i={\mathbf X}_i{\boldsymbol\theta} = {\mathbf C}^T\log({\mathbf L}{\boldsymbol\pi}_i)$. The formulae towards calculating $\pi_{ij}$'s are listed as follows:\\
(1) Baseline-category logit model
\begin{equation*}
\log({\boldsymbol\pi}_i)=
\begin{pmatrix}
1 & & & & -1\\
& 1 & & & -1\\
& & \ddots  & & \vdots\\
& & & 1 & -1\\
& & & & -1
\end{pmatrix}_{J \times J}
\cdot \log \left(
\begin{pmatrix}
1 & & & & 0\\
& 1 & & & 0\\
& & \ddots  & & \vdots\\
& & & 1 & 0\\
1 & 1 & \cdots & 1 & 1
\end{pmatrix}_{J \times J}
\cdot \exp({\boldsymbol\eta}_i)
\right)
\end{equation*}
(2) Adjacent-categories logit model
\begin{eqnarray*}
	\log({\boldsymbol\pi}_i)&=&
	\begin{pmatrix}
		1 & & & & -1\\
		& 1 & & & -1\\
		& & \ddots  & & \vdots\\
		& & & 1 & -1\\
		& & & & -1
	\end{pmatrix}_{J \times J}
	\cdot \\
	& & \log \left(
	\begin{pmatrix}
		1 & & & & 0\\
		& 1 & & & 0\\
		& & \ddots  & & \vdots\\
		& & & 1 & 0\\
		1 & 1 & \cdots & 1 & 1
	\end{pmatrix}_{J \times J}
	\cdot \exp \left(
	\begin{pmatrix}
		1 & 1 & \cdots & 1 & 0\\
		& 1 & \cdots & 1 & 0\\
		& & \ddots  & \vdots & \vdots\\
		& & & 1 & 0\\
		& & & & 1
	\end{pmatrix}_{J \times J}
	\cdot{\boldsymbol\eta}_i \right)
	\right)
\end{eqnarray*}
(3) Continuation-ratio logit model
\begin{eqnarray*}
	\log({\boldsymbol\pi}_i)&=&{\boldsymbol\eta}_i-
	\begin{pmatrix}
		1 & & & & \\
		1 & 1 & & & \\
		\vdots & & \ddots  & & \\
		1 & 1 & \cdots & 1 & \\
		1 & 1 & \cdots & 1 & 1
	\end{pmatrix}_{J \times J}
	\cdot\\
	& & \log \left(
	\begin{pmatrix}
		1 & 0 & \cdots & 0 & 1\\
		0 & 1 & \cdots & 0 & 1\\
		& & \ddots  & & \\
		0 & 0 & \cdots & 1 & 1\\
		0 & 0 & \cdots & 0 & 1
	\end{pmatrix}_{J \times J}
	\cdot \exp({\boldsymbol\eta}_i)
	\right)
\end{eqnarray*}
(4) Cumulative logit model
\begin{align*}
\log \left(
\begin{pmatrix}
1 & & & & \\
1 & 1 & & & \\
\vdots & & \ddots  & & \\
1 & 1 & \cdots & 1 & \\
0 & 0 & \cdots & 0 & 1
\end{pmatrix}_{J \times J}
\cdot {\boldsymbol\pi}_i \right)
&=\begin{pmatrix}
1 & & &-1 & & \\
&\ddots & & &\ddots & \\
& &1  & & &-1 \\
0 &\cdots &0 & 0&\cdots &-1
\end{pmatrix}_{J \times 2(J-1)}\\
&\cdot \log \left(
\begin{pmatrix}
1 & & & 0 \\
& \ddots & & \vdots \\
& & 1 & 0 \\
1 & & & 1 \\
& \ddots & & \vdots \\
& & 1 & 1
\end{pmatrix}_{ 2(J-1) \times J}
\cdot \exp({\boldsymbol\eta}_i)
\right)
\end{align*}
Note that ${\mathbf X}_i{\boldsymbol\theta}$ in the above models could be {\it po}, {\it npo}, or {\it ppo}.

\subsection{Reparametrization and D-optimality}\label{subsec:reparametrization}

In general, let $\boldsymbol\theta = (\theta_1,$ $\ldots,$ $\theta_p)^T$ be one set of parameters and $\boldsymbol\vartheta=(\vartheta_1,$ $\ldots,$ $\vartheta_p)^T$ be another set of parameters, such that, $\theta_l=h_l(\boldsymbol\vartheta)$, $l=1,\cdots,p$; the map ${\boldsymbol\theta} = \boldsymbol\theta(\boldsymbol\vartheta) = (h_1(\boldsymbol\vartheta), \ldots, h_p(\boldsymbol\vartheta))^T$ is one-to-one; $h_l$'s are differentiable; and the $p \times p$ Jacobian matrix ${\mathbf J} = (h_i(\boldsymbol\vartheta)/\partial \vartheta_j)_{ij}$ is nonsingular.

Consider a design $\xi = \{({\mathbf x}_i, w_i), i=1, \ldots, m\}$ with the distinct experimental settings ${\mathbf x}_i$'s and the corresponding proportions $w_i \in [0,1]$. According to \citeauthor{schervish1995} (1995, page 115), the Fisher information matrix ${\mathbf F}_{\xi}(\boldsymbol\vartheta)$ at $\boldsymbol\vartheta$ and the Fisher information matrix ${\mathbf F}_{\xi}(\boldsymbol\theta)$ at ${\boldsymbol\theta} = \boldsymbol\theta(\boldsymbol\vartheta)$ satisfy
${\mathbf F}_{\xi}(\boldsymbol\vartheta)={\mathbf J}^T{\mathbf F}_{\xi}(\boldsymbol\theta(\boldsymbol\vartheta)){\mathbf J}$.
Then $|{\mathbf F}_{\xi}(\boldsymbol\vartheta)| = |{\mathbf J}|^2 \cdot |{\mathbf F}_{\xi}(\boldsymbol\theta(\boldsymbol\vartheta))|$, where ${\mathbf J}$ contains no design points but parameters. A locally D-optimal design maximizing $|{\mathbf F}_{\xi}(\boldsymbol\vartheta)|$ also maximizes $|{\mathbf F}_{\xi}(\boldsymbol\theta(\boldsymbol\vartheta))|$. That is, it is mathematically equivalent to find D-optimal designs for parameters $\boldsymbol\vartheta$ or $\boldsymbol\theta$.

In terms of Bayesian D-optimal criterion, if a prior distribution of $\boldsymbol\vartheta$ is available, it induces a prior distribution of $\boldsymbol\theta$ since $\boldsymbol\theta=\boldsymbol\theta(\boldsymbol\vartheta)$ is one-to-one.
Then $E_{\boldsymbol\vartheta}\log|{\mathbf F}_{\xi}(\boldsymbol\vartheta)|
= E_{\boldsymbol\vartheta}\log\left|{\mathbf J}^T{\mathbf F}_{\xi}(\boldsymbol\theta(\boldsymbol\vartheta)){\mathbf J}\right|
= E_{\boldsymbol\vartheta}\log|{\mathbf J}|^2+E_{\boldsymbol\vartheta}\log|{\mathbf F}_{\xi}(\boldsymbol\theta(\boldsymbol\vartheta))|$
$=$ $E_{\boldsymbol\vartheta}\log|{\mathbf J}|^2+E_{\boldsymbol\theta}\log|{\mathbf F}_{\xi}(\boldsymbol\theta)|
$.
Therefore, a Bayesian D-optimal design that maximizes $E_{\boldsymbol\theta}\log|{\mathbf F}_{\xi}(\boldsymbol\theta)|$ also maximizes $E_{\boldsymbol\vartheta}\log|{\mathbf F}_{\xi}(\boldsymbol\vartheta)|$.

\begin{example}\label{example:perevozskaya2003}{\rm
\cite{perevozskaya2003} considered the {\it po} model:
\begin{equation}\label{J=3model1}
\log\frac{ \gamma_j(x)}{1- \gamma_j(x)}= {\frac{x-\alpha_j^\prime}{ \beta^\prime}}    \qquad  j=2,\dots,J
\end{equation}
where $\gamma_j(x)=P(Y \geq j|x)$.
Let us reparametrize this model as
\begin{equation}\label{J=3model2}
\log\frac{ \gamma_j(x)}{1- \gamma_j(x)}= \alpha_j+ \beta x   \qquad  j=2,\dots,J
\end{equation}
Let $\boldsymbol\theta=(\alpha_2, \alpha_3, \beta)^T$ be the parameters in \eqref{J=3model1}, and $\boldsymbol\vartheta = (\alpha_2', \alpha_3', \beta')^T$ be the parameters in \eqref{J=3model2}. Then $\beta=1/\beta', \alpha_2=-\alpha_2'/\beta', \alpha_3=-\alpha_3'/\beta'$, and the Jacobian matrix
\begin{equation*}
{\mathbf J}=\begin{pmatrix}
-\frac{1}{\beta^\prime} & 0 & \frac{\alpha_2^\prime}{\beta^{\prime2}}\\
0 & -\frac{1}{\beta^\prime} & \frac{\alpha_3^\prime}{\beta^{\prime2}}\\
0 & 0 &  -\frac{1}{\beta^{\prime2}}
\end{pmatrix}
\end{equation*}
Based on Theorem~\ref{thm:fisher}, the Fisher information $I_i(\boldsymbol\theta)$ at $x_i$ is
\[
\begin{pmatrix}
\frac{\pi_{i1}\pi_{i2,3}^2\pi_{i1,2}}{\pi_{i2}} & -\frac{\pi_{i1}\pi_{i1,2}\pi_{i2,3}\pi_{i3}}{\pi_{i2}} & \pi_{i1}\pi_{i1,2}\pi_{i2,3}x_i\\
-\frac{\pi_{i1}\pi_{i1,2}\pi_{i2,3}\pi_{i3}}{\pi_{i2}} & \frac{\pi_{i1,2}^2\pi_{i2,3}\pi_{i3}}{\pi_{i2}} & \pi_{i3}\pi_{i1,2}\pi_{i2,3}x_i\\
\pi_{i1}\pi_{i1,2}\pi_{i2,3}x_i & \pi_{i3}\pi_{i1,2}\pi_{i2,3}x_i & (\pi_{i1}\pi_{i2,3}^2+\pi_{i2}(\pi_{i1}-\pi_{i3})^2+\pi_{i1,2}^2\pi_{i3})x_i^2
\end{pmatrix}
\]
where $\pi_{ij,k}=\pi_{ij}+\pi_{ik}$~.
It can be verified that $I_i(\boldsymbol\vartheta) = {\mathbf J}^TI_i(\boldsymbol\theta){\mathbf J}$ equals to
the corresponding one given by \cite{perevozskaya2003}. For any given design $\xi = \{({\mathbf x}_i, w_i), i=1, \ldots, m\}$ with proportions $w_i \in [0,1]$, the Fisher information matrix $I_{\xi}(\boldsymbol\vartheta) = \sum_{i=1}^m w_i I_i(\boldsymbol\vartheta)$ $=$ ${\mathbf J}^TI_{\xi}(\boldsymbol\theta){\mathbf J}$.
Then $|I_{\xi}(\boldsymbol\vartheta)| = |{\mathbf J}|^2\cdot |I_{\xi}(\boldsymbol\theta)|$ and the D-optimal design maximizing $|I_{\xi}(\boldsymbol\theta)|$ also maximizes $|I_{\xi}(\boldsymbol\vartheta)|$. That is, the D-optimal designs for Models~(\ref{J=3model1}) and (\ref{J=3model2}) are the same.
}\hfill{$\Box$}
\end{example}

\subsection{More discussion on D-optimality of uniform designs}\label{subsec:uniformDoptimal}

\begin{theorem}\label{thm:J=2}
	Consider Multinomial logit model \eqref{logitunifiedmodel} with only two response categories ($J=2$). In this case, the minimum number of support points is $m=p$. The objective function $f({\mathbf w}) \propto w_1 \cdots w_m$ and the D-optimal allocation among minimally supported designs is ${\mathbf w} = (1/m, \ldots, 1/m)^T$.
\end{theorem}

It can be verified that with $J=2$ all of the four logit models are equivalent to the usual logistic model for binary response. In this case, {\it po}, {\it npo}, or {\it ppo} are essentially the same. Theorem~\ref{thm:J=2} confirms the corresponding results for binary responses in the literature (see, for example, \cite{ym2015}). We provide an independent proof in Section~\ref{sec:proofs}.

Besides the cases with $J=2$, for certain {\it npo} models with $J\geq 3$, uniform allocations could still be D-optimal among minimally supported designs if $p_1 = \cdots = p_{J-1}$. 

\begin{corollary}\label{col:npouniform}
	Consider multinomial logit models \eqref{logitunifiedmodel} with {\it npo} assumption. Suppose $p_1 = \cdots = p_{J-1}$ and there exist $p_1$ distinct experimental settings such that ${\rm rank}({\mathbf H}_1) = \cdots ={\rm rank}({\mathbf H}_{J-1}) = p_1$. Then the minimal number of experimental settings is $m=p_1$ and the uniform allocation is D-optimal among minimally supported designs.
\end{corollary}

According to Corollary~\ref{col:npouniform}, for ``regular" {\it npo} models (that is, $p_1= \cdots = p_{J-1}$), uniform allocations are still D-optimal among minimally supported designs even with $J\geq 3$. However, the following lemma and example further represent that, if the condition $p_1 = \cdots = p_{J-1}$ is violated, uniform allocations are not D-optimal in general even for {\it npo} models.

\begin{lemma}\label{lemma:solutiontoflies}
	Given $0 < c_1 \leq c_2 \leq c_3$, we consider the maximization problem $f(w_1, w_2,$ $w_3) = w_1w_2w_3(c_1w_2w_3 + c_2 w_1 w_3 + c_3 w_1 w_2)$ with respect to $0\leq w_i \leq 1$ and $w_1 + w_2 + w_3=1$.
	Then the solution is $w_1 = w_2 = w_3 = 1/3$ if and only if $c_1=c_2=c_3$.
\end{lemma}

The proof of Lemma~\ref{lemma:solutiontoflies} is relegated to Section~\ref{sec:proofs}, where analytical solutions are provided for $(w_1, w_2, w_3)$ for general values of $c_1$, $c_2$ and $c_3$.

\begin{example}\label{example:atkinson1999} {\rm
		Consider the {\it npo} model adopted by \cite{atkinson1999} with ${\mathbf h}_1(x_i) = (1, x_i, x_i^2)^T$, ${\mathbf h}_2(x_i) = (1, x_i)^T$, $J=3$, $p_1=3, p_2=2$, and $p=5$. According to Corollary~\ref{col:calpha>0minimalm}, the minimum number of support points is $m=\max\{p_1, p_2\} = 3$, which is feasible. The objective function $f({\mathbf w})$ is an order-$5$ polynomial with terms $c_{\alpha_1, \alpha_2, \alpha_3} w_1^{\alpha_1} w_2^{\alpha_2} w_3^{\alpha_3}$. Lemma~\ref{lemma:max>=J} implies that $\alpha_i\in \{0, 1, 2\}, i=1,2,3$ in order to keep $c_{\alpha_1, \alpha_2, \alpha_3} \neq 0$. Combined with Corollary~\ref{col:calpha>0minimalm}, we further know $\alpha_i \in \{1,2\}, i=1,2,3$. According to Theorem~\ref{thm:FdetGeneral}, the objective function is
		\begin{equation}\label{eq:fw123}
		f(w_1, w_2, w_3) = w_1w_2w_3(c_{122}w_2w_3 + c_{212} w_1 w_3 + c_{221} w_1 w_2)
		\end{equation}
		for all the four logit models. Rewriting $(c_{122}, c_{212}, c_{221}) = C\cdot (c_1, c_2, c_3)$, it can be verified that for the continuation-ratio logit model adopted by \cite{atkinson1999} for the house flies experiment (Example~\ref{example:flies}), $C = (x_1 - x_2)^2 (x_1 - x_3)^2 (x_2 - x_3)^2 \prod_{i=1}^3\prod_{j=1}^3 \pi_{ij}$, $c_1 = (x_2-x_3)^2 (\pi_{12}^{-1}+\pi_{13}^{-1})$, $c_2= (x_1-x_3)^2 (\pi_{22}^{-1} + \pi_{23}^{-1})$, $c_3 = (x_1-x_2)^2 (\pi_{32}^{-1} + \pi_{33}^{-1})$; for a cumulative logit model (see, for example, Example~\ref{example:trauma}), $C = (x_1 - x_2)^2 (x_1 - x_3)^2 (x_2 - x_3)^2 \prod_{i=1}^3 \pi_{i1}\pi_{i2}^{-1}\pi_{i3} (\pi_{i1} + \pi_{i2})^2 (\pi_{i2} + \pi_{i3})^2$, $c_1 = (x_2 - x_3)^2 \pi_{13}^{-1} (\pi_{11} + \pi_{12})^{-1}$, $c_2 = (x_1 - x_3)^2 \pi_{23}^{-1} (\pi_{21} + \pi_{22})^{-1}$, and $c_3 = (x_1 - x_2)^2 \pi_{33}^{-1} (\pi_{31} + \pi_{32})^{-1}$.
		According to Lemma~\ref{lemma:solutiontoflies}, $w_1=w_2=w_3=1/3$ is D-optimal if and only if $c_1 = c_2 = c_3$, which is in general not true for both continuation-ratio and cumulative logit models with non-proportional odds.
	}\hfill{$\Box$}
\end{example}

\subsection{More examples}\label{subsec:moreexamples}

\begin{example}\label{ex:maineffects} (For Section~\ref{sec:positive}) {\rm
Consider an experiment with a main-effects multinomial logistic model with $d$ factors and $m$ distinct experimental settings ${\mathbf x}_1, \ldots, {\mathbf x}_m$, where ${\mathbf x}_i = (x_{i1}, \ldots, x_{id})^T$, $i=1, \ldots, m$. 

For a main-effects model, the linear predictors may take the form of
\begin{equation}\label{eq:main-effects}
\eta_{ij} = \beta_{j1} + \beta_{j2} x_{i1} + \cdots + \beta_{j,k+1} x_{ik} + \zeta_1 x_{i,k+1} + \cdots + \zeta_{d-k} x_{id}
\end{equation}
where $i=1, \ldots, m$, $j=1, \ldots, J-1$. In other words, the intercept and the coefficients of the first $k$ factors depend on $j$, while the coefficients of the last $d-k$ factors do not.

We claim that the minimum number of experimental settings is simply $d+1$ for the main-effects multinomial logistic model~\eqref{eq:main-effects} with $0\leq k\leq d$, regardless of $J$.

Actually, first we consider $1\leq k\leq d-1$. It is a {\it ppo} model. In this case, $p_1 = \cdots = p_{J-1} = k+1$, $p_c = d-k$, 
\[
{\mathbf H}_1 = \cdots = {\mathbf H}_{J-1}
=\left(\begin{array}{ccc}
1 & \cdots & 1\\
x_{11} & \cdots & x_{m1}\\
\vdots & \vdots & \vdots\\
x_{1k} & \cdots & x_{mk}
\end{array}\right),\>
{\mathbf H}_c 
=\left(\begin{array}{ccc}
x_{1,k+1} & \cdots & x_{m,k+1}\\
\vdots & \vdots & \vdots\\
x_{1d} & \cdots & x_{md}
\end{array}\right)
\]
According to the special case of Theorem~\ref{thm:minimalm}, the Fisher information matrix ${\mathbf F}$ is positive definite if and only if $m\geq p_c+p_1 = d+1$ and the matrix
\[
({\mathbf H}^T_1, {\mathbf H}^T_c) = 
=\left(\begin{array}{cccc}
1 & x_{11} & \cdots & x_{1d}\\
\vdots & \vdots & \vdots & \vdots\\
1 & x_{m1} & \cdots & x_{md}
\end{array}\right)
\]
is of full rank $d+1$.

Now we let $k=0$. The model~\eqref{eq:main-effects} leads to a {\it po} model. By applying Theorem~\ref{thm:minimalm_po}, we obtain the same conditions as for the {\it ppo} model. Similarly, if we let $k=d$ and apply Theorem~\ref{thm:minimalm_npo}, we get the same conditions for {\it npo} models.
\hfill{$\Box$}  
}\end{example}

\begin{example}\label{ex:m=maxpj} (For Section~\ref{sec:positive}) {\rm
		Consider an experiment with four factors $(d=4)$, three response categories $(J=3)$, and four distinct experimental settings $(m=4)$. Then the experimental settings are ${\mathbf x}_i = (x_{i1}, x_{i2}, x_{i3}, x_{i4})^T$, $i=1,2,3,4$. Consider a multinomial logistic model with {\it ppo} such that
		\[
		{\mathbf H}_1^T = \left(\begin{array}{cccc}
		1 & x_{11} & x_{12} & x_{13}\\
		1 & x_{21} & x_{22} & x_{23}\\
		1 & x_{31} & x_{32} & x_{33}\\
		1 & x_{41} & x_{42} & x_{43}\\
		\end{array}\right),\>
		{\mathbf H}_2^T = \left(\begin{array}{cc}
		1 & x_{11}\\
		1 & x_{21}\\
		1 & x_{31}\\
		1 & x_{41}\\
		\end{array}\right),\>
		{\mathbf H}_c^T = \left(\begin{array}{c}
		x_{14}\\
		x_{24}\\
		x_{34}\\
		x_{44}
		\end{array}\right)
		\]
		That is, $p_1=4, p_2=2, p_c=1, p_H=2$, $\max\{p_1, p_2, p_c+p_H\}= p_1 = 4$, and there are $p=p_1 + p_2 + p_c =7$ parameters. In this case,
		\[
		{\mathbf H} =
		\left(\begin{array}{cc}
		{\mathbf H}_1 & \\
		& {\mathbf H}_2 \\
		{\mathbf H}_c & {\mathbf H}_c
		\end{array}\right)
		\]
		is $7 \times 8$ with rank $7$. That is, the minimum number in Theorem~\ref{thm:minimalm}, $m=\max\{p_1, \ldots, p_{J-1},$ $p_c+p_H\}=4$, is attained in this case.
		\hfill{$\Box$}
}\end{example}

\begin{example}\label{ex:m=maxpc+pH} (For Section~\ref{sec:positive}) {\rm
Consider an experiment with three factors $(d=3)$, three response categories $(J=3)$, and three distinct experimental settings $(m=3)$. Denote the experimental settings as ${\mathbf x}_i = (x_{i1}, x_{i2}, x_{i3})^T$, $i=1,2,3$. Consider a multinomial logistic model with {\it ppo} such that
\[
{\mathbf H}_1^T = \left(\begin{array}{cc}
1 & x_{11} \\
1 & x_{21} \\
1 & x_{31} \\
\end{array}\right),\>
{\mathbf H}_2^T = \left(\begin{array}{c}
1 \\
1 \\
1 \\
\end{array}\right),\>
{\mathbf H}_c^T = \left(\begin{array}{cc}
 x_{12} & x_{13}\\
 x_{22} & x_{23}\\
 x_{32} & x_{33}
\end{array}\right)
\]
That is, $p_1=2, p_2=1, p_c=2, p_H=1$, $\max\{p_1, p_2, p_c+p_H\} = p_c+p_H = 3$, and there are $p=p_1 + p_2 + p_c =5$ parameters. In this case,
\[
{\mathbf H} =
\left(\begin{array}{cccccc}
1 & 1 & 1 & 0 & 0 & 0\\
x_{11} & x_{21} & x_{31} & 0 & 0 & 0\\
0 & 0 & 0 & 1 & 1 & 1\\
x_{12} & x_{22} & x_{32} & x_{12} & x_{22} & x_{32}\\
x_{13} & x_{23} & x_{33} & x_{13} & x_{23} & x_{33}
\end{array}\right)
\]
is $5\times 6$. It can be verified that ${\rm rank}({\mathbf H})= 5$ using Theorem~\ref{thm:rankH}. That is, the minimal number of experimental settings in this case is $m=\max\{p_1, \ldots, p_{J-1}, p_c+p_H\}=3$.
\hfill{$\Box$}
}\end{example}

\begin{example}\label{example:ppopc} (For Section~\ref{sec:determinant}) {\rm
Consider an example with responses in $J=4$ categories, $d=5$ factors, and $m=5$ distinct experimental settings ${\mathbf x}_i = (x_{i,1}, \ldots, x_{i,5})^T$, $i=1, \ldots, 5$. Suppose a multinomial logistic model with
\[
{\mathbf H}_1^T = \left(\begin{array}{ccc}
1 & x_{11} & x_{12}\\
1 & x_{21} & x_{22}\\
\vdots & \vdots & \vdots\\
1 & x_{51} & x_{52}
\end{array}\right),\>
{\mathbf H}_2^T = \left(\begin{array}{cc}
1 & x_{11} \\
1 & x_{21}\\
\vdots & \vdots\\
1 & x_{51}
\end{array}\right),\>
{\mathbf H}_3^T = \left(\begin{array}{c}
1 \\
1\\
\vdots \\
1
\end{array}\right),\>
{\mathbf H}_c^T = \left(\begin{array}{ccc}
x_{13} & x_{14} & x_{15} \\
x_{23} & x_{24} & x_{25} \\
\vdots & \vdots & \vdots \\
x_{53} & x_{54} & x_{55} \\
\end{array}\right)
\]
is used. That is, $p_1=3, p_2=2, p_3 = 1$, $p_H=1$, $p_c=3$, and $p=9$. In this case, ${\mathbf G}$ defined in Theorem~\ref{thm:Fdiagonal} is $20 \times 9$ and $p_c + p_H = 4$ is the minimum number of $\#\{i\mid \alpha_i > 0\}$ to keep $|G[i_1,\ldots,i_p]|\neq 0$ if $(i_1, \ldots, i_p) \in \Lambda(\alpha_1, \ldots, \alpha_m)$. Actually, $(i_1, \ldots, i_9) = (1,2,3,6,7,8,10,11,12) \in \Lambda (3,3,3,0,0)$ leads to $rank({\mathbf G}[i_1, \ldots, i_9])=8$, while $(1,2,5,6,$ $9,$ $10,13,14,15) \in \Lambda (2,2,2,3,0)$ leads to $rank({\mathbf G}[i_1, \ldots, i_9])=9$. Therefore, $|G[i_1,\ldots,i_9]|$ $\neq$ $0$ in general if $(i_1,$ $\ldots,$ $i_9) \in \Lambda(2,2,2,3,0)$ for such a {\it ppo} model.
}\hfill{$\Box$}
\end{example}

\medskip\noindent
{\bf Example~\ref{example:trauma}.} ({\it continued, for Section~\ref{sec:designspace}}) 
Recall that there are eight parameters with fitted values $\hat{\boldsymbol\beta} = (\hat{\beta}_{11}, \hat{\beta}_{21}, \hat{\beta}_{31}, \hat{\beta}_{41}, \hat{\beta}_{12},$ $\hat{\beta}_{22},$ $\hat{\beta}_{32}, \hat{\beta}_{42})^T = (-0.865, -0.094, 0.706, 1.909,$ $-0.113,$ $-0.269, -0.182, -0.119)^T$. If we treat the fitted parameter values as the assumed values, the design space is ${\cal X} = \{ x \geq 0 \mid \beta_{11} + \beta_{12} x < \beta_{21} + \beta_{22} x < \beta_{31} + \beta_{32} x < \beta_{41} + \beta_{42} x\} = \{ x\geq 0 \mid -9.195 < x <  4.942\} = [0, 4.942)$. It is not a surprise that the four levels $\{1, 2, 3, 4\}$ in the original dataset are included in the design space.
\hfill{$\Box$}

\begin{example}\label{example:d1J3m2po} (For Section~\ref{section:minimallysupported}) {\rm
Consider a multinomial logistic model with proportional odds for responses with $J=3$ categories, $d=1$ factors, and $m=2$ distinct experimental settings $x_1, x_2$. Same as in Example~\ref{example:perevozskaya2003}, the parameters are $\beta_1, \beta_2, \zeta_1$ and the linear predictors
$$\eta_{i1} = \beta_1 + \zeta_1 x_i,\>\> \eta_{i2} = \beta_2 + \zeta_1 x_i,\>\> i=1,2.$$
According to Theorem~\ref{thm:FdetGeneral}, the objective function of allocation $(w_1, w_2)$ is an order-$3$ homogeneous polynomial of $w_1, w_2$ consisting of monomials $c_{\alpha_1, \alpha_2} w_1^{\alpha_1} w_2^{\alpha_2}$ with coefficients $c_{\alpha_1, \alpha_2} \geq 0$. Based on Lemma~\ref{lemma:max>=J} and Corollary~\ref{col:calpha>0minimalm}, $c_{\alpha_1, \alpha_2} \neq 0$ only if $\max\{\alpha_1, \alpha_2\} \leq 2$ and $\#\{i\mid \alpha_i > 0\}=2$, which implies $(\alpha_1, \alpha_2)$ is either $(2,1)$ or $(1,2)$. That is, the objective function is
$$f(w_1, w_2) = w_1 w_2 (c_{21} w_1 + c_{12} w_2),$$
which takes the same form as in Corollary~5.2 in \cite{ytm2016}. If we rewrite $c_{21} = C\cdot c_2$ and $c_{12} = C \cdot c_1$, that is, $f(w_1, w_2) = C\cdot w_1 w_2 (c_2 w_1 + c_1 w_2)$, then for a baseline-category logit model, $C = \pi_{13}\pi_{23} (x_1 - x_2)^2$, $c_2 = \pi_{11}\pi_{12}(1-\pi_{23})$, $c_1 = \pi_{21} \pi_{22} (1-\pi_{13})$; for a cumulative logit model, $C = \pi_{12}^{-1} (1-\pi_{13})(1-\pi_{11}) \pi_{22}^{-1} (1-\pi_{23})(1-\pi_{21}) (x_1 - x_2)^2$, $c_2 = \pi_{11}(1-\pi_{11}) \pi_{13} (1-\pi_{13}) \pi_{22} (1-\pi_{22})$, $c_1 = \pi_{12} (1-\pi_{12}) \pi_{21} (1-\pi_{21}) \pi_{23} (1-\pi_{23})$;
for an adjacent-categories logit model, $C = (x_1 -x_2)^2$,  $c_2 = \pi_{11} \pi_{12} \pi_{13} (\pi_{21} \pi_{22} + \pi_{22}\pi_{23} + 4 \pi_{21} \pi_{23})$, $c_1 = \pi_{21}\pi_{22}\pi_{23} (\pi_{11}\pi_{12} + \pi_{12} \pi_{13} + 4 \pi_{11}\pi_{13})$; for a continuation-ratio logit model, $C = (1-\pi_{11})^{-1} (1-\pi_{21})^{-1} (x_1 - x_2)^2$, $c_2 = \pi_{11}\pi_{12}\pi_{13} (1-\pi_{11}) [\pi_{22}\pi_{23} + \pi_{21} (1-\pi_{21})^2]$, $c_1 = \pi_{21} \pi_{22} \pi_{23} (1-\pi_{21}) [ \pi_{12} \pi_{13} + \pi_{11} (1-\pi_{11})^2]$.
According to Corollary~5.2 in \cite{ytm2016}, the uniform allocation $w_1^*=w_2^*=1/2$ is D-optimal if and only if $c_1=c_2$, which is not true in general for all the four logit models.
}\hfill{$\Box$}
\end{example}

\medskip\noindent
{\bf Example~\ref{example:trauma}.} ({\it continued, for Section~\ref{sec:discussion}}) 
In practice, we may use designs not as extreme as the D-optimal design. Here are some alternative allocations of subjects, along with efficiencies:

\begin{table}[h]
	\caption{Alternative Designs for Trauma Clinical Trial}
	\centering 
	\footnotesize
	\begin{tabular}{lrrrrr} \hline\hline
Design point $x$  &     1   &     2   &    3   &     4  &     Efficiency(\%)\\ \hline         
D-optimal design  &   401   &     0   &    0   &   401  &                100.0\\
Original design   &   210   &   190   &  207   &   195  &                 74.7\\
1\% reallocated   &   397   &     4   &    4   &   397  &                 99.4\\
2.5\% reallocated &   391   &    10   &   10   &   391  &                 98.8\\
5\% reallocated   &   381   &    20   &   20   &   381  &                 97.6\\
10\% reallocated  &   361   &    40   &   40   &   361  &                 95.3\\ \hline\hline
\end{tabular}
	\label{tab:trauma.list}
\normalsize
\end{table}

We may recommend 2.5\% or 5\% reallocated design, which is not so extreme but still highly efficient.

\subsection{Proofs}\label{sec:proofs}

$\ $

\medskip\noindent
{\bf Proof of Theorem~\ref{thm:fisher}:}

Suppose for distinct ${\mathbf x}_i, i=1,\cdots,m$, we have independent multinomial responses
\[
{\mathbf Y}_i=(Y_{i1},\cdots,Y_{iJ})^T \sim {\rm Multinomial}(n_i; \pi_{i1},\cdots,\pi_{iJ})
\]
where $n_i=\sum_{j=1}^J Y_{ij}$. Then the log-likelihood for the multinomial model is
\begin{align*}
l(\boldsymbol\theta)&=\log L(\boldsymbol\theta)\\
&=\log \prod_{i=1}^m \frac{{n_i}!}{{Y_{i1}}!\cdots{Y_{iJ}}!} \pi_{i1}^{Y_{i1}}\cdots\pi_{iJ}^{Y_{iJ}}\\
&={\rm constant}+\sum_{i=1}^m {\mathbf Y}_i^T\log \boldsymbol\pi_i
\end{align*}
where $\log \boldsymbol\pi_i=(\log \pi_{i1}, \cdots,\log \pi_{iJ})^T$.  Then the score vector
\begin{eqnarray*}
	\frac{\partial l}{\partial \boldsymbol\theta^T}&=&\sum_{i=1}^{m}{\mathbf Y}_i^T{\rm diag}(\boldsymbol\pi_i)^{-1}\frac{\partial \boldsymbol\pi_i}{\partial \boldsymbol\theta^T}\\
	\frac{\partial l}{\partial \boldsymbol\theta}&=&(\frac{\partial l}{\partial \boldsymbol\theta^T})^T=\sum_{i=1}^{m}(\frac{\partial \boldsymbol\pi_i}{\partial \boldsymbol\theta^T})^T{\rm diag}(\boldsymbol\pi_i)^{-1}{\mathbf Y}_i
\end{eqnarray*}

Using the formulae of matrix differentiation, we get
\begin{align*}
\frac{\partial \boldsymbol\pi_i}{\partial \boldsymbol\theta^T}
&= \frac{\partial \boldsymbol\pi_i}{\partial \boldsymbol\eta_i^T} \cdot \frac{\partial \boldsymbol\eta_i}{\partial \boldsymbol\theta^T}\\
&= \left(\frac{\partial \boldsymbol\eta_i}{\partial \boldsymbol\pi_i^T}\right)^{-1} \cdot {\mathbf X}_i\\
&= \left(\frac{\partial [{\mathbf C}^T\log({\mathbf L}\boldsymbol\pi_i)]}{\partial [\log({\mathbf L}\boldsymbol\pi_i)]^T} \cdot \frac{\partial [\log({\mathbf L}\boldsymbol\pi_i)]}{\partial [{\mathbf L}\boldsymbol\pi_i]^T} \cdot \frac{\partial [{\mathbf L}\boldsymbol\pi_i]}{\partial \boldsymbol\pi_i^T}\right)^{-1} \cdot {\mathbf X}_i\\
&=\left( {\mathbf C}^T[{\rm diag}({\mathbf L}\boldsymbol\pi_i)]^{-1}{\mathbf L} \right)^{-1}{\mathbf X}_i
\end{align*}

\begin{lemma}\label{lemma:pii}
	\[ {\boldsymbol\pi_i}^T {\rm diag}(\boldsymbol\pi_i)^{-1}({\mathbf C}^T{\mathbf D}_i^{-1}{\mathbf L})^{-1}{\mathbf X}_i=\boldsymbol 0^T  \]
\end{lemma}

\noindent
{\bf Proof of Lemma~\ref{lemma:pii}:}
Recall that $\boldsymbol1^T\boldsymbol\pi_i=\pi_{i1}+\cdots + \pi_{iJ}=1$ for each $i$; the last row of ${\mathbf X}_i$ is all $0$; and
\begin{equation*}
{\mathbf C}^T=\begin{pmatrix}
*&*&\cdots&0\\
*&*&\cdots&0\\
\ &\ddots\\
0&0&\cdots&1
\end{pmatrix},\qquad
{\mathbf L}=\begin{pmatrix}
*&*&\cdots&*\\
*&*&\cdots&*\\
\ &\ddots\\
1&1&\cdots&1
\end{pmatrix}\qquad
\end{equation*}
Then
\begin{equation*}
{\mathbf D}_i^{-1}={\rm diag}({\mathbf L}\boldsymbol\pi_i)^{-1}=\begin{pmatrix}
*&0&\cdots&0\\
0&*&\cdots&0\\
\ &\ddots\\
0&0&\cdots&\frac{1}{\boldsymbol1^T\boldsymbol\pi_i}
\end{pmatrix}
=\begin{pmatrix}
*&0&\cdots&0\\
0&*&\cdots&0\\
\ &\ddots\\
0&0&\cdots& 1
\end{pmatrix}
\end{equation*}
and
\begin{equation*}
{\mathbf D}_i^{-1}{\mathbf L}=\begin{pmatrix}
*&\cdots&*\\
*&\cdots&*\\
&\ddots\\
&\boldsymbol1^T
\end{pmatrix}\mbox{  and  }
{\mathbf C}^T{\mathbf D}_i^{-1}{\mathbf L}=\begin{pmatrix}
*&\cdots&*\\
*&\cdots&*\\
&\ddots\\
&\boldsymbol1^T
\end{pmatrix}
\end{equation*}
Rewrite $({\mathbf C}^T{\mathbf D}_i^{-1}{\mathbf L})^{-1}=(\boldsymbol c_{i1},\cdots,\boldsymbol c_{iJ})$.
Then $\boldsymbol 1^T\boldsymbol c_{i1}=\cdots=\boldsymbol 1^T\boldsymbol c_{i,J-1}=0$ and  $\boldsymbol 1^T\boldsymbol c_{iJ}=1$ (just check the last row of ${\mathbf C}^T{\mathbf D}_i^{-1}{\mathbf L}$).
Since $\boldsymbol\pi_i^T {\rm diag}(\boldsymbol\pi_i)^{-1}=(1,\cdots,1)$, then
\[{\boldsymbol\pi_i}^T {\rm diag}(\boldsymbol\pi_i)^{-1} ({\mathbf C}^T{\mathbf D}_i^{-1}{\mathbf L})^{-1}=(1,\cdots,1)(\boldsymbol c_{i1},\cdots,\boldsymbol c_{iJ})=(0,\cdots,0,1)
\]
Since the last row of ${\mathbf X}_i$ is all $0$, then
${\boldsymbol\pi_i}^T {\rm diag}(\boldsymbol\pi_i)^{-1}({\mathbf C}^T{\mathbf D}_i^{-1}{\mathbf L})^{-1}{\mathbf X}_i=\boldsymbol 0^T$.
\hfill{$\Box$}

\medskip
As a direct conclusion of Lemma~\ref{lemma:pii},
\[
E(\frac{\partial l}{\partial \boldsymbol\theta^T})
=\sum_{i=1}^m n_i{\boldsymbol\pi_i}^T {\rm diag}(\boldsymbol\pi_i)^{-1}({\mathbf C}^T{\mathbf D}_i^{-1}{\mathbf L})^{-1}{\mathbf X}_i
=\boldsymbol0^T
\]
Then the Fisher information matrix (see, for example, \citeauthor{schervish1995} (1995, Section~2.3.1))
\begin{eqnarray*}
	{\mathbf F} &=& {\rm Cov}\left(\frac{\partial l}{\partial \boldsymbol\theta},\> \frac{\partial l}{\partial \boldsymbol\theta}\right)
	= E\left(\frac{\partial l}{\partial \boldsymbol\theta}\cdot\frac{\partial l}{\partial \boldsymbol\theta^T}\right)\\
	&=& E\left(\sum_{i=1}^{m}(\frac{\partial \boldsymbol\pi_i}{\partial \boldsymbol\theta^T})^T {\rm diag}(\boldsymbol\pi_i)^{-1}{\mathbf Y}_i\cdot\sum_{j=1}^{m}{\mathbf Y}_j^T{\rm diag}(\boldsymbol\pi_j)^{-1}\frac{\partial \boldsymbol\pi_j}{\partial \boldsymbol\theta^T}\right)\\
	&=& E\left(\sum_{i=1}^{m}\sum_{j=1}^{m}(\frac{\partial \boldsymbol\pi_i}{\partial \boldsymbol\theta^T})^T{\rm diag}(\boldsymbol\pi_i)^{-1}{\mathbf Y}_i{\mathbf Y}_j^T{\rm diag}(\boldsymbol\pi_j)^{-1}\frac{\partial \boldsymbol\pi_j}{\partial \boldsymbol\theta^T}\right)
\end{eqnarray*}
Since ${\mathbf Y}_i$'s follow independent multinomial distributions, then
\begin{align*}
E({\mathbf Y}_i{\mathbf Y}_i^T)
&=\begin{pmatrix}
n_i(n_i-1)\pi_{i1}^2+n_i\pi_{i1}&\cdots&n_i(n_i-1)\pi_{is}\pi_{it}\\
\vdots &\ddots &\vdots\\
n_i(n_i-1)\pi_{is}\pi_{it}&\cdots&n_i(n_i-1)\pi_{iJ}^2+n_i\pi_{iJ}
\end{pmatrix}\\
&=n_i(n_i-1)\boldsymbol\pi_i\boldsymbol\pi_i^T+n_i{\rm diag}(\boldsymbol\pi_i)
\end{align*}
On the other hand, for $i \neq j$,
\[
E({\mathbf Y}_i{\mathbf Y}_j^T)=E({\mathbf Y}_i) \cdot E({\mathbf Y}_j^T)=n_in_j\boldsymbol\pi_i\boldsymbol\pi_j^T
\]
Then the Fisher information matrix
\begin{align*}
{\mathbf F} &=\sum_{i=1}^{m}(\frac{\partial \boldsymbol\pi_i}{\partial \boldsymbol\theta^T})^T {\rm diag}(\boldsymbol\pi_i)^{-1}n_i(n_i-1)\boldsymbol\pi_i\boldsymbol\pi_i^T {\rm diag}(\boldsymbol\pi_i)^{-1}\frac{\partial \boldsymbol\pi_i}{\partial \boldsymbol\theta^T}    \\
&+\sum_{i=1}^{m}(\frac{\partial \boldsymbol\pi_i}{\partial \boldsymbol\theta^T})^T{\rm diag}(\boldsymbol\pi_i)^{-1}n_i {\rm diag}(\boldsymbol\pi_i) {\rm diag}(\boldsymbol\pi_i)^{-1}\frac{\partial \boldsymbol\pi_i}{\partial \boldsymbol\theta^T}   \\
&+\sum_{i \neq j}(\frac{\partial \boldsymbol\pi_i}{\partial \boldsymbol\theta^T})^T {\rm diag}(\boldsymbol\pi_i)^{-1}n_in_j\boldsymbol\pi_i\boldsymbol\pi_j^T{\rm diag}(\boldsymbol\pi_j)^{-1}\frac{\partial \boldsymbol\pi_j}{\partial \boldsymbol\theta^T} \\
&\triangleq (a)+(b)+(c)
\end{align*}
where
\begin{align*}
(b)&=\sum_{i=1}^{m}(\frac{\partial \boldsymbol\pi_i}{\partial \boldsymbol\theta^T})^T {\rm diag} (\boldsymbol\pi_i)^{-1}\frac{\partial \boldsymbol\pi_i}{\partial \boldsymbol\theta^T} n_i  \\
(a)+(c)&=\left[\sum_{i=1}^{m}(\frac{\partial \boldsymbol\pi_i}{\partial \boldsymbol\theta^T})^T {\rm diag}(\boldsymbol\pi_i)^{-1}\boldsymbol\pi_i n_i \right]\left[\sum_{i=1}^{m}(\frac{\partial \boldsymbol\pi_i}{\partial \boldsymbol\theta^T})^T {\rm diag}(\boldsymbol\pi_i)^{-1}\boldsymbol\pi_i n_i \right]^T\\
&-\sum_{i=1}^{m}(\frac{\partial \boldsymbol\pi_i}{\partial \boldsymbol\theta^T})^T {\rm diag}(\boldsymbol\pi_i)^{-1}n_i \boldsymbol\pi_i{\boldsymbol\pi_i}^T {\rm diag}(\boldsymbol\pi_i)^{-1}\frac{\partial \boldsymbol\pi_i}{\partial \boldsymbol\theta^T}
\end{align*}
Actually, let
\[{\mathbf E}_i={\boldsymbol\pi_i}^T {\rm diag}(\boldsymbol\pi_i)^{-1}\frac{\partial \boldsymbol\pi_i}{\partial \boldsymbol\theta^T}
={\boldsymbol\pi_i}^T {\rm diag}(\boldsymbol\pi_i)^{-1}({\mathbf C}^T{\mathbf D}_i^{-1}{\mathbf L})^{-1}{\mathbf X}_i
\]
which is $\boldsymbol0^T$ for each $i$ according to Lemma~\ref{lemma:pii}.
Then
\[
(a)+(c)=\left[\sum_{i=1}^{m}n_i{\mathbf E}_i^T\right]\left[\sum_{i=1}^{m}n_i{\mathbf E}_i^T\right]^T-\sum_{i=1}^{m}n_i{\mathbf E}_i^T{\mathbf E}_i=\boldsymbol0_{J \times J}
\]
The arguments above have proved Theorem~\ref{thm:fisher}.
\hfill{$\Box$}

\medskip\noindent
{\bf Proof of Theorem~\ref{thm:fisher_all}:}
Because the last row of ${\mathbf X}_i$ consists of all zeros, the entries in the last row and last column of ${\mathbf U}_i$ actually won't make any difference. In order to simplify the notations in this proof, we rewrite
\begin{eqnarray*}
	{\mathbf h}_{ji}  & \triangleq & {\mathbf h}_j({\mathbf x}_i)  \qquad j=1,\ldots,J-1;\quad i=1,\ldots,m \\
	{\mathbf h}_{ci}  & \triangleq & {\mathbf h}_c({\mathbf x}_i)  \qquad  i=1,\ldots,m \\
	u_{sti} & \triangleq & u_{st}({\boldsymbol\pi}_i)  \qquad s,t=1, \ldots, J-1; \quad i=1,\ldots,m \\
	u_{s \cdot i} & \triangleq & \sum_{t=1} ^{J-1}  u_{sti}  \qquad s=1, \ldots, J-1; \quad i=1,\ldots,m \\
	u_{\cdot ti} & \triangleq & \sum_{s=1} ^{J-1}  u_{sti}  \qquad t=1, \ldots, J-1; \quad i=1,\ldots,m \\
	u_{\cdot \cdot i} & \triangleq & \sum_{s=1} ^{J-1} \sum_{t=1} ^{J-1}  u_{sti}  \qquad i=1,\ldots,m \\
\end{eqnarray*}
Based on Corollary $\ref{col:Fi}$, when ${\mathbf X}_i$ takes partial proportional odds form (\ref{eqn:Xi_ppo}), the Fisher information ${\mathbf F}_i={\mathbf X}_i^T {\mathbf U}_i {\mathbf X}_i=$
\[
\left(
\begin{array}{cccc}
u_{11i}{\mathbf h}_{1i}{\mathbf h}_{1i}^T & \cdots & u_{1,J-1,i}{\mathbf h}_{1i}{\mathbf h}_{J-1,i}^T & u_{1 \cdot i}{\mathbf h}_{1i}{\mathbf h}_{ci}^T\\
\vdots & \ddots & \vdots & \vdots\\
u_{J-1,1,i}{\mathbf h}_{J-1,i}{\mathbf h}_{1i}^T & \cdots & u_{J-1,J-1,i}{\mathbf h}_{J-1,i}{\mathbf h}_{J-1,i}^T & u_{J-1 \cdot i}{\mathbf h}_{J-1,i}{\mathbf h}_{ci}^T\\
u_{\cdot 1i}{\mathbf h}_{ci}{\mathbf h}_{1i}^T & \cdots & u_{\cdot J-1,i}{\mathbf h}_{ci}{\mathbf h}_{J-1,i}^T & u_{\cdot \cdot i}{\mathbf h}_{ci}{\mathbf h}_{ci}^T
\end{array}\right)
\]
Then the Fisher information matrix ${\mathbf F}=\sum _{i=1}^m n_i{\mathbf F}_i =$
\[
\left(
\begin{array}{cccc}
\sum\limits _{i=1}^m n_i u_{11i}{\mathbf h}_{1i}{\mathbf h}_{1i}^T & \cdots & \sum\limits _{i=1}^m  n_i u_{1,J-1,i}{\mathbf h}_{1i}{\mathbf h}_{J-1,i}^T & \sum\limits _{i=1}^m  n_i u_{1 \cdot i}{\mathbf h}_{1i}{\mathbf h}_{ci}^T\\
\vdots & \ddots & \vdots & \vdots\\
\sum\limits _{i=1}^m  n_i u_{J-1,1,i}{\mathbf h}_{J-1,i}{\mathbf h}_{1i}^T & \cdots & \sum\limits _{i=1}^m  n_i u_{J-1,J-1,i}{\mathbf h}_{J-1,i}{\mathbf h}_{J-1,i}^T & \sum\limits _{i=1}^m  n_i u_{J-1 \cdot i}{\mathbf h}_{J-1,i}{\mathbf h}_{ci}^T\\
\sum\limits _{i=1}^m  n_i u_{\cdot 1i}{\mathbf h}_{ci}{\mathbf h}_{1i}^T & \cdots & \sum\limits _{i=1}^m  n_i u_{\cdot J-1,i}{\mathbf h}_{ci}{\mathbf h}_{J-1,i}^T & \sum\limits _{i=1}^m  n_i u_{\cdot \cdot i}{\mathbf h}_{ci}{\mathbf h}_{ci}^T
\end{array}\right)
\]
or simply
\[
\left(\begin{array}{ccc}
{\mathbf H}_{1} & &\\  & \ddots & \\& &{\mathbf H}_{J-1} \\{\mathbf H}_{c} & \cdots & {\mathbf H}_{c}
\end{array}\right)
\left(\begin{array}{ccc}
{\mathbf U}_{11} & \cdots & {\mathbf U}_{1,J-1}\\
\vdots & \ddots & \vdots\\
{\mathbf U}_{J-1,1} & \cdots & {\mathbf U}_{J-1,J-1}
\end{array}\right)
\left(\begin{array}{cccc}
{\mathbf H}_{1}^T & & & {\mathbf H}_{c}^T\\
& \ddots & & \vdots \\
& &{\mathbf H}_{J-1}^T & {\mathbf H}_{c}^T
\end{array}\right)
\]
\hfill{$\Box$}

\medskip\noindent
{\bf Proof of Theorem~\ref{thm:Upd}:}
Recall that $({\mathbf C}^T {\mathbf D}_i^{-1}{\mathbf L})^{-1} = ({\mathbf c}_{i1} \cdots {\mathbf c}_{iJ})$ and $u_{st}({\boldsymbol\pi}_i)$ $=$ ${\mathbf c}_{is}^T {\rm diag}(\boldsymbol\pi_i)^{-1} {\mathbf c}_{it}$, for $s,t=1, \ldots, J-1$ and $ i=1, \ldots, m$. Denote
\[
\tilde{\mathbf C}=
\begin{pmatrix}
{\mathbf c}_{11}^T  & &\\
&  \ddots  &\\
& & {\mathbf c}_{m1}^T  \\
{\mathbf c}_{12}^T  & &\\
&  \ddots  &\\
& & {\mathbf c}_{m2}^T  \\
\vdots  &  \ddots  & \vdots\\
{\mathbf c}_{1,J-1}^T  & &\\
&  \ddots  &\\
& & {\mathbf c}_{m,J-1}^T
\end{pmatrix}_{m(J-1)\times mJ}
\]
\[
\mbox{ and }\tilde{\mathbf W}=
\begin{pmatrix}
{n_1\rm diag}(\boldsymbol\pi_1)^{-1}  & &\\
&  \ddots  &\\
& & {n_m \rm diag}(\boldsymbol\pi_m)^{-1}
\end{pmatrix}_{mJ \times mJ}
\]
We claim that $ {\mathbf U} = \tilde{\mathbf C} \tilde{\mathbf W} \tilde{\mathbf C}^T $.
Actually
\[
\tilde{\mathbf C} \tilde{\mathbf W} =
\begin{pmatrix}
n_1 {\mathbf c}_{11}^T \rm diag(\boldsymbol\pi_1)^{-1}  & &\\
&  \ddots  &\\
& & n_m {\mathbf c}_{m1}^T \rm diag(\boldsymbol\pi_m)^{-1}\\
\vdots  &  \ddots  & \vdots\\
n_1 {\mathbf c}_{1,J-1}^T \rm diag(\boldsymbol\pi_1)^{-1}  & &\\
&  \ddots  &\\
& & n_m {\mathbf c}_{m,J-1}^T \rm diag(\boldsymbol\pi_m)^{-1}
\end{pmatrix}
\]
and
\begin{align*}
\tilde{\mathbf C} \tilde{\mathbf W} \tilde{\mathbf C}^T &= \tilde{\mathbf C} \tilde{\mathbf W}
\begin{pmatrix}
{\mathbf c}_{11} & & & \cdots & {\mathbf c}_{1,J-1} & &\\
& \ddots & & \ddots &  & \ddots &\\
& & {\mathbf c}_{m1} & \cdots & & & {\mathbf c}_{m,J-1}
\end{pmatrix}\\
&=\begin{pmatrix}
{\mathbf U}_{11} & \cdots & {\mathbf U}_{1,J-1}\\
\vdots & \ddots & \vdots\\
{\mathbf U}_{J-1,1} & \cdots & {\mathbf U}_{J-1, J-1}
\end{pmatrix}={\mathbf U}
\end{align*}
Note that $\tilde{\mathbf W}$ is diagonal with positive diagonal entries. Thus $\tilde{\mathbf W}$ is positive definite. By adjusting the rows, we can verify that $rank(\tilde{\mathbf C})$ is the same as $rank(\tilde{\mathbf C}')$, where
\[
\tilde{\mathbf C}'=
\begin{pmatrix}
{\mathbf c}_{11}^T & & &\\
\vdots & & &\\
{\mathbf c}_{1,J-1}^T & & &\\
& {\mathbf c}_{21}^T & & \\
& \vdots & & \\
& {\mathbf c}_{2,J-1}^T & & \\
& & \ddots &  \\
& & & {\mathbf c}_{m1}^T  \\
& & & \vdots  \\
& & & {\mathbf c}_{m,J-1}^T  \\
\end{pmatrix}
\]
That is, $\tilde{\mathbf C}$ has full row rank and thus ${\mathbf U}$ is positive definite.
\hfill{$\Box$}

\medskip\noindent
{\bf Proof of Theorem~\ref{thm:vdet}:}

\begin{lemma} \label{lemma:udet} $|{\mathbf U}| = (\prod_{i=1}^{m} n_i)^{J-1}|{\mathbf V}|$,
	where
	\begin{align*}
	{\mathbf V} &=\begin{pmatrix}
	{\mathbf V}_{11} & \cdots & {\mathbf V}_{1,J-1}\\
	\vdots & \ddots & \vdots\\
	{\mathbf V}_{J-1,1} & \cdots & {\mathbf V}_{J-1,J-1}
	\end{pmatrix}\\
	&=\begin{pmatrix}
	u_{111} & & &\cdots & u_{1,J-1,1} & &\\
	&\ddots & & &  &\ddots &\\
	& & u_{11m} &  &  & \cdots & u_{1,J-1,m}\\
	\vdots  & &  & \ddots &   & & \vdots \\
	u_{J-1,1,1} & & &\cdots & u_{J-1,J-1,1} & &\\
	&\ddots & & &  &\ddots &\\
	& & u_{J-1,1,m} &  &  & \cdots & u_{J-1,J-1,m}
	\end{pmatrix}
	\end{align*}
\end{lemma}
\noindent

\cite{kovacs1999} generalized Schur's Formula (\cite{gantmacher1960}) as follows:
\begin{lemma}\label{lemma:g_schur}
	(\citeauthor{kovacs1999}, 1999, Theorem~1) \\
	Assume that ${\mathbf M}$ is a $k \times k$ block matrix with each block element ${\mathbf A}_{ij}$ as an $n \times n$ matrix.
	\begin{equation*}
	{\mathbf M} =\begin{pmatrix}
	{\mathbf A}_{11} & \cdots & {\mathbf A}_{1k}\\
	\vdots & \ddots & \vdots\\
	{\mathbf A}_{k1} & \cdots & {\mathbf A}_{kk}
	\end{pmatrix}
	\end{equation*}
	If all of ${\mathbf A}_{ij}$'s commute pairwise, that is, ${\mathbf A}_{ij}{\mathbf A}_{lm} = {\mathbf A}_{lm}{\mathbf A}_{ij}$ for all possible pairs of indices $i$, $j$ and $l$, $m$. Then
	\begin{equation}\label{eqn:g_schur}
	|\mathbf M|=\left| \sum_{\pi \in S_k} (sgn \pi){\mathbf A}_{1\pi(1)}{\mathbf A}_{2\pi(2)} \cdots {\mathbf A}_{k\pi(k)} \right|
	\end{equation}
	Here the sum is computed over all permutations $\pi$ of $\{1, 2, ..., k\}$.
\end{lemma}
In our case, all of ${\mathbf V}_{ij}$'s are diagonal matrices, so they commute pairwise. Moreover, the sum of product matrices in Equation~\eqref{eqn:g_schur} is a diagonal matrix, in which each element is the sum of products of the corresponding elements in those matrices. If we apply the above lemma, we get
\begin{align*}
|\mathbf V|&=\left| \sum_{\pi \in S_{J-1}} (sgn \pi){\mathbf V}_{1\pi(1)}{\mathbf V}_{2\pi(2)} \cdots {\mathbf V}_{J-1,\pi(J-1)} \right|\\
&= \prod_{i=1}^m \left| \sum_{\pi \in S_{J-1}} (sgn \pi)u_{1\pi(1)i}u_{2\pi(2)i} \cdots u_{J-1,\pi(J-1),i} \right|
\end{align*}
Then the following result is obtained:
\begin{lemma}\label{lemma:v=vi}
	$|{\mathbf V}| = \prod_{i=1}^{m} |{\mathbf V}_{i}|$, where
	\[
	{\mathbf V}_i=\left(
	\begin{array}{ccc}
	u_{11}({\boldsymbol\pi}_i) & \cdots & u_{1,J-1}({\boldsymbol\pi}_i) \\
	\vdots & \ddots & \vdots \\
	u_{J-1,1}({\boldsymbol\pi}_i) & \cdots & u_{J-1,J-1}({\boldsymbol\pi}_i)
	\end{array}\right)
	\]
\end{lemma}
\noindent
Note that ${\mathbf V}_i$ defined above is very similar to ${\mathbf U}_i$ define in equation~\eqref{eqn:Ui}.

\begin{lemma}\label{lemma:vi=CDL}
	$|{\mathbf V}_i| = \left(\prod_{j=1}^{J} \pi_{ij}\right)^{-1} \cdot |{\mathbf C}^T {\mathbf D}_i^{-1} {\mathbf L}|^{-2}$.
\end{lemma}

\noindent
{\bf Proof of Lemma~\ref{lemma:vi=CDL}:}
It can be verified that ${\mathbf c}_{iJ} = {\boldsymbol\pi}_i$.\\
Since ${\bf c}_{ij}^T {\rm diag}({\boldsymbol\pi}_i)^{-1} {\bf c}_{iJ}$ $=$ ${\bf c}_{ij}^T {\bf 1} = 0$ for $j=1,\ldots, J-1$ and $1$ for $j=J$, then
\[
\left[\left({\mathbf C}^T {\mathbf D}_i^{-1} {\mathbf L}\right)^{-1}\right]^T
{\rm diag}({\boldsymbol\pi}_i)^{-1}
\left[\left({\mathbf C}^T {\mathbf D}_i^{-1} {\mathbf L}\right)^{-1}\right]
=\left[\begin{array}{cc}
{\mathbf V}_i & {\bf 0}\\
{\bf 0}^T & 1
\end{array}\right]
\]
\hfill{$\Box$}

Combining Lemmas~\ref{lemma:udet}, \ref{lemma:v=vi}, and \ref{lemma:vi=CDL}, we obtain Theorem~\ref{thm:vdet}.

\begin{remark}{\rm Actually, we provide an explicit formula for $|{\mathbf C}^T {\mathbf D}_i^{-1} {\mathbf L}|$ in \eqref{eq:CDiLdet}, which can further clarify Lemma~\ref{lemma:vi=CDL} as (1) $|{\mathbf V}_i| = \prod_{j=1}^J \pi_{ij}$ for baseline-category, adjacent-categories, and continuation-ratio logit models; (2) $|{\mathbf V}_i| = \pi_{iJ}^{-1} \prod_{j=1}^{J-1} \pi_{ij}^{-1} \gamma_{ij}^2 (1-\gamma_{ij})^2$ for cumulative logit models.}
\end{remark}

\medskip\noindent
{\bf Proof of Theorem~\ref{thm:rankH}:}

The simplest case is the {\it npo} model whose conclusion is straightforward.

The {\it ppo} model is the most general case. In this case, we consider a sequence of linear subspaces
\[
\{0\} \subset {\cal M}(\mathbf H_c^T)\cap ( \cap_{j=1}^{J-1}{\cal M}(\mathbf H_j^T)) \subset  M(\mathbf H_c^T)
\]
with corresponding dimensions $0 \leq r_c-r_0 \leq r_c \triangleq rank(\mathbf H_c)$, where $r_0 = rank({\mathbf H_{c}})-dim [{\cal M}(\mathbf H_c^T)\cap ( \cap_{j=1}^{J-1}{\cal M}(\mathbf H_j^T))]$.
Then there exist ${\boldsymbol\alpha}_1, \cdots, {\boldsymbol\alpha}_{r_c-r_0},$ ${\boldsymbol\alpha}_{r_c-r_0+1},$ $\cdots, {\boldsymbol\alpha}_{r_c} \in \mathds{R}^m$ s.t. $\{{\boldsymbol\alpha}_1, \cdots, {\boldsymbol\alpha}_{r_c-r_0}\}$ forms a basis of ${\cal M}(\mathbf H_c^T)\cap ( \cap_{j=1}^{J-1}{\cal M}(\mathbf H_j^T))$ and $\{{\boldsymbol\alpha}_1, \cdots, {\boldsymbol\alpha}_{r_c}\}$ forms a basis of ${\cal M}(\mathbf H_c^T)$. By simple operations $\mathbf H_c$ can be transformed into $\mathbf H_c^*=({\boldsymbol\alpha}_1, \cdots, {\boldsymbol\alpha}_{r_c}, \mathbf 0,$ $\cdots,$ $\mathbf 0)^T$ and $\mathbf H_j$ can be transformed into
\[ \mathbf H_j^*= ({\boldsymbol\alpha}_1, \cdots, {\boldsymbol\alpha}_{r_c-r_0}, {\boldsymbol\alpha}_{r_c-r_0+1}^{(j)}, \cdots, {\boldsymbol\alpha}_{r_j}^{(j)}, \mathbf 0, \cdots, \mathbf 0 )^T \]
where $r_j=rank(\mathbf H_j)$, $j=1,2, \cdots, J-1$. Then $ rank(\mathbf H_{ppo})= rank(\mathbf H_{ppo}^*)$
with
\[ {\mathbf H_{ppo}^*}=\left(\begin{array}{ccc}
{\mathbf H}^*_{1} & &\\  & \ddots & \\& &{\mathbf H}^*_{J-1} \\{\mathbf H}^*_{c} & \cdots & {\mathbf H}^*_{c}
\end{array}\right)_{p \times m(J-1)} \]

Since the first $r_c-r_0$ rows of $({\mathbf H}^*_{c}, \cdots, {\mathbf H}^*_{c})$ can be eliminated by applying row operations of ${\mathbf H}^*_{j}$ onto it separately, then $ rank(\mathbf H_{ppo}^*)= rank(\mathbf H_{ppo}^{**})$ where
\[ {\mathbf H_{ppo}^{**}}=\left(\begin{array}{ccc}
{\mathbf H}^*_{1} & &\\  & \ddots & \\& &{\mathbf H}^*_{J-1} \\{\mathbf H}^{**}_{c} & \cdots & {\mathbf H}^{**}_{c}
\end{array}\right)_{p \times m(J-1)} \]
and $\mathbf H_c^{**}=(\mathbf 0, \cdots, \mathbf 0, {\boldsymbol\alpha}_{r_c-r_0+1}, \cdots, {\boldsymbol\alpha}_{r_c}, \mathbf 0, \cdots, \mathbf 0)^T$. Therefore, $rank(\mathbf H_{ppo})$ $=$ $rank(\mathbf H_{ppo}^{**})$ $\leqslant$ $r_1+\cdots+r_{J-1}+r_0 $.

We claim that the nonzero rows of $\mathbf H_{ppo}^{**}$ are linearly independent which will lead to the final conclusion. Actually, let's denote those nonzero rows of $\mathbf H_{ppo}^{**}$ as ${\boldsymbol\Lambda}_i^{(j)}, i=1,2, \cdots, r_j, j=1,2,\cdots, J-1$ and ${\boldsymbol\Lambda}_{r_c-r_0+1}, \cdots, {\boldsymbol\Lambda}_{r_c}$, where ${\boldsymbol\Lambda}_i^{(j)}$ is the $i$th row of $(\mathbf 0, \cdots, \mathbf 0,$ ${\mathbf H}^*_{j},$ $\mathbf 0,$ $\cdots,$ $\mathbf 0)$, and ${\boldsymbol\Lambda}_i$ is the $i$th row  of $({\mathbf H}^{**}_{c}, \cdots ,{\mathbf H}^{**}_{c})$.
Suppose there exist $a_i^{(j)} \in \mathds{R}, i=1,2, \cdots, r_j, j=1,2,\cdots, J-1$ and $a_i \in \mathds{R}, i=r_c-r_0+1,\cdots,r_c $ s.t.
\[ {\bf 0}=\sum_{j=1}^{J-1}\sum_{i=1}^{r_j}a_i^{(j)}{\boldsymbol\Lambda}_i^{(j)}+ \sum_{i=r_c-r_0+1}^{r_c}a_i{\boldsymbol\Lambda}_i\]
then for $j=1, \ldots, J-1$,
\[ {\bf 0}=\sum_{i=1}^{r_c-r_0} a_i^{(j)} {\boldsymbol\alpha}_i + \sum_{i=r_c - r_0 + 1}^{r_j} a_i^{(j)} {\boldsymbol\alpha}_i^{(j)} + \sum_{i=r_c-r_0+1}^{r_c} a_i {\boldsymbol\alpha}_i\]
which implies for $j=1, \ldots, J-1$,
\[ \sum_{i=r_c-r_0+1}^{r_c} a_i {\boldsymbol\alpha}_i = -\sum_{i=1}^{r_c-r_0} a_i^{(j)} {\boldsymbol\alpha}_i - \sum_{i=r_c - r_0 + 1}^{r_j} a_i^{(j)} {\boldsymbol\alpha}_i^{(j)} \in {\cal M}({\mathbf H}_c^T) \cap {\cal M}({\mathbf H}_j^T)\]
Thus, $\sum_{i=r_c-r_0+1}^{r_c} a_i {\boldsymbol\alpha}_i \in {\cal M}({\mathbf H}_c^T) \bigcap \left(\cap_{j=1}^{J-1} {\cal M}({\mathbf H}_j^T)\right)$. Then we must have $\sum_{i=r_c-r_0+1}^{r_c} a_i {\boldsymbol\alpha}_i$ $=$ ${\bf 0}$ since $\{{\boldsymbol\alpha}_{r_c-r_0+1}, \ldots, {\boldsymbol\alpha}_{r_c}\}$ and $\{{\boldsymbol\alpha}_1, \ldots, {\boldsymbol\alpha}_{r_c-r_0}\}$ are linearly independent. Therefore, $a_i=0$ for $i=r_c-r_0+1, \ldots, r_c$ and thus
\[ {\bf 0}=\sum_{i=1}^{r_c-r_0} a_i^{(j)} {\boldsymbol\alpha}_i + \sum_{i=r_c - r_0 + 1}^{r_j} a_i^{(j)} {\boldsymbol\alpha}_i^{(j)}\]
It implies $a_i^{(j)} =0$, $i=1, \ldots, r_c-r_0, r_c-r_0+1, \ldots, r_j$ since $\{{\boldsymbol\alpha}_1, \ldots, {\boldsymbol\alpha}_{r_c-r_0},$ ${\boldsymbol\alpha}_{r_c-r_0+1}^{(j)},$ $\ldots,$ ${\boldsymbol\alpha}_{r_j}^{(j)}\}$ are linear independent.

Therefore, the conclusion on {\it ppo} models is justified.

Since {\it po} models are special cases of {\it ppo} models, the corresponding result is a direct conclusion.
\hfill{$\Box$}

\medskip\noindent
{\bf Proof of Theorem~\ref{thm:rankintersection}:}

Recall that $dim({\cal M}({\mathbf H}_i^T))=rank({\mathbf H}_i^T)=r_i$ and
$dim({\cal M}({\mathbf H}_{i_1}^T)+\cdots+{\cal M}({\mathbf H}_{i_k}^T))
= dim({\cal M}(({\mathbf H}_{i_1}^T, \cdots, {\mathbf H}_{i_k}^T)))
= rank(({\mathbf H}_{i_1}^T, \cdots,$ ${\mathbf H}_{i_k}^T))$
$=$ $r_{i_1, \ldots, i_k}$, for $i_1 < \cdots < i_k$ and $k=2, \ldots, n$, where ``$+$" stands for the sum of two linear subspaces.
\\
First of all, $dim({\cal M}({\mathbf H}_1^T)\cap {\cal M}({\mathbf H}_2^T))
=dim({\cal M}({\mathbf H}_1^T))+ dim({\cal M}({\mathbf H}_2^T))-dim({\cal M}({\mathbf H}_1^T)+{\cal M}({\mathbf H}_2^T))
=r_1+r_2-r_{12}$. That is, (\ref{eq:nintersections}) is true for $n=2$.
\\
Suppose (\ref{eq:nintersections}) is true for $n=k$. Then for $n=k+1$,
\begin{eqnarray*}
	& & dim(\cap_{i=1}^{k+1}{\cal M}({\mathbf H}_i^T)) = dim(\cap_{i=1}^{k}{\cal M}({\mathbf H}_i^T) \cap {\cal M}({\mathbf H}_{k+1}^T))\\
	&=& dim(\cap_{i=1}^{k}{\cal M}({\mathbf H}_i^T)) + dim({\cal M}({\mathbf H}_{k+1}^T))- dim(\cap_{i=1}^{k}{\cal M}({\mathbf H}_i^T) + {\cal M}({\mathbf H}_{k+1}^T))\\
	&=& \sum_{i=1}^k r_i-\sum_{1 \leqslant i_1 < i_2 \leqslant k} r_{i_1i_2}+ \cdots +(-1)^{k-1}r_{12\cdots k}
	+r_{k+1}-\triangle
\end{eqnarray*}
where
\begin{eqnarray*}
	& &\triangle = dim(\cap_{i=1}^{k}{\cal M}({\mathbf H}_i^T) + {\cal M}({\mathbf H}_{k+1}^T)) = dim(\cap_{i=1}^{k}{\cal M}(({\mathbf H}_i^T, {\mathbf H}_{k+1}^T)))\\
	&=& \sum_{i=1}^k rank(({\mathbf H}_i^T, {\mathbf H}_{k+1}^T))- \sum_{1 \leqslant i_1 < i_2 \leqslant k} rank(({\mathbf H}_{i_1}^T, {\mathbf H}_{k+1}^T, {\mathbf H}_{i_2}^T, {\mathbf H}_{k+1}^T))\\
	& & + \cdots+(-1)^{k-1}rank(({\mathbf H}_1^T, {\mathbf H}_{k+1}^T, \cdots, {\mathbf H}_k^T, {\mathbf H}_{k+1}^T))\\
	&=& \sum_{i=1}^k r_{i,k+1}-\sum_{1 \leqslant i_1 < i_2 \leqslant k} r_{i_1,i_2,k+1}
	+ \cdots + (-1)^{k-1} r_{1,2, \ldots, k+1}
\end{eqnarray*}
Therefore,
\begin{eqnarray*}
	& & dim(\cap_{i=1}^{k+1}{\cal M}({\mathbf H}_i^T))\\
	&=& \sum_{i=1}^k r_i-\sum_{1 \leqslant i_1 < i_2 \leqslant k} r_{i_1i_2}+ \cdots +(-1)^{k-1}r_{12\cdots k}
	+r_{k+1}\\
	& & -  \sum_{i=1}^k r_{i,k+1} + \sum_{1 \leqslant i_1 < i_2 \leqslant k} r_{i_1,i_2,k+1}
	+ \cdots + (-1)^{k} r_{1,2, \ldots, k+1}\\
	&=& \sum_{i=1}^{k+1} r_i- \sum_{1 \leqslant i_1 < i_2 \leqslant k+1}r_{i_1i_2}+ \cdots +(-1)^{(k+1)-1} r_{1,2, \ldots, k+1}
\end{eqnarray*}
That is, (\ref{eq:nintersections}) is true for $n=k+1$ as well.
By mathematical induction, (\ref{eq:nintersections}) is true for general $n$.
\hfill{$\Box$}

\medskip\noindent
{\bf Proof of Corollary~\ref{lemma:ppo_m}:}

Suppose $p_H >0$. Then there exist $m\times 1$ vectors ${\boldsymbol\alpha}_1, \cdots, {\boldsymbol\alpha}_{p_H}$, which form a basis of $\cap_{j=1}^{J-1}{\cal M}({\mathbf H}_j^T)$. Write ${\mathbf H}_c=({\boldsymbol\gamma}_1, \cdots, {\boldsymbol\gamma}_{p_c})^T$. According to Theorem~\ref{thm:rankH}, if $|{\mathbf F}|>0$, then $r_0 = rank({\mathbf H}_c) = p_c$, or equivalently, ${\cal M}({\mathbf H}_c^T) \cap \left(\cap_{j=1}^{J-1}{\cal M}({\mathbf H}_j^T)\right) = \{{\mathbf 0}\}$. Then ${\boldsymbol\alpha}_1, \cdots, {\boldsymbol\alpha}_{p_H}, {\boldsymbol\gamma}_1, \cdots, {\boldsymbol\gamma}_{p_c}$ are linearly independent. Thus $m \geq p_c+p_H$.
\hfill{$\Box$}

\medskip\noindent
{\bf Proof of Theorem~\ref{thm:Fdiagonal}:}

Actually, according to Theorem~\ref{thm:fisher_all}, ${\mathbf F}={\mathbf H}{\mathbf U}{\mathbf H}^T$. From the proof of Theorem~\ref{thm:Upd}, ${\mathbf U} = \tilde{\mathbf C} \tilde{\mathbf W} \tilde{\mathbf C}^T$, where $\tilde{\mathbf W}$ is a diagonal matrix. Therefore, ${\mathbf F} = {\mathbf H} \tilde{\mathbf C} \tilde{\mathbf W} \tilde{\mathbf C}^T {\mathbf H}^T$. Let ${\mathbf W} = \tilde{\mathbf W}/n$ and ${\mathbf G} = \tilde{\mathbf C}^T {\mathbf H}^T$. Then ${\mathbf F}=n {\mathbf G}^T {\mathbf W} {\mathbf G}$, which leads to the final result.
\hfill{$\Box$}

\medskip\noindent
{\bf Proof of Lemma~\ref{lemma:max>=J}:}
Actually, $\max_{1\leq i\leq m} \alpha_i \leq J$. Suppose $\max_{1\leq i\leq m} \alpha_i \geq J$, which means $\max_{1\leq i\leq m} \alpha_i = J$. Without any loss of generality, we assume $\alpha_1=J$. Then $i_j=j$ for $j=1, \ldots, J$.

According to the proof of Lemma~\ref{lemma:pii}, we have ${\mathbf 1}^T {\mathbf c}_{ij} = 0$ for $i=1, \ldots, m$ and $j=1, \ldots, J-1$. Then ${\mathbf 1}^T ({\mathbf c}_{11} + \cdots + {\mathbf c}_{1,J-1})=0$ and thus ${\mathbf 1}^T {\mathbf G}[i_1, \ldots, i_J]$ $=$ $0$. That is, rank$({\mathbf G}[i_1, \ldots, i_J])\leq J-1$. Therefore, rank$({\mathbf G}[i_1, \ldots, i_p]) \leq p-1$ and $|{\mathbf G}[i_1, \ldots, i_p]|=0$.
\hfill{$\Box$}

\medskip\noindent
{\bf Proof of Theorem~\ref{thm:alpha>0}:} Suppose $c_{\alpha_1, \ldots, \alpha_m} \neq 0$ for some $(\alpha_1, \dots, \alpha_m)$. Therefore, there exist $(i_1, \ldots, i_p) \in (\alpha_1, \ldots, \alpha_m)$ such that ${\mathbf G}[i_1, \ldots, i_p]$ is of full rank $p$. Without any loss of generality, we assume $\alpha_1 \geq \cdots \geq \alpha_k  >0=\alpha_{k+1} = \cdots = \alpha_m$, that is, $\{i\mid \alpha_i > 0\} = \{1, \ldots k\}$. Consider the submatrix $\tilde{\mathbf G} := {\mathbf G}[1, \ldots, kJ]$ which is $kJ\times p$ and contains ${\mathbf G}[i_1, \ldots, i_p]$ as a submatrix. Then $\tilde{\mathbf G}$ is of rank $p$ or $\tilde{\mathbf G}^T$ is of full row rank $p$. Write $\tilde{\mathbf W} = k^{-1} {\rm diag}\{{\rm diag}({\boldsymbol\pi}_1)^{-1}, \ldots, {\rm diag} ({\boldsymbol\pi}_k)^{-1}\}$. Then the restricted  matrix ${\mathbf F}$ $:=$ $n$ $\tilde{\mathbf G}^T \tilde{\mathbf W}\tilde{\mathbf G}$ is positive definite. On the other hand, ${\mathbf F}$ is the Fisher information matrix $n {\mathbf G}^T {\mathbf W} {\mathbf G}$ as defined in Theorem~\ref{thm:Fdiagonal} with $w_1 = \cdots = w_k = 1/k$ and $w_{k+1} = \cdots = w_m = 0$. According to Theorem~\ref{thm:Fdiagonal} and Theorem~\ref{thm:fisher}, ${\mathbf F} = nk^{-1} \sum_{i=1}^k {\mathbf F}_i$. Therefore, ${\mathbf F}_{res} := \sum_{i=1}^k {\mathbf F}_i$ is positive definite.
\hfill{$\Box$}

\medskip\noindent
{\bf Proof of Theorem~\ref{thm:designspace}:}

\noindent
{\it Case 1: Baseline-category logit model}

The baseline-category logit model for nominal response (Agresti, 2013; Zocchi and
Atkinson, 1999) can be extended in general as follows
\begin{equation}\label{model:baseline}
\log\left(\frac{\pi_{ij}}{\pi_{iJ}}\right) = {\mathbf h}_j^T({\mathbf x}_i)\boldsymbol\beta_j+{\mathbf h}_c^T({\mathbf x}_i)\boldsymbol\zeta~, \quad j=1, \ldots, J-1
\end{equation}
\begin{lemma}\label{claim4.1}
	Fixing ${\mathbf x}_i$, ${\boldsymbol\beta_j, j=1,\cdots,J-1}$ and $\boldsymbol\zeta$ in Model~(\ref{model:baseline}), let $a_j={\mathbf h}_j^T({\mathbf x}_i)\boldsymbol\beta_j+{\mathbf h}_c^T({\mathbf x}_i)\boldsymbol\zeta$, $j=1, \ldots, J-1$. Then $0<\pi_{ij}<1, j=1, \ldots, J$ exist uniquely if and only if
	$-\infty < a_j < \infty$, $j=1, \ldots, J-1$. In this case,
	\begin{equation}\label{sol:baseline}
	\pi_{ij} = \left\{
	\begin{array}{ll}
	\frac{e^{a_j}}{e^{a_1}+\cdots + e^{a_{J-1}} + 1} & 1 \leq j \leq J-1\\
	\frac{1}{e^{a_1}+\cdots + e^{a_{J-1}} + 1} & j=J
	\end{array}\right.
	\end{equation}
\end{lemma}

\noindent
{\bf Proof of Lemma~\ref{claim4.1}:}
Write $ y_j=\log \pi_{ij}$, $j=1,\ldots, J$. Then $0<\pi_{ij}<1, j=1, \ldots, J$ if and only if $y_j \in (-\infty, 0)$, $j=1, \ldots, J$. In this case, Model~(\ref{model:baseline}) implies
$a_j = y_j-y_J \in (-\infty, \infty)$, $j=1, \ldots, J-1$.\\
On the other hand, for any given $a_1, \ldots, a_{J-1} \in (-\infty, \infty)$, $y_j = a_j + y_J$, $j=1, \ldots, J-1$. Note that
\begin{eqnarray*}
	1 &=& \pi_{i1} + \pi_{i2} + \cdots + \pi_{i,J-1} + \pi_{iJ}\\
	&=& e^{y_1} + e^{y_2} + \cdots + e^{y_{J-1}} + e^{y_J}\\
	&=& e^{a_1 + y_J} + e^{a_2 + y_J} + \cdots + e^{a_{J-1}+y_J} + e^{y_J}\\
	&=& e^{y_J} \left(e^{a_1} + e^{a_2} + \cdots + e^{a_{J-1}} + 1\right)
\end{eqnarray*}
Since $\pi_{ij}=e^{y_j}$, we get solutions of $\pi_{ij}$ given in~(\ref{sol:baseline}), and thus
$\pi_{ij} \in (0, 1)$ exists and is unique, $j=1, \ldots, J$.
\hfill{$\#$}

\noindent
{\it Case 2: Cumulative logit model}

The cumulative logit model for ordinal responses \citep{pmcc1980, christensen2013} can be described in general as follows:
\begin{equation}\label{model:cumulative}
\log\left(\frac{\pi_{i1}+\cdots + \pi_{ij}}{\pi_{i,j+1}+\cdots + \pi_{iJ}}\right) = {\mathbf h}_j^T({\mathbf x}_i)\boldsymbol\beta_j+{\mathbf h}_c^T({\mathbf x}_i)\boldsymbol\zeta~, \quad j=1, \ldots, J-1
\end{equation}

\begin{lemma}\label{claim4.2}
	Fixing ${\mathbf x}_i$, ${\boldsymbol\beta_j, j=1,\cdots,J-1}$ and $\boldsymbol\zeta$ in Model~(\ref{model:cumulative}), let $a_j = {\mathbf h}_j^T({\mathbf x}_i)\boldsymbol\beta_j$ $+$ ${\mathbf h}_c^T({\mathbf x}_i)\boldsymbol\zeta$, $j=1, \ldots, J-1$. Then
	$0< \pi_{ij} <1, j=1, \ldots, J$ exist and are unique if and only if
	$-\infty < a_1 < a_2 < \cdots < a_{J-1} < \infty$. In this case,
	\begin{equation}\label{sol:cumulative}
	\pi_{ij} = \left\{
	\begin{array}{ll}
	\frac{\exp(a_1)}{1+\exp(a_1)} & j=1\\
	\frac{\exp(a_j)}{1+\exp(a_j)} - \frac{\exp(a_{j-1})}{1+\exp(a_{j-1})} & 1<j<J\\
	\frac{1}{1+\exp(a_{J-1})} & j=J
	\end{array}\right.
	\end{equation}
\end{lemma}

\noindent
{\bf Proof of Lemma~\ref{claim4.2}:}
Taking $j=1$ in Model~(\ref{model:cumulative}), then $\log\left(\pi_{i1}/(1-\pi_{i1})\right)$ $=$ $a_1$ and  $\pi_{i1} = \exp(a_1)/[1+\exp(a_1)]$. Then $0<\pi_{i1}<1$ if and only if $-\infty < a_1 < \infty$.
For $j=2, \cdots, J-1$,
\[\pi_{ij}=\frac{\exp(a_j)}{1+\exp(a_j)} - \frac{\exp(a_{j-1})}{1+\exp(a_{j-1})}\]
which implies that $\pi_{ij}> 0$ if and only if  $a_j > a_{j-1}$~. Therefore,
$\pi_{iJ}= 1- (\pi_{i1} + \cdots + \pi_{i,J-1}) = 1 - \exp(a_{J-1})/[1+\exp(a_{J-1})] = 1/[1+\exp(a_{J-1})]$,
which indicates $0< \pi_{iJ} < 1$ if and only if $-\infty < a_{J-1} < \infty$.
Given $\pi_{i1} + \cdots + \pi_{iJ} =1$, we have
\[ -\infty < a_1 < a_2 < \cdots < a_{J-1} < \infty \Leftrightarrow \pi_{ij} \in (0,1), \quad  j=1, \ldots, J\]
\hfill{$\#$}

\begin{corollary}\label{corollary:po}
	For the cumulative logit model with proportional odds
	\begin{equation}
	\log\left(\frac{\pi_{i1}+\cdots+\pi_{ij}}{\pi_{i,j+1}+\cdots+\pi_{iJ}}\right) =\beta_j+ \mathbf h_c^T({\mathbf x}_i)\boldsymbol\zeta~, \quad j=1, \ldots, J-1
	\end{equation}
	The design space has no restriction since $-\infty < \beta_1 < \beta_2 < \cdots < \beta_{J-1} < \infty$ is part of the model assumptions, which implies
$\pi_{ij} \in (0,1)$, $j=1, \ldots, J$.
\end{corollary}

\noindent
{\it Case 3: Adjacent-categories logit model}

The adjacent-categories logit model for ordinal responses \citep{agresti2005, agresti2013} can be extended as follows:
\begin{equation}\label{model:adjacent}
\log\left(\frac{\pi_{ij}}{\pi_{i,j+1}}\right) = {\mathbf h}_j^T({\mathbf x}_i)\boldsymbol\beta_j+{\mathbf h}_c^T({\mathbf x}_i)\boldsymbol\zeta~, \quad j=1, \ldots, J-1
\end{equation}

\begin{lemma}\label{claim4.3}
	Fixing ${\mathbf x}_i$, ${\boldsymbol\beta_j, j=1,\cdots,J-1}$ and $\boldsymbol\zeta$ in Model~(\ref{model:adjacent}), let $a_j = {\mathbf h}_j^T({\mathbf x}_i)\boldsymbol\beta_j$ $+$ ${\mathbf h}_c^T({\mathbf x}_i)\boldsymbol\zeta$, $j=1, \ldots, J-1$. Then $0<\pi_{ij}<1, j=1, \ldots, J$ exist uniquely if and only if
	$-\infty < a_j < \infty$, $j=1, \ldots, J-1$. In this case,
	\begin{equation}\label{sol:adjacent}
	\pi_{ij} = \left\{
	\begin{array}{ll}
	\frac{\exp(a_{J-1}+\cdots+a_j)}{\exp(a_{J-1}+\cdots+a_1)+\exp(a_{J-1}+\cdots+a_2)+\cdots+\exp(a_{J-1})+1} & j=1, \ldots, J-1\\
	\frac{1}{\exp(a_{J-1}+\cdots+a_1)+\exp(a_{J-1}+\cdots+a_2)+\cdots+\exp(a_{J-1})+1} & j=J
	\end{array}\right.
	\end{equation}
\end{lemma}

\noindent
{\bf Proof of Lemma~\ref{claim4.3}:}
Let $ y_j=\log \pi_{ij}$. Then $0<\pi_{ij}<1, j=1, \ldots, J$ if and only if $y_j \in (-\infty, 0)$. In this case, Model~(\ref{model:adjacent}) implies
$a_j = y_j-y_{j+1} \in (-\infty, \infty)$, $j=1, \ldots, J-1$.\\
On the other hand, for any given $a_1, \ldots, a_{J-1} \in (-\infty, \infty)$, $y_j = (a_{J-1} + \cdots +a_j) + y_J$, $j=1, \ldots, J-1$. Note that
\begin{eqnarray*}
	1 &=& \pi_{i1} + \pi_{i2} + \cdots + \pi_{i,J-1} + \pi_{iJ}\\
	&=& e^{y_1} + e^{y_2} + \cdots + e^{y_{J-1}} + e^{y_J}\\
	&=& e^{y_J} \left(e^{a_{J-1} + \cdots +a_1} + e^{a_{J-1} + \cdots +a_2} + \cdots + e^{a_{J-1}} + 1\right)
\end{eqnarray*}
Since $\pi_{ij}=e^{y_j}$, we get solutions of $\pi_{ij}$ given in~(\ref{sol:adjacent}), and thus
$\pi_{ij} \in (0, 1)$ exists and is unique, $j=1, \ldots, J$.
\hfill{$\#$}

\noindent
{\it Case 4: Continuation-ratio logit model}

The continuation-ratio logit model for hierarchical responses (Agresti, 2013; Zocchi and
Atkinson, 1999) can be rewritten in general as follows:
\begin{equation}\label{model:continuation}
\log\left(\frac{\pi_{ij}}{\pi_{i,j+1} + \cdots + \pi_{iJ}}\right) = {\mathbf h}_j^T({\mathbf x}_i)\boldsymbol\beta_j+{\mathbf h}_c^T({\mathbf x}_i)\boldsymbol\zeta~, \quad j=1, \ldots, J-1
\end{equation}

\begin{lemma}\label{claim4.4}
	Fixing ${\mathbf x}_i$, ${\boldsymbol\beta_j, j=1,\cdots,J-1}$ and $\boldsymbol\zeta$ in Model~(\ref{model:continuation}), let $a_j={\mathbf h}_j^T({\mathbf x}_i)\boldsymbol\beta_j$ $+$ ${\mathbf h}_c^T({\mathbf x}_i)\boldsymbol\zeta$, $j=1, \ldots, J-1$. Then $0<\pi_{ij}<1, j=1, \ldots, J$ exist uniquely if and only if
	$-\infty < a_j < \infty$, $j=1, \ldots, J-1$. In this case,
	\begin{equation}\label{sol:continuation}
	\pi_{ij} = \left\{
	\begin{array}{ll}
	e^{a_j}\prod_{s=1}^j\left(e^{a_s}+1\right)^{-1} & j=1, \ldots, J-1\\
	\prod_{s=1}^{J-1}\left(e^{a_s}+1\right)^{-1} & j=J
	\end{array}\right.
	\end{equation}
\end{lemma}

\noindent
{\bf Proof of Lemma~\ref{claim4.4}:}
Let $ y_j=\log \pi_{ij}$. Then $0<\pi_{ij}<1, j=1, \ldots, J$ if and only if $y_j \in (-\infty, 0)$. In this case, Model~(\ref{model:continuation}) implies
$a_j = y_j-\log(e^{y_{j+1}} + \cdots e^{y_J}) \in (-\infty, \infty)$, $j=1, \ldots, J-1$.\\
On the other hand, for any given $a_1, \ldots, a_{J-1} \in (-\infty, \infty)$, it can be verified by induction that
\begin{eqnarray*}
	e^{y_{J-1}} &=& e^{y_J} e^{a_{J-1}}\\
	e^{y_{J-2}} &=& e^{y_J} e^{a_{J-2}}\left(e^{a_{J-1}}+1\right)\\
	e^{y_j} &=& e^{y_J} e^{a_j}\left(e^{a_{j+1}}+1\right)\cdots \left(e^{a_{J-1}}+1\right), \> j=J-3, J-4, \cdots, 1
\end{eqnarray*}
Therefore, it can be verified that
\begin{eqnarray*}
	1 &=& \pi_{i1} + \pi_{i2} + \cdots + \pi_{i,J-1} + \pi_{iJ}\\
	&=& e^{y_1} + e^{y_2} + \cdots + e^{y_{J-1}} + e^{y_J}\\
	&=& e^{y_J} \left(e^{a_1} + 1\right) \left(e^{a_2} + 1\right) \cdots \left(e^{a_{J-1}} + 1\right)
\end{eqnarray*}
Since $\pi_{ij}=e^{y_j}$, we get solutions of $\pi_{ij}$ given in~(\ref{sol:continuation}), and thus
$\pi_{ij} \in (0, 1)$ exists and is unique, $j=1, \ldots, J$.
\hfill{$\#$}

Theorem~\ref{thm:designspace} is obtained as a summary of Lemmas~\ref{claim4.1}, \ref{claim4.2}, \ref{claim4.3}, and \ref{claim4.4}.
\hfill{$\Box$}

\medskip\noindent
{\bf Proof of Corollary~\ref{col:fwu>0}:}
We only need to verity the ``only if" part. According to Theorem~\ref{thm:Fpd}, if $f({\mathbf w})>0$ for some ${\mathbf w} = (w_1, \ldots, w_m)^T = (n_1, \ldots,$ $n_m)^T/n$, then the corresponding ${\mathbf H}^*$ is of full row rank. Note that ${\mathbf H}^*$ can be obtained from ${\mathbf H}$ after removing the columns of ${\mathbf H}$ corresponding to $n_i=0$. Thus ${\mathbf H}$ is of full row rank too, which corresponds to the uniform allocation. That is, $f({\mathbf w}_u)>0$.\\
In this case, any ${\mathbf w} = (w_1, \ldots, w_m)^T$ such that $0 < w_i < 1, i=1, \ldots, m$ leads to $f({\mathbf w}) > 0$ since it corresponds to the same ${\mathbf H}$ matrix.
\hfill{$\Box$}

\medskip\noindent
{\bf Proof of Theorem~\ref{thm:f_ij(z)}:}
According to Theorem~\ref{thm:FdetGeneral},
\[
f_{ij}(z) = \sum_{\alpha_i\geq 0, \alpha_j\geq 0, \alpha_i+\alpha_j\leq p} {\rm coefficient} \cdot z^{\alpha_i} (n_i + n_j - z)^{\alpha_j}
\]
is a polynomial with nonnegative coefficients, whose order depends on the largest possible $\alpha_i+\alpha_j$. Lemma~\ref{lemma:max>=J} implies that $\max\{\alpha_i, \alpha_j\} \leq J-1$ for positive coefficients and Corollary~\ref{col:calpha>0minimalm} further implies that $\alpha_i + \alpha_j \leq p - (k_{\min}-2) = p - k_{\min} +2$ for positive coefficients. Therefore, $f_{ij}(z)$ is at most an order-$q$ polynomial of $z$.
\hfill{$\Box$}

\medskip\noindent
{\bf Proof of Theorem~\ref{thm:J=2}:}
In this case, the model is essentially a generalized linear model for binomial response with logit link. Theorem~\ref{thm:FdetGeneral} says that the objective function $f({\mathbf w}) = |{\mathbf G}^T {\mathbf W} {\mathbf G}|$ is an order-$p$ polynomial consisting of terms $c_{\alpha_1, \ldots, \alpha_m} w_1^{\alpha_1} \cdots w_m^{\alpha_m}$. According to Lemma~\ref{lemma:max>=J}, $c_{\alpha_1, \ldots, \alpha_m} \neq 0$ only if $\alpha_i \in \{0, 1\}, i=1, \ldots, m$. Therefore, in order to keep $f({\mathbf w})>0$, we must have $m\geq p$. In other words, a minimally supported design may contain exactly $m=p$ distinct design points or experimental settings. In this case, the objective function $f({\mathbf w}) \propto w_1 \cdots w_m$ and the D-optimal allocation is ${\mathbf w} = (1/m, \ldots, 1/m)^T$.
\hfill{$\Box$}

\medskip\noindent
{\bf Proof of Corollary~\ref{col:npouniform}:}
According to Theorem~\ref{thm:fisher_all}, ${\mathbf F} = {\mathbf H} {\mathbf U} {\mathbf H}^T$. In this case, there exist $m=p_1$ experimental settings such that ${\rm rank}({\mathbf H}) = p_1(J-1) = p$. On the other hand, the minimum number of experimental settings is at least $\max\{p_1, \ldots, p_{J-1}\} = p_1$ based on Corollary~\ref{col:calpha>0minimalm}. Therefore, the minimal number is $m=p_1$. In this case, ${\mathbf H}$ is a square matrix and
\[
|{\mathbf F}| =  |{\mathbf H}|^2\cdot |{\mathbf U}| \propto \left(\prod_{i=1}^m w_i\right)^{J-1}
\]
according to Theorem~\ref{thm:vdet}. Thus, the uniform allocation ${\mathbf w}_u = (1/m, \ldots,$ $1/m)^T$ is D-optimal in this case. Note that $m=p_1 < p_1(J-1) = p$.
\hfill{$\Box$}

\medskip\noindent
{\bf Proof of Lemma~\ref{lemma:solutiontoflies}:} We actually claim more detailed conclusions as follows:
\begin{itemize}
	\item[(i)] If $c_1=c_2=c_3$, then the solution is $w_1 = w_2 = w_3 = 1/3$.
	\item[(ii)] If $c_1 = c_2 < c_3$, then $w_1 = w_2 > w_3 > 0$. Actually, $w_1 = w_2 = (-2 c_1 + c_3 + \Delta_1)/D_1$ and $w_3 = c_3/D_1$, where $\Delta_1 = \sqrt{4 c_1^2 - c_1 c_3 + c_3^2}$ and $D_1 = -4 c_1 + 3 c_3 + 2 \Delta_1$.
	\item[(iii)] If $c_1 < c_2 = c_3$, then $w_1 > w_2 = w_3 > 0$. Actually, $w_1 = (- c_1 + 2 c_3 + \Delta_2)/D_2$ and $w_2 = w_3 = 3 c_3/D_2$, where $\Delta_2 = \sqrt{c_1^2 - c_1 c_3 + 4 c_3^2}$ and $D_2 = - c_1 + 8 c_3 +  \Delta_2$.
	\item[(iv)] If $c_1 < c_2 < c_3$, then $w_1 > w_2 > w_3 > 0$. The procedure of obtaining analytic solutions of $w_1, w_2, w_3$ is as follows: (1) obtain $y_1$ from \eqref{y1solution}; (2) obtain $y_2$ from \eqref{eq:y2solution}; (3) $w_1 = y_1/(y_1 + y_2 + 1)$, $w_2 = y_2/(y_1 + y_2 + 1)$, $w_3 = 1/(y_1 + y_2 + 1)$.
\end{itemize}

First of all, we only need to consider the cases of $0<w_i<1$, $i=1,2,3$ (otherwise, $f(w_1, w_2, w_3)=0$). It can also be verified that $0 < c_1 \leq c_2 \leq c_3$ implies that $w_1 \geq w_2 \geq w_3 > 0$ (otherwise, for example, if $w_1 < w_2$, one may replace $w_1, w_2$ both with $(w_1 + w_2)/2$ and strictly increase $f$). The same argument implies that if $c_i = c_j$, then $w_i = w_j$ in the solution.

According to Theorem~5.10 in \cite{ytm2016}, $(w_1, w_2, w_3)^T$ maximizes $f(w_1, w_2, w_3)$ if and only if
$$\frac{\partial f}{\partial w_1} = \frac{\partial f}{\partial w_2} = \frac{\partial f}{\partial w_3}$$
which is equivalent to $\partial f/\partial w_1 = \partial f/\partial w_3$ and $\partial f/\partial w_2 = \partial f/\partial w_3$ and thus equivalent to
\begin{equation}\label{eq:f1=f3}
c_3 w_1 w_2 (w_1 - 2 w_3) +  2 c_2 w_1 w_3 (w_1 - w_3) = c_1 w_2 w_3 (-2 w_1 + w_3)
\end{equation}
\begin{equation}\label{eq:f2=f3}
c_3 w_1 w_2 (w_2 - 2 w_3) +  2 c_1 w_2 w_3 (w_2 - w_3) = c_2 w_1 w_3 (-2 w_2 + w_3)
\end{equation}
Following \citeauthor{ytm2016} (2016b, Section 5.2), we denote $y_1 = w_1/w_3 > 0$ and $y_2 = w_2/w_3 > 0$. Actually, $w_1 \geq w_2 \geq w_3 > 0$ implies $y_1 \geq y_2 \geq 1$. Since $w_1 + w_2 +w _3 = 1$, it implies $w_3 = 1/(y_1 + y_2 + 1)$, $w_1 = y_1/(y_1 + y_2 + 1)$, and $w_2 = y_2/(y_1 + y_2 + 1)$. Then \eqref{eq:f1=f3} and \eqref{eq:f2=f3} are equivalent to
\begin{eqnarray}
c_3 y_1 y_2 (y_1 - 2) +  2 c_2 y_1 (y_1 - 1) &=& c_1 y_2 (-2 y_1 + 1) \label{eq:f1=f3y}\\
c_3 y_1 y_2 (y_2 - 2) +  2 c_1 y_2 (y_2 - 1) &=& c_2 y_1 (-2 y_2 + 1) \label{eq:f2=f3y}
\end{eqnarray}
From \eqref{eq:f1=f3y} we get $y_2 [c_3 y_1^2 - 2(c_3 - c_1) y_1 - c_1] = 2 c_2 y_1 (1 - y_1)$. If $y_1 = 1$, then we must have $y_2 = 1$ and $c_3  - 2(c_3 - c_1) - c_1 = 0$, which implies $w_1 = w_2 = w_3 = 1/3$ and $c_1 = c_2 = c_3$. Actually, we can also verify that $c_1=c_3$ implies $y_1=1$.

Now we assume $y_1 > 1$, which implies $c_1 < c_3$. Then
\begin{equation}\label{eq:y2solution}
y_2 = \frac{2 c_2 (1 - y_1) y_1}{c_3 y_1^2 - 2(c_3 - c_1) y_1 - c_1}
\end{equation}
After plugging \eqref{eq:y2solution} into \eqref{eq:f2=f3y}, we get
\begin{equation}\label{eq:y1equation}
a_0 + a_1 y_1 + a_2 y_1^2 + a_3 y_1^3 + y_1^4=0
\end{equation}
where $a_0 = c_1^2/c_3^2 > 0$, $a_1 = 4c_1(-2 c_1 + c_2 + 2 c_3)/(3c_3^2) > 0$, $a_2 = 2(2 c_1^2 - 2 c_1 c_2 -
7 c_1 c_3 - 2 c_2 c_3 + 2 c_3^2)/(3c_3^2)$, and $a_3 = 4(2 c_1 + c_2 - 2 c_3)/(3 c_3)$.

Denote $h(y_1) = a_0 + a_1 y_1 + a_2 y_1^2 + a_3 y_1^3 + y_1^4$. Note that $h(\infty) = \infty$, $h(-c_1/c_3) = - c_1^2 (c_1^2 + 8 c_1 c_2 - 2 c_1 c_3 + 8 c_2 c_3 + c_3^2)/(3 c_3^4) < 0$, $h(0) = c_1^2/c_3^2 > 0$, $h(1) = -(c_1 - c_3)^2/(3 c_3^2) < 0$, and $h(\infty) = \infty$. Then $h(y_1)=0$ yields four real roots in $(\infty, -c_1/c_3)$, $(-c_1/c_3, 0)$, $(0, 1)$, and $(1, \infty)$, respectively. That is, there is one and only one $y_1 \in (1, \infty)$.

According to \citeauthor{tong2014} (2014, equation~(12)), \begin{equation}\label{y1solution}
y_1=-\frac{a_3}{4}+\frac{\sqrt{A_1}}{2} + \frac{\sqrt{C_1}}{2}\ ,
\end{equation}
where
\begin{eqnarray*}
	A_1 &=& -\frac{2  a_2}{3}+\frac{ a_3^2}{4}+\frac{G_1}{3\times 2^{1/3}}\ ,\\
	C_1 &=& -\frac{4 a_2}{3}+\frac{a_3^2}{2}-\frac{G_1}{3\times 2^{1/3}} + \frac{-8 a_1+4 a_2 a_3-a_3^3}{4 \sqrt{A_1}}\ ,\\
	G_1 &=& \left(F_1-\sqrt{F_1^2-4 E_1^3}\right)^{1/3}+ \left(F_1+\sqrt{F_1^2-4 E_1^3}\right)^{1/3}\ ,\\
	E_1 &=& 12 a_0+a_2^2-3 a_1 a_3\ ,\\
	F_1 &=& 27 a_1^2-72 a_0 a_2+2 a_2^3-9 a_1 a_2 a_3+27 a_0 a_3^2\ .
\end{eqnarray*}
The calculation of $G_1$, $A_1$, $C_1$, and $y_1$ are operations among complex numbers, while $y_1$ at the end would be a real number.

The procedure of obtaining analytic solutions of $w_1, w_2, w_3$ would be, (1) obtain $y_1$ from \eqref{y1solution}; (2) obtain $y_2$ from \eqref{eq:y2solution}; (3) $w_1 = y_1/(y_1 + y_2 + 1)$, $w_2 = y_2/(y_1 + y_2 + 1)$, $w_3 = 1/(y_1 + y_2 + 1)$.

Now we discuss some special cases.
\begin{itemize}
	\item[(i)] If $c_1 = c_2 < c_3$, then $w_1 = w_2$ and thus $y_1 = y_2$. Both \eqref{eq:f1=f3y} and \eqref{eq:f2=f3y} yield
	$y_1 = c_3^{-1}(-2 c_1 + c_3 + \sqrt{4 c_1^2 - c_1 c_3 + c_3^2})$, which implies
	$$w_1 = w_2 = \frac{-2 c_1 + c_3 + \Delta_1}{-4 c_1 + 3 c_3 + 2 \Delta_1}, \quad w_3 = \frac{c_3}{-4 c_1 + 3 c_3 + 2 \Delta_1}$$
	where $\Delta_1 = \sqrt{4 c_1^2 - c_1 c_3 + c_3^2}$. Note that $w_1 > w_3$ since $\Delta_1 > 2 c_1$.
	\item[(ii)] If $c_1 < c_2 = c_3$, then $w_2 = w_3$ and thus $y_2 = 1$. From \eqref{eq:f1=f3y} we get
	$y_1 = 3 c_3^{-1}(- c_1 + 2 c_3 + \sqrt{c_1^2 - c_1 c_3 + 4 c_3^2})$, which implies
	$$w_1 = \frac{- c_1 + 2 c_3 + \Delta_2}{- c_1 + 8 c_3 +  \Delta_2}, \quad w_2 = w_3 = \frac{3 c_3}{- c_1 + 8 c_3 +  \Delta_2}$$
	where $\Delta_2 = \sqrt{c_1^2 - c_1 c_3 + 4 c_3^2}$. Note that $w_1 > w_2$ since $\Delta_2 > c_1 + c_3$.
	\item[(iii)] If $c_1 < c_2 < c_3$, then $y_1, y_2$ and thus $w_1, w_2, w_3$ can be obtained analytically. We have proven $y_1 \geq y_2 \geq 1$. Using \eqref{eq:f1=f3y} and \eqref{eq:f2=f3y}, it can be verified that $y_1 \neq y_2$ unless $c_1=c_2$; and $y_2 \neq 1$ unless $c_2 = c_3$. That is, $y_1 > y_2 > 1$ and $w_1 > w_2 > w_3$.
\end{itemize}
\hfill{$\Box$}

\medskip\noindent
{\bf Proof of Theorem~\ref{thm:ustpi}:}
(i) is straightforward. (ii) follows from the facts in the proof of Lemma~\ref{lemma:pii}, ${\mathbf c}_{iJ} \equiv {\boldsymbol\pi}_i$; ${\mathbf 1}^T {\mathbf c}_{ij} = 0, j=1, \ldots, J-1$; and ${\mathbf 1}^T {\mathbf c}_{iJ} = 1$. (iii) and (iv) can be verified using the formulae of ${\mathbf c}_{ij}$ in Section~\ref{subsec:cij}.
\hfill{$\Box$}

\end{document}